\documentclass[11pt]{amsart}
\usepackage{amsmath,amssymb,amsthm}
\usepackage[latin1]{inputenc}
\usepackage{version,tabularx,multicol}
\usepackage{graphicx,float}
\usepackage[usenames]{color}

\newcommand{\reff}[1]{(\ref{#1})}

\theoremstyle{plain}
\newtheorem{theo}{Theorem}[section]
\newtheorem*{theo*}{Theorem}

\newtheorem{cor}[theo]{Corollary}
\newtheorem{prop}[theo]{Proposition}

\newtheorem{lem}[theo]{Lemma}
\newtheorem{defi}[theo]{Definition}
\theoremstyle{remark}
\newtheorem{rem}[theo]{Remark}

\newcommand{\Hm}{H_{\text{max}}}
\newcommand{\Tm}{{T_{\text{max}}}}
\newcommand{\rP}{{\rm P}}
\newcommand{\rL}{{\rm L}}
\newcommand{\mL}{{\mathcal L}}
\newcommand{\rE}{{\rm E}}
\newcommand{\cm}{{\mathcal{M}}}
\newcommand{\cx}{{\mathcal{X}}}

\def\Bc{\ensuremath{\mathcal{B}}}

\def\C{\ensuremath{\mathcal{C}}}

\def\D{\mathcal{D}}
\def\Dc{\mathcal{D}}
\def\Ec{\ensuremath{\mathcal{E}}}

\def\E{\ensuremath{\mathbb{E}}}

\def\F{\ensuremath{\mathcal{F}}}

\def\G{\ensuremath{\mathcal{G}}}

\def\I{\ensuremath{\mathcal{I}}}
\def\J{\ensuremath{\mathcal{J}}}

\def\M{\ensuremath{\mathcal{M}}}

\def\N{\ensuremath{\mathbb{N}}}
\def\NN{\ensuremath{\mathbf{N}}}
\def\cn{\ensuremath{\mathcal{N}}}

\def\P{\ensuremath{\mathbb{P}}}
\def\PP{\mathbf{P}}

\def\qq{q}

\def\R{\ensuremath{\mathbb{R}}}

\def\W{\ensuremath{\mathcal{W}}}

\def\to{\rightarrow}

\newcommand{\ind}{{\bf 1}}

\newcommand{\norm}[1]{\mathop{\parallel\! #1 \! \parallel}\nolimits}
\newcommand{\val}[1]{\mathop{\left| #1 \right|}\nolimits}
\newcommand{\inv}[1]{\mathop{\frac{1}{ #1}}\nolimits}
\newcommand{\expp}[1]{\mathop {\mathrm{e}^{ #1}}}

\begin{document}

\title{A Williams' decomposition for spatially dependent superprocesses}

\date{\today}

\author{Jean-Fran\c cois Delmas}

\address{
Jean-Fran\c cois Delmas,
Universit\' e Paris-Est,  CERMICS, 6-8
av. Blaise Pascal, 
  Champs-sur-Marne, 77455 Marne La Valle, France.}

\email{delmas@cermics.enpc.fr}

\author{Olivier H\' enard} 

\address{
Olivier H\' enard,
Universit\' e Paris-Est,  CERMICS, 6-8
av. Blaise Pascal, 
  Champs-sur-Marne, 77455 Marne La Valle, France.}

\email{henardo@cermics.enpc.fr}

\thanks{This work is partially supported by the ``Agence Nationale de
  la Recherche'', ANR-08-BLAN-0190.}

\begin{abstract}
  We  present  a genealogy  for  superprocesses  with a  non-homogeneous
  quadratic branching  mechanism, relying on  a weighted version  of the
  superprocess and a Girsanov  theorem. We then decompose this genealogy
  with respect  to the last individual  alive (William's decomposition).
  Letting the extinction time tend  to infinity, we get the Q-process by
  looking at the superprocess from  the root, and define another process
  by  looking from  the top.   Examples including  the  multitype Feller
  diffusion and the superdiffusion are provided.
\end{abstract}



\maketitle

\section{Introduction}

% Superprocesses with non homogeneous branching mechanisms are known for a long time. The study of their fine properties in the case of the superdiffusions is more recent, see Engl\"ander and Pinsky \cite{EP99}.

Even if superprocesses with very general branching mechanisms are known,
most of the works devoted to  the study of their genealogy are concerned
with  homogeneous  branching   mechanisms,  that  is,  populations  with
identical individuals.  Four distinct approaches have  been proposed for
describing  these  genealogies.   When   there  is  no  spatial  motion,
superprocesses  are  reduced to  continuous  state branching  processes,
whose  genealogy can  be  understood  by a  flow  of subordinators,  see
Bertoin  and Le  Gall \cite{BLG00},  or by  growing discrete  trees, see
Duquesne and Winkel \cite{DW07}.  With a spatial motion, the description
of the genealogy can be done  using the lookdown process of Donnelly and
Kurtz \cite{DK99}  or the  snake process of  Le Gall  \cite{LG93}.  Some
works  generalize   both  constructions  to   non-homogeneous  branching
mechanisms:  Kurtz  and   Rodriguez  \cite{KR11}  recently  extended  the
lookdown process  in this direction whereas Dhersin  and Serlet proposed
in \cite{DS00} modifications of the snake.

Using the genealogy,  it is natural to consider the corresponding
Williams' decomposition,  which  is  named  after the  work  of  Williams
\cite{WI74}  on  the  Brownian  excursion. After  Aldous  recognized  in
\cite{A91} the genealogy of a  branching process in this excursion, they
also  designate decompositions  of branching  processes with  respect to
their  height,  see  Serlet  \cite{SE96}  for  the  quadratic  branching
mechanism  or  Abraham  and  Delmas \cite{AD08}  for  general  branching
mechanism.   Their  interest is  twice:  they  allow  to understand  the
behavior of processes  at the top, see Goldschmidt  and Haas \cite{GH10}
for  an application  of this  approach, and  to investigate  the process
conditioned  on  non  extinction,  or  Q-process,  see  \cite{SE96}  and
Overbeck \cite{OV94}.

For Markov processes with absorbing  states, the Q-process is defined as
the process  conditioned on non  absorption in remote time,  see Darroch
and  Seneta \cite{DS65}.  Lamperti and  Ney \cite{LN68}  found  a simple
construction  in the  case of  discrete branching  processes.  Later on,
Roelly and  Rouault \cite{RR89} provided a superprocess  version of this
result.  Q-processes have intrinsic  interest as  a model  of stochastic
population, see Chen and  Delmas \cite{CD10}. They also find application
in the study of the associated martingale, see Lyons, Pemantle and Peres
\cite{LPP95} in a discrete setting. Understanding this martingale allows
to better understand the original process, see Engl\"ander and Kyprianou
\cite{EK04} for superprocesses with non homogeneous branching mechanism.

Our primary interest  is to present a genealogy  for superprocess with a
non-homogeneous  quadratic  branching mechanism,  to  condition it  with
respect  to its  height (this  is the  William's decomposition),  and to
study the associated Q-process.

\medskip

Let $X=(X_t, t  \geq 0)$ be an $(\mL,\beta,\alpha)$  superprocess over a
Polish space $E$. The underlying  spatial motion $Y=(Y_t, t\geq 0)$ is a
Markov process  with infinitesimal generator $\mL$ started  at $x$ under
$\rP_{x}$.  The non-homogeneous quadratic branching mechanism is denoted
by  $\psi(x,\lambda)=\beta(x)   \lambda  +  \alpha(x)   \lambda^2$,  for
suitable  functions $\beta$  and  $\alpha$ (explicit  conditions can  be
found in Section \ref{section0}).  Let $\P_{\nu}$ be the distribution of
$X$ started  from the finite measure  $\nu$ on $E$, and  $\N_{x}$ be the
corresponding  canonical measure  of  $X$ with  initial  state $x$.   In
particular, the process $X$  under $\P_{\nu}$ is distributed as $\sum_{i
  \in \I}  X^i$, where $\sum_{i  \in \I} \delta_{X^i}(dX)$ is  a Poisson
Point  measure with  intensity $\int_{x\in  E} \nu(dx)  \N_{x}(dX)$.  We
define  the extinction time  of $X$:  $\Hm=\inf\{t>0, X_t=0\}$,  and
assume 
that $X$  suffers almost sure extinction, that is  $\N_{x} \left[ \Hm=\infty
\right]=0$   for  all   $x\in   E$.   Using an $h$-transform from
Engl\"ander     and    Pinsky     \cite{EP99} and a Girsanov
transformation from Perkins \cite{PE99}, we provide a genealogical
structure for the superprocess $X$, see Proposition
\ref{proplinkgenealogy}, by transferring the genealogical
structure of an homogeneous superprocess.  

We    define  the   function
$v_h(x)=\N_{x}\left[X_h \neq 0 \right]=\N_{x}\left[\Hm\geq h\right]$ and
a family of probability measures by setting:
\[ \forall \ 0 \leq t <h, \ \ 
  \frac{d\rP^{(h)}_{x \ |\D_t}} {d\rP_{x  \ |\D_t}}= \frac{\partial_{h} v_{h-t}(Y_t)}{\partial_{h} v_{h}(x)}
\expp{- \int_{0}^{t}  ds \ \partial_{\lambda} \psi(Y_s,v_{h-s}(Y_s)) },
\]
where $\D_t=\sigma(Y_s,  0 \leq s\leq  t)$ is the natural  filtration of
$Y$, see Lemma  \ref{AC2}.  Using the genealogical structure  of $X$, we
give  a  decomposition  of  the  superprocess $X$  with  respect  to  an
individual  chosen at  random, also  called a  Bismut  decomposition, in
Proposition  \ref{propbismut3}.  The  following  Theorem, see  Corollary
\ref{cor:Nh} for a precise statement and Corollary \ref{corwilliams} for
a  statement  under  $\P_{\nu}$,   gives  a    Williams'
decomposition of $X$, that is a spine decomposition  with respect to its
extinction time $\Hm$.
\begin{theo*}{(Williams' decomposition under $\N_{x}$)} 
Assume that the $(\mL,\beta,\alpha)$ superdiffusion $X$ suffers almost
sure extinction and some regularities on $\alpha$ and $\beta$. 
\begin{itemize}
 \item[(i)] The distribution of $\Hm$ under $\N_x$ is characterized
   by: $\N_{x}[\Hm > h]= v_h(x)$. 
 \item[(ii)] Conditionally on $\{\Hm=h_0\}$, the $(\mL,\beta,\alpha)$ superdiffusion $X$ under 
$\N_x$ is distributed as $X^{(h_0)}$ constructed as follows. Let  $x\in  E$  and
$Y_{[0,h_0)}$  be distributed  according  to 
  $\rP^{(h_0)}_{x}$.  Consider  the Poisson point  measure $\cn=\sum_{j\in
    J}  \delta_{(s_j,X^{j})}$ on  $[0,h_0) \times \Omega$  with intensity:  
\[
2\ind_{[0,  h_0)}(s)   ds  \;  \ind_{   \{  \Hm(X)  <   h_0-s  \}}
  \alpha(Y_s)\; \N_{Y_s}[ dX].
\]
 The process $X^{(h_0)}=(X^{(h_0)}_t, t\geq
  0)$ is then defined  for all $t\geq 0$ by:
\[
X^{(h_0)}_t= \sum_{j\in J, \, s_j<t} X_{t-s_j }^{j}.
\]
\end{itemize}
\end{theo*}
The proof  of this  Theorem relies on  a William's decomposition  of the
genealogy  of  $X$,  see  Theorem \ref{theowilliams}.   Notice  it  also
implies the existence  of a measurable family $(\N_{x}^{(h)},  h >0)$ of
probabilities such that $\N_{x}^{(h)}$  is the distribution of $X$ under
$\N_{x}$ conditionally on $\{ \Hm=h \}$.

\medskip
We shall from now on consider the case of $Y$ a diffusion on $\R^{K}$ or a pure jump process on a finite state space.
The generalized eigenvalue $\lambda_{0}$ of the operator $\beta - \mL$ is defined in Pinsky \cite{PI95} for diffusion on $\R^{d}$. For finite state space, it reduces to the Perron Frobenius eigenvalue, see Seneta \cite{SE06}. In both cases, we have:
\begin{align*}
\lambda_{0} = \sup{\{ \ell \in \R, \exists u\in \D(\mL), u > 0 \mbox{ such that } (\beta - \mL) u = \ell \  u \}} \cdot
\end{align*}
We   assume  that  the   space  of   positive  harmonic   functions  for
$(\beta-\lambda_{0}) - \mL$ is  one dimensional, generated by a function
$\phi_{0}$.  From these  assumptions, we    have that  the space  of
positive harmonic functions of the adjoint of $(\beta-\lambda_{0})- \mL$
is one dimensional,  and we denote by $\tilde{\phi_{0}}$  a generator of
this   space. We also assume   that $\phi_{0}$ is bounded from below
  and   above   by   positive    constants   and   that   the   operator
  $(\beta-\lambda_{0})-   \mL$  is   product  critical, that is  $\int_{E} dx \:
\phi_{0}(x)  \:  \tilde{\phi_{0}}(x)  <   \infty$. Thanks to the
product-critical property, the 
probability measure $\rP^{\phi_0}$, given by:
\[
\forall t \geq 0, \quad \frac{d \rP_{x \ | \D_t}^{\phi_{0}}} {d \rP_{x \ |
    \D_t}} = \frac{\phi_{0}(Y_t)}{\phi_{0}(Y_0)} \expp{- \int_{0}^{t} ds
  \ (\beta(Y_s)- \lambda_{0})},
\]  defines  a   recurrent  Markov  process  (in  the   sense  given  by
\reff{H6''}).  Since $\phi_{0}$ is bounded  from below and from above by
two positive  constants, the  non-negativity of $\lambda_0$  implies the
weak  convergence  of  the  spine   that  is  the  weak  convergence  of
$\rP_x^{(h)}$   towards   $\rP_x^{(\infty   )}$   which  is   given   by
$\rP_x^{\phi_0}$,  see Proposition  \ref{prop:H4+Pf=Pinf}.   An explicit
expression  is given  for $\rP_x^{(h)}$  and $\rP_x^{\phi_0}$  in Lemmas
\ref{lem:spinemulti}  and  \ref{lem:spinediff}.   The non-negativity  of
$\lambda_0$  implies  the  almost  sure  extinction of  $X$,  see  Lemma
\ref{nonnegativelambdaH1}.   Under  very  general conditions,  the  weak
convergence  of the spine  implies the  convergence of  the superprocess
(Corollary     \ref{cor:cvNh})     and     its    genealogy     (Theorem
\ref{theo-convN}). We can easily state them in the particular case of an
underlying motion being  a diffusion or a pure jump  process on a finite
state space.  We also see that $\N_{x}^{(\infty)}$, defined below,
is actually the  law of the Q-process, defined as the  weak limit of the
probability measures $\N_{x}^{( \geq  h)} = \N_{x} \left[\; \cdot \; |
  \Hm \geq h\right]$, see also Lemma \ref{lemmaconv}.

\begin{theo*}[Q-process under $\N_{x}$]
  Assume that $\lambda_0 \geq 0$.    Let   $Y$  be
  distributed  according to  $\rP^{\phi_0}_{x}$,  and, conditionally  on
  $Y$, let $\cn=\sum_{j\in \I}  \delta_{(s_j,X^{j})}$ be a Poisson point
  measure with intensity:
\[
2  \ind_{\R^+}(s)   ds   \;  \alpha(Y_s)
   \N_{Y_s}[dX]. 
\]
Consider the process $X^{(\infty )}=(X^{(\infty )}_t, t\geq
  0)$, which is  defined  for all $t\geq 0$ by:
\[
X^{(\infty )}_t= \sum_{j\in J, \, s_j<t} X_{t-s_j }^{j}, 
\]
and denote by $\N_{x}^{(\infty)}$ its distribution. Then, for all $t\geq
0$, the distribution of $(X_s, s\in[0,t])$ under $\N_x^{(h)}$ or
$\N_x^{(\geq h)}$ converges weakly to
  $(X_s^{(\infty)}, s\in[0,t])$.
\end{theo*}
Notice  those  results  also   hold  under  $\P_{\nu}$  (see  Corollary
\ref{cor:cvNh}). It is  interesting to notice that the  law of the spine
$\rP^{\phi_0}$ is  quite different  from that of  the backbone  given in
\cite{EP99},    see     also    Remark
\ref{backbone-spine}.

\begin{rem}
  As  noticed by  Li  \cite{LI92},  the  multitype Dawson  Watanabe
  superprocess   can   be  understood   as   a  single   non-homogeneous
  superprocess  on an  extended space.  The above  Theorem  on Q-process
  provides  a construction of  the Q-process  associated to  a multitype
  Dawson  Watanabe  superprocess  considered  in Champagnat  and  Roelly
  \cite{CH08}, and  this construction gives  a precise meaning  to ``the
  interactive immigration'' introduced in Remark 2.8 of \cite{CH08}.
\end{rem}

\begin{rem}
  In  \cite{EK04}, a spinal decomposition of
  the semi-group of a  Doob $h$-transform  of the  $(\mL,\beta,\alpha)$ superdiffusion is  provided, see Theorem 5  of  \cite{EK04}. The second item of Corollary
  \ref{cor:cvNh},  together with  Lemma \ref{lemmaconv}  gives a pathwise decomposition of the genealogy and establishes  that the
  process considered actually is the Q-process.
\end{rem}
We   also   prove  the weak   convergence   of   the  probability   measures
$(\N_{x}^{(h)},   h>0)$   backward  from   the   extinction  time.   Let
$\rP^{(-h)}$ denote the push forward probability measure of $\rP^{(h)}$,
defined by:
\begin{align*}
\label{defP-h}
\rP^{(-h)}((Y_{s}, s \in [-h,0])  \in  \bullet)  =
\rP^{(h)} ((Y_{h+s}, s \in [-h,0]) \in \bullet).
\end{align*}
The product criticality assumption yields the existence of a probability
measure $\rP^{(-\infty)}$  such that for all  $x \in E$,  $t\geq 0$, and
$f$ bounded measurable:
\begin{align*}
  \rE_{x}^{(-h)}    \big[f(Y_{s}, s \in[-t,0])\big]   \;    \xrightarrow[h   \to
  +\infty]{} \; \rE^{(-\infty)} \big[f(Y_{s}, s \in [-t,0])\big].
\end{align*}
Once again, the convergence of  the spine implies the convergence of the
superprocess.  The  following result corresponds  to the second  item of
Theorem \ref{convtop}.
\begin{theo*}[Asymptotic distribution at the extinction time]
Assume that $\lambda_0 > 0$. Then the   process $(X_{h+s},  s\in [-t,0])$
under $\N_{x}^{(h)}$ weakly  converges towards  $X^{(-\infty
  )}_{[-t,0]}$, where for $s\leq 0$: 
\[
X^{(-\infty )}_s=\sum_{j\in J,  \; s_j<s} X^j_{s-s_j},
\]
and conditionally on $Y$ with distribution $\rP^{(-\infty)}$, $\sum_{j\in J}
\delta_{(s_j,X^{j})}$ is a Poisson point measure with
intensity:
 \[ 
2 \;  \ind_{\{s<0\}} \alpha(Y_s) \: ds \; \ind_{\{\Hm(X) < -s\}} \;
\N_{Y_s}[ dX].
\] 
\end{theo*}

\begin{rem}
  Considering a  superprocess with homogeneous  branching mechanism, the
  Q-process may be easily defined  from the well known Q-process for the
  total  mass process (see  for instance  \cite{CD10} in  the case  of a
  general branching  mechanism).  Thus the  recurrence condition imposed
  on the spatial motion is not necessary for Williams decomposition, but
  it seems more  natural in order to get  the asymptotic distribution at
  the extinction time.
\end{rem}

\begin{rem}
   \label{rem:penalisation}
The genealogy  of $X$ defined in Proposition  \ref{proplinkgenealogy} allows us
to interpret the following  probability measure $\rP^{(B,t)}_{x}$ as the
law  of the  ancestral lineage  of an  individual sampled  at  random at
height    $t$    (see     the    Bismut    decomposition,    Proposition
\ref{propbismut3}):
\begin{align*}
\frac{d \rP_{x \ | \D_{t}}^{(B,t)}}{d \rP_{x \   | \D_{t}}} =
\frac{\expp{-\int_{0}^{t} ds \;
    \beta(Y_s)}}{\rE_{x}\left[\expp{-\int_{0}^{t} ds \;
      \beta(Y_s)}\right]}\cdot
\end{align*}
We prove in  Lemma \ref{lem:convbismutspine} that, if $\phi_{0}$ is
bounded from below and above by positive constants and that the operator
$(\beta-\lambda_{0})-  \mL$  is  product  critical, then  the  ancestral
lineage of  an individual sampled at  random at height  $t$ under $\N_x$
converges  as  $t  \to  \infty$  to  the  law  of  the  spine,  that  is
$\rP_{x}^{(B,t)}$   converges  weakly   to   $\rP_{x}^{\phi_0}$.  This
Feynman-Kac type penalization result (see  Chapter 2 of Roynette and Yor
\cite{RY09})  heavily relies  on the  product criticality  assumption, but
holds  without  restriction  on  the  sign of  $\lambda_0$.  It  may  be
interpreted  as an example  of the  so called  globular state  in random
polymers, investigated in Cranston, Koralov and Molchanov \cite{CK10}.
\end{rem}

\textbf{Outline}.   We give  some  background on  superprocesses with  a
non-homogeneous branching mechanism in  the Section 2.  Section 3 begins
with the definition of the $h$-transform in the sense of Engl\"ander and
Pinsky, Definition  \ref{htransform}, goes  on with a  Girsanov Theorem,
Proposition  \ref{girsanov}, and  ends  up with  the  definition of  the
genealogy, Proposition \ref{proplinkgenealogy}, by combining both tools.
Section 4 is mainly devoted to the proof of the William's decomposition,
Theorem  \ref{theowilliams}. By the  way, we  give a  decomposition with
respect  to  a  randomly  chosen  individual, also  known  as  a  Bismut
decomposition, in  Proposition \ref{propbismut2}.  Section  5 gives some
applications  of the Williams'  decomposition: We  first prove  in Lemma
\ref{lemmaconv}  that the  limit  of the  superprocesses conditioned  to
extinct  \textit{at} a  remote  time coincide  with  the Q-process  (the
superprocess conditioned  to extinct  \textit{after} a remote  time) and
actually show in  Theorem \ref{theo-convN} that such a  limit exists. We
also consider  in Theorem \ref{convtop}  the convergence of  the process
seen from the top (so, backward from the extinction time).  All previous
results are provided with a set  of assumptions. We then give in Section
\ref{sec:H} sufficient  conditions for these assumptions to  be valid in
term of the generalized eigenvector and eigenvalue, then check they hold
in Section \ref{sec:twoexamples} in two examples: the finite state space
superprocess (with mass process  the multitype Feller diffusion) and the
superdiffusion.

% All results are stated at the level of genealogy and at the level of measure valued processes, both under the excursion measure and under the usual probability measure.

\section{Notations and definitions}
\label{section0}

This section, based on the lecture notes of Perkins \cite{PE99}, provides us with basic material about superprocesses, relying on their characterization via the Log Laplace equation.

We first introduce some definitions:

\begin{itemize}
\item $(E, \delta)$ is a Polish space, $\Bc$ its Borel sigma-field.
\item $\Ec$ is the set of real valued measurable functions and $b\Ec \subset \Ec$ the subset of bounded functions.
\item $\C(E,\R)$, or simply $\C$, is the set of continuous real valued functions on $E$, $\C_b \subset \C$ the subset of continuous bounded functions.
\item $D(\R^{+},E)$, or simply $D$, is the set of c\`adl\`ag  paths of
  $E$ equipped with the Skorokhod topology, $\Dc$ is the Borel sigma
  field on $D$,  and $\Dc_{t}$ the canonical right continuous filtration
  on $D$. 
\item For each set of functions, the superscript $.^{+}$ will denote the subset of the non-negative functions: For instance, $b\Ec^{+}$ stands for the subset of non negative functions of $b\Ec$.
\item $\M_{f}(E)$ is  the space of finite measures  on $E$. The standard
  inner product notation  will be used: for $g  \in \Ec$ integrable with
  respect to $M \in \M_f(E)$, $M(g)= \int_{E} M(dx) g(x)$.
\end{itemize}

We can now introduce the two main ingredients which enter in the definition of a superprocess, the spatial motion and the branching mechanism:
\begin{itemize}
\item  Assume  $Y=(D,\Dc,\Dc_{t},Y_t, \rP_{x})$ is  a  Borel strong  Markov
  process. ``Borel'' means that $x \to \rP_{x}(A)$ is $\Bc$ measurable for
  all  $A \in  \Bc$. Let  $\rE_{x}$ denote  the expectation  operator, and
  $(P_{t},     t    \geq     0)$    the     semi-group     defined    by:
  $P_t(f)(x)=\rE_{x}[f(Y_t)]$. We  impose the additional  assumption that
  $P_{t}:  \C_{b} \to  \C_{b}$.  In  particular the  process $Y$  has no
  fixed discontinuities. The generator  associated to the semi-group will
  be denoted $\mL$. Remember $f$ belongs to the domain $\Dc(\mL)$ of
  $\mL$ if $f \in \C_{b}$ and for some $g \in \C_{b}$,
\begin{equation}
\label{generator}
f(Y_t)-f(x)-\int_{0}^{t} ds \ g(Y_s) \mbox{ is a } P_{x} \mbox{ martingale for all } x  \mbox{ in } E, 
\end{equation}
in which case $g = \mL(f)$.
\item The functions $\alpha$ and $\beta$ being elements of $\C_{b}$,
  with $\alpha$ bounded from below by a positive constant, the non-homogeneous
  quadratic branching mechanism $\psi^{\beta,\alpha}$ is
  defined by: 
\begin{equation}
\label{eq:branchingmechanism}
\psi^{\beta,\alpha}(x,\lambda) =  \beta(x) \lambda + \alpha(x)
\lambda^2, 
\end{equation}
for all $x  \in E$ and $\lambda  \in \R$. We will just  write $\psi$ for
$\psi^{\beta,\alpha}$ when there is  no possible confusion.  If $\alpha$
and $\beta$ are constant functions, we will call the branching mechanism
(and by extension, the corresponding superprocess) homogeneous.
\end{itemize}

The mild form of the Log Laplace equation is given by the integral
equation, for $\phi, f \in b\Ec^{+}$, $t\geq 0$, $x \in E$:

\begin{equation}
\label{mildequation}
\quad u_t(x) + \rE_{x} \left[ \int_{0}^{t} ds \ \psi(Y_s,u_{t-s}(Y_s))\right]= \rE_{x} \left[f(Y_t)+\int_{0}^{t} ds \ \phi(Y_{s}) \right] \cdot
\end{equation}

\begin{theo}({\cite{PE99}, Theorem II.5.11}) 
\label{theo1}
Let $\phi, f \in b\Ec^{+}$.
There is a unique jointly (in $t$ and $x$) Borel measurable solution $u_t^{f,\phi}(x)$ of equation (\ref{mildequation}) such that $u_t^{f,\phi}$ is bounded on $[0,T]\times E$ for all $T>0$. Moreover, $u_{t}^{f,\phi} \geq 0$ for all $t \geq 0$.
\end{theo}
We shall  write $u^{f}$ for $u^{f,0}$ when $\phi$ is null.

We introduce the canonical space of continuous applications from
$[0,\infty)$ to $\M_{f}(E)$, denoted by $\Omega \mathrel{\mathop:}=
\C(\R^{+}, \M_{f}(E))$, endowed with its Borel sigma field $\F$, and the
canonical right continuous filtration $\F_t$. Notice that $\F=\F_\infty
$. 

\begin{theo}({\cite{PE99}, Theorem II.5.11}) 
\label{theo1a}
Let  $u_t^{f,\phi}(x)$  denote   the  unique  jointly  Borel  measurable
solution  of equation (\ref{mildequation})  such that  $u_t^{f,\phi}$ is
bounded on $[0,T]\times E$ for  all $T>0$.  There exists a unique Markov
process $X=(\Omega, \F, \F_{t}, X_t, (\P^{(\mL,\beta,\alpha)}_\nu, \nu \in
\M_{f}(E)))$ such that:
\begin{equation}
\label{laplaceequation}
\forall \phi, f \in b\Ec^{+}, \  \qquad \E_{\nu}^{(\mL,\beta,\alpha)} \left[
\expp{-X_t(f)-\int_{0}^{t} ds\ X_s(\phi)} \right]=\expp{-\nu(u_{t}^{f,\phi})}. 
\end{equation}
$X$ is called the $(\mL,\beta,\alpha)$-superprocess.
\end{theo}
We now state the existence theorem of the canonical measures:
\begin{theo}({\cite{PE99}, Theorem II.7.3})
\label{canonical}
There   exists   a  measurable   family   of  $\sigma$-finite   measures
$(\N_{x}^{(\mL,\beta,\alpha)},x\in E)$  on $(\Omega, \F)$  which satisfies
the following properties: If $\sum_{j \in \J} \delta_{(x^j,X^j)}$ is a
Poisson point  measure on $E  \times \Omega$ with intensity  $\nu(dx) \;
\N^{(\mL,\beta,\alpha)}_{x}$,   then  $\sum_{j  \in   \J}  X^{j}$   is  an
$(\mL,\beta,\alpha)$-superprocess started at $\nu$.
\end{theo}
We will often abuse notation by denoting $\P_{\nu}$ (resp. $\N_x$)
instead of $\P^{(\mL,\beta,\alpha)}_{\nu}$
(resp. $\N^{(\mL,\beta,\alpha)}_{x}$), and $\P_{x}$ instead of
$\P_{\delta_{x}}$ when starting from $\delta_x$ the Dirac mass at point $x$. 

Let $X$  be a $(\mL,\beta,\alpha)$-superprocess.   The exponential formula
for Poisson point measures yields the following equality:
\begin{equation}
\label{eq:masterformula}
\forall f \in b\Ec^{+} , \ \N_{x_0}\big[1- \expp{-X_{t}(f)}\big]= -\log
\E_{x_0} \big[ \expp{-X_{t}(f)}\big] = u_t^{f}(x_0),  
\end{equation}
where $u_t^{f}$ is (uniquely) defined by equation (\ref{laplaceequation}).

Denote  $\Hm$ the extinction time of $X$:
\begin{equation}
   \label{eq:Hmax}
\Hm = \inf \{t >0; \; X_t =0 \}.
\end{equation}

\begin{defi}[Global extinction]
  The superprocess $X$ suffers global extinction if $\P_\nu(\Hm
  <\infty )=1$  for all $\nu \in \M_{f}(E)$.
\end{defi}

We will need the the following assumption:

\medskip
 \noindent
$(H1)$ \textbf{The $(\mL,\beta,\alpha)$-superprocess satisfies the global
  extinction property.}
\medskip

We shall be interested in the function
\begin{equation}
   \label{eq:def-vtx}
v_t(x) = \N_{x}[\Hm > t].
\end{equation}
We set $v_{\infty}(x)= \lim_{t \to \infty} \downarrow v_{t}(x)$. The
global extinction property is easily stated using $v_\infty $. 
\begin{lem}
\label{lemmaextinction}
The global extinction property holds if and only if $v_{\infty}=0$. 
\end{lem}
See also Lemma \ref{EDPV} for other properties of the function $v$. 

\begin{proof}
The exponential formula for Poisson point measures yields:
\[
\P_{\nu}(\Hm \leq t) = \expp{-\nu(v_t)}.
\]
To conclude, let $t$ goes to infinity in the previous equality to get:
\[
\P_{\nu}(\Hm < \infty) = \expp{-\nu(v_{\infty})}.
\]
\end{proof}

For  homogeneous  superprocesses ($\alpha$  and  $\beta$ constant),  the
function $v$ is  easy to compute and the global  extinction holds if and
only  $\beta$   is  non-negative.  Then,   using  stochastic  domination
argument, one  get that a  $(\mL,\beta,\alpha)$-superprocess, with $\beta$
non-negative,  exhibits  global  extinction  (see \cite{EK04}  p.80  for
details).

\section{A genealogy for the spatially dependent superprocess}
\label{sectiongenealogy}

We first  recall (Section \ref{subsectionhtransform})  the $h$-transform
for   superprocess   introduced  in   \cite{EP99}   and  then   (Section
\ref{subsectiongirsanov})  a Girsanov  theorem previously  introduced in
\cite{PE99} for  interactive superprocesses.  Those  two transformations
allow us  to give  a Radon-Nikodym derivative  of the distribution  of a
superprocess  with non-homogeneous branching  mechanism with  respect to
the  distribution  of  a  superprocess  with  an  homogeneous  branching
mechanism.   The  genealogy  of  the superprocess  with  an  homogeneous
branching  mechanism  can  be  described  using a  Brownian  snake,  see
\cite{DLG02}.   Then, in Section  \ref{subsectionTransport}, we  use the
Radon-Nikodym derivative to transport  this genealogy and get a genealogy
for the superprocess with non-homogeneous branching mechanism.

\subsection{$h$-transform for superprocesses}
\label{subsectionhtransform}

We first introduce a new probability measure on $(D, \D)$ using the next Lemma. 

\begin{lem}
\label{AClemma}
Let  $g$  be  a  positive   function  of  $\D(\mL)$  such that $g$ is bounded from below by a positive constant. Then, the process  $\big(\frac{g(Y_t)}{g(x)} \expp{-\int_{0}^{t}
  ds \ (\mL g/g)(Y_s)}, t \geq 0\big)$ is a positive martingale under $\rP_x$.
\end{lem}  

We set $\D_g(\mL)=\{ v \in \C_{b}, \; gv \in \D(\mL) \}$. 

\begin{proof}
Let $g$ be as in Lemma \ref{AClemma} and $f \in \D_g(\mL)$. The process:
\begin{equation*}
\bigg((fg)(Y_t)-(fg)(x)-\int_{0}^{t} ds \ \mL(fg)(Y_s), \ t \geq 0\bigg) 
\end{equation*}
is a $\rP_{x}$ martingale by definition of the generator $\mL$. Thus,
the process:
\begin{equation*}
\bigg(\frac{(fg)(Y_t)}{g(x)}-f(x)-\int_{0}^{t} ds \ \frac{\mL(fg)(Y_s)}{g(x)}, \ t \geq 0\bigg) 
\end{equation*}
is a $\rP_{x}$ martingale. We set: 
\begin{multline}
   \label{generatorLh1}
M^{f,g}_t=\expp{-\int_{0}^{t}  ds \ (\mL g/g)(Y_s)}
\frac{(fg)(Y_t)}{g(x)}-f(x)\\
-\int_{0}^{t} ds \ \expp{-\int_{0}^{s}
  dr \
  (\mL g/g)(Y_r)}  \left[ \frac{\mL(fg)(Y_s)}{g(x)} - \frac{\mL(g)(Y_s)}{g(Y_s)}
\frac{(fg)(Y_s)}{g(x)} \right].
\end{multline}
It\^o's lemma  then yields  that the process  $(M^{f,g}_t, t\geq  0)$ is
another $\rP_{x}$ martingale. Take $f$  constant equal to $1$ to get the
result. 
\end{proof}

Let $\rP_{x}^{g}$ denote the probability measure on $(D, \D)$ defined by:
\begin{equation}
\label{ACequation}
\forall t \geq 0, \quad \frac{d\rP^{g}_{x \ |\D_t}} {d\rP_{x  \ |\D_t}} = \frac{g(Y_t)}{g(x)}  \expp{-\int_{0}^{t} ds \ (\mL g/g)(Y_s)}.
\end{equation}
Note that in the case where  $g$ is harmonic for the linear operator $\mL$
(that  is  $\mL g=0$),  the probability distribution $\rP^{g}$  is the  usual  Doob
$h$-transform of $\rP$ for $h=g$. 

We also introduce the generator $\mL^{g}$ of the canonical process $Y$
under $\rP^g$ and the expectation operator $\rE^g$ associated to $\rP^{g}$. 

\begin{lem} 
\label{generatorLh}
Let  $g$  be  a  positive   function  of  $\D(\mL)$  such that $g$ is bounded from below by a positive constant. Then, we have $\D_g(\mL)  \subset \D(\mL^{g})$ and 
\begin{equation*}
% \label{defLh}
\forall u \in \D_g(\mL), \quad  \mL^{g}(u)= \frac{\mL(gu)-\mL(g)u}{g}\cdot
\end{equation*}
\end{lem}

\begin{proof}
  As, for $f\in \D_g(\mL)$, the  process $(M^{f,g}_t, t\geq 0)$ defined by
  \reff{generatorLh1}  is a martingale  under $\rP_x$,  we get  that the
  process:
\begin{equation*}
f(Y_t)-f(x)-\int_{0}^{t} ds \
\bigg(\frac{\mL(fg)(Y_s)-\mL(g)(Y_s)f(Y_s)}{g(Y_s)}  \bigg) , \quad t \geq 0
\end{equation*}
is a $\rP_{x}^{g}$ martingale. This gives the result. 
\end{proof}

\begin{rem}
\label{remPgt}
Let $\left((t,x) \to g(t,x)\right)$ be a function bounded from below by a positive constant, differentiable in $t$, such that $g(t,.)\in \D(\mL)$ for each $t$  and $\left((t,x) \to \partial_t g(t,x)\right)$ is bounded from above.
By considering the process $(t,Y_t)$ instead of $Y_t$, we have the immediate counterpart of Lemma  \ref{AClemma} for time dependent function $g(t, .)$. In particular, we may define the following probability measure on $(D, \D)$ (still denoted $\rP^g_{x}$ by a small abuse of notations):
\begin{equation}
\label{ACequation2}
\forall t \geq 0, \quad \frac{d\rP^{g}_{x \ |\D_t}} {d\rP_{x  \ |\D_t}} = \frac{g(t,Y_t)}{g(0,x)}  \expp{-\int_{0}^{t} ds \ \frac{\mL g+\partial_t g}{g}(s,Y_s)},
\end{equation}
where $\mL$ acts on $g$ as a function of $x$.
\end{rem}

We  now define the  $h$-transform for  superprocesses, as  introduced in
\cite{EP99} (notice  this does not correspond to  the Doob $h$-transform
for superprocesses).

\begin{defi}
\label{htransform}
Let $X=(X_t, t\geq 0)$ be an $(\mL,\beta,\alpha)$ superprocess. For $g \in
b\Ec^{+}$,  we   define  the  $h$-transform  of  $X$   (with  $h=g$)  as
$X^{g}=(X^g_t, t\geq 0)$ the measure valued process given for all $t\geq
0$ by:
\begin{equation}
   \label{eq:defXg}
X^{g}_t (dx) = g(x) X_t (dx). 
\end{equation}
\end{defi}

Note that \reff{eq:defXg} holds point-wise, and that the law of the $h$-transform of a superprocess may be singular with respect to the law of the initial superprocess.

We first give an easy generalization of a result in section 2 of
\cite{EP99} for a  general spatial motion.

\begin{prop}
\label{prophtransform}
Let $g$ be a positive function of $\D(\mL)$ such that $g$ is bounded from below by a positive constant. Then
the  process  $X^{g}$  is  a  $\left(\mL^{g}  ,  \frac{(-\mL+\beta)g}{g}  ,\alpha
g\right)$-superprocess.
\end{prop}

\begin{proof}
The Markov property of $X^{g}$ is clear. We compute, for $f \in b\Ec^{+}$ :
\[
\E_{{x}}[ \expp{-X_{t}^{g}(f)}]
=\E_{\delta_x/g(x)}[ \expp{-X_{t}(f g)}]
= \expp{-u_t(x)/g(x)},
\]
where, by Theorem \ref{theo1a}, $u$ satisfies:
\begin{equation}
\label{mildequation2}
u_t(x) + \rE_{x}\bigg[\int_{0}^{t} dr \ \psi(Y_r,u_{t-r}(Y_r))\bigg]= \rE_{x}\big[(fg)(Y_t) \big],
\end{equation}
which can also be written:
\begin{equation*}
u_t(x) + \rE_{x}\bigg[\int_{0}^{s} dr \ \psi(Y_r,u_{t-r}(Y_r))\bigg] + \rE_{x}\bigg[\int_{s}^{t} dr \ \psi(Y_r,u_{t-r}(Y_r))\bigg] = \rE_{x}\big[(fg)(Y_t) \big].
\end{equation*}
But \reff{mildequation2} written at time $t-s$ gives:
\begin{equation*}
u_{t-s}(x) + \rE_{x}\bigg[\int_{0}^{t-s} dr \ \psi(Y_r,u_{t-s-r}(Y_r))\bigg]= \rE_{x}\big[(fg)(Y_{t-s}) \big]. 
\end{equation*}
By comparing the two previous equations, we get:
\begin{equation*}
u_t(x) + \rE_{x}\bigg[\int_{0}^{s} dr \ \psi(Y_r,u_{t-r}(Y_r))\bigg]= \rE_{x}\big[u_{t-s}(Y_s) \big],
\end{equation*}
and the Markov property now implies that the process:
\begin{equation*}
u_{t-s}(Y_s) - \int_{0}^{s} dr \ \psi(Y_r,u_{t-r}(Y_r)) 
\end{equation*}
with $s\in [0,t]$ is a $\rP_x$ martingale.
It\^o's lemma now yields that the process:
\begin{equation*}
u_{t-s}(Y_s)\expp{-\int_{0}^{s} dr (\mL g/g)(Y_r)}  - \int_{0}^{s}  dr \
\expp{-\int_{0}^{r} du \ (\mL g/g)(Y_u)}
\big(\psi(Y_r,u_{t-r}(Y_r))-(\mL g/g)(Y_r) \ u_{t-r}(Y_r)\big) 
\end{equation*}
with $s\in [0,t]$ is another $\rP_{x}$ martingale (the integrability comes
from  the  assumption  $\mL g\in   \C_b$  and  $1/g  \in  \C_{b}$).  Taking
expectations at time $s=0$ and at time $s=t$, we have:
\begin{multline*}
 u_{t}(x) + \rE_{x} \bigg[ \int_{0}^{t} ds \ \expp{-\int_{0}^{s} dr
   (\mL g/g)(Y_r)} \big(\psi(Y_s,u_{t-s}(Y_s))-(\mL g/g)(Y_s)
 u_{t-s}(Y_s)\big) \bigg] \\ 
=  \rE_{x}\bigg[\expp{-\int_{0}^{t} dr (\mL g/g)(Y_r)} (fg)(Y_t)\bigg].   
\end{multline*}
We divide both sides by $g(x)$ and expand $\psi$ according to its definition:
\begin{multline*}
\big(\frac{u_{t}}{g}\big)(x) + \rE_{x} \bigg[ \int_{0}^{t} ds \
\frac{g(Y_s)}{g(x)} \expp{-\int_{0}^{s} dr (\mL g/g)(Y_r)}  \bigg((\alpha g)
(Y_s) \big( \frac{u_{t-s}}{g}\big)^2 (Y_s) + (\beta-\frac{\mL g}{g})(Y_s)
\big(\frac{u_{t-s}}{g}\big)(Y_s) \bigg)  \bigg] \\ 
=  \rE_{x}\bigg[\frac{g(Y_t)}{g(x)} \expp{-\int_{0}^{t} dr (\mL g/g)(Y_r)}  f(Y_t) \bigg].
\end{multline*}
By definition of $\rP_{x}^g$ from \reff{ACequation}, we get that:
\begin{equation*}
\big(\frac{u_{t}}{g}\big)(x) + \rE_{x}^{g} \bigg[ \int_{0}^{t} ds \
\bigg( (\alpha g) (Y_s) \big( \frac{u_{t-s}}{g}\big)^2 (Y_s) +
(\beta-\frac{\mL g}{g})(Y_s) \big(\frac{u_{t-s}}{g} \big)(Y_s) \bigg]
\bigg)  =  \rE_{x}^{g}\big[f(Y_t) \big]. 
\end{equation*}
We conclude from Theorem \ref{theo1a} that $X^{g}$ is a $(\mL^{g} ,
\frac{(-\mL+\beta)g}{g} ,\alpha g)$-superprocess. 
\end{proof}

In order to perform the $h$-transform of interest, we shall consider the
following assumption. 

\medskip
 \noindent
$(H2)$ \textbf{$1/\alpha$ belongs to $ \D(\mL)$.}
\medskip

Notice that $(H2)$ implies that  $\alpha \mL(1/\alpha) \in \C_{b}$. 
Proposition \ref{prophtransform} and Lemma \ref{AClemma} then yield the
following Corollary. 
\begin{cor}
\label{corhtransform}
Let $X$ be an $(\mL,\beta,\alpha)$-superprocess. Assume $(H2)$. The process
$X^{1/\alpha}$ is an $(\tilde{\mL}, \tilde{\beta}, 1)$-superprocess with: 
\begin{equation}
\label{betatilde}
\tilde{\mL}= \mL^{1/\alpha} \quad\text{and}\quad 
 \tilde{\beta}=\beta - \alpha \mL( 1/\alpha).
\end{equation}
Moreover, for all $t\geq 0$, the law $\tilde{\rP}_{x}$ of the process
$Y$ with
generator $\tilde{\mL}$ is absolutely continuous on $\D_t$ with respect to
$\rP_{x}$ and its Radon-Nikodym derivative is given by: 
\begin{equation}
\label{ACequationspec}
\frac{d\tilde{\rP}_{x \ |\D_t}} {d\rP_{x  \ |\D_t}}  =
\frac{\alpha(x)}{\alpha(Y_{t})}  \expp{\int_{0}^{t} ds \
  (\tilde{\beta}-\beta) (Y_s)}. 
\end{equation}
\end{cor}

We will note $\tilde{\P}$ for the law of  $X^{1/\alpha}$ on the
canonical space (that is $\tilde{\P}=\P^{(\tilde \mL, \tilde \beta, 1)}$)  and
$\tilde{\N}$ for its canonical measure. 
Observe that the branching mechanism of $X$ under $\tilde \P$, which we
shall write $\tilde \psi$, is given by:
\begin{equation}
\label{tildepsi}
\tilde \psi(x,\lambda)=  \tilde{\beta}(x) \ \lambda + \lambda^2,
\end{equation}
and the quadratic  coefficient is no more dependent  on $x$. Notice that
$\P_{\alpha \nu} (X \in \cdot)=\tilde{\P}_\nu (\alpha X\in \cdot)$. This
implies the following relationship on the canonical measures (use
Theorem \ref{canonical} to check it): 
\begin{equation}
\label{eqnormalisation2}
\alpha(x) \N_{x}[ X\in \cdot]= \tilde{\N}_{x}[\alpha X\in \cdot].  
\end{equation}
Recall that $v_{t}(x) = \N_{x}[\Hm>t]= \N_x[X_t\neq 0]$. We set
$\tilde v_t(x) = \tilde \N_x[X_t\neq 0]$.
As $\alpha$ is positive, equality \reff{eqnormalisation2} implies in
particular that, for all $t 
>0$ and $x \in E$: 
\begin{equation}
\label{eqnormalisation}
\alpha(x) v_{t}(x)= \tilde{v}_{t}(x).  
\end{equation}

\subsection{A Girsanov type theorem}
\label{subsectiongirsanov}

The following assumption  will be used to perform the Girsanov
change of measure. 

\medskip
 \noindent
$(H3)$   \textbf{Assume $(H2)$ holds. The function $\tilde{\beta}$ defined in 
    \reff{betatilde} is in  $\D(\tilde{\mL})$, with $\tilde \mL$
   defined  in 
  \reff{betatilde}.}
\medskip
 
For $z\in \R$, we set $z_+=\max(z, 0)$. Under $(H2)$ and $(H3)$, we define:
\begin{equation}
   \label{eq:defbeta0}
\beta_{0} = \sup_{x \in E} \max\bigg(\tilde{\beta}(x), \sqrt{ (\tilde{\beta}^2(x)-2
  \tilde{\mL}(\tilde{\beta})(x))_+}\bigg)\quad
\text{and}\quad 
\qq(x) = \frac{\beta_0-\tilde{\beta}(x)}{2}\cdot
\end{equation}
Notice that $\qq\geq 0$. 

We shall consider the distribution of the homogeneous
$(\tilde{\mL},\beta_0,1)$-superprocess, which we will denote by $\P^{0}$
($\P^{0} = \P^{(\tilde{\mL},\beta_0,1)}$) and its canonical measure
 $\N^{0}$. Note that the branching mechanism of $X$ under $\P^0$ is
 homogeneous (the branching mechanism does not depend on $x$). 
We set $\psi^0$ for $\psi^{\beta_0,1}$. Since
$\psi^0$ does not depend anymore on $x$ we shall also write
$\psi^0(\lambda)$ for $\psi^0(x,\lambda)$:
\begin{equation}
\label{psi0}
\psi^{0}(\lambda)=\beta_0\lambda+\lambda^2.
\end{equation}

\medskip Proposition  \ref{girsanov} below is a  Girsanov's type theorem
which allows  us to finally reduce  the distribution $\tilde  \P$ to the
homogeneous  distribution  $\P^0$.  We  introduce  the process  $M=(M_t,
t\geq 0)$ defined by:
\begin{equation} 
\label{aq:defM}
M_t= \exp{\bigg( X_{0}(\qq) -X_t(\qq) - \int_{0}^{t} ds \ X_s(\varphi)
  \bigg) }, 
\end{equation}
where the function $\varphi$ is defined by:
\begin{equation}
   \label{eq:deftildephi}
 \varphi (x)=\tilde \psi (x,\qq(x))- \tilde{\mL}(\qq)(x), \quad
x\in E. 
\end{equation}

\begin{prop}{A Girsanov's type theorem.}
\label{girsanov}
Assume  $(H2)$ and $(H3)$ hold.
Let $X$ be a $(\tilde{\mL}, \tilde{\beta}, 1)$-superprocess.
\begin{itemize}
 \item[(i)] The process $M$ is a bounded $\F$-martingale under
   $\tilde{\P}_{\nu}$ which converges a.s. to 
\[
M_{\infty}=\expp{X_0(q) -\int_0^{+\infty
    } ds\; X_s (\varphi)} \ind_{\{\Hm<+\infty \}}.
\]
 \item[(ii)] We have:
\[
   \frac{d\P^{0}_{\nu}} {d\tilde{\P}_{\nu }}  =
M_{\infty}.
\]
\item[(iii)]  If   moreover  $(H1)$  holds,   then $\P_{\nu}^{0}$-a.s.  we
  have $M_{\infty} >0$, the probability measure 
 $\tilde{\P}_{\nu}$  is
  absolutely continuous  with respect  to $\P^{0}_{\nu}$ on  $\F$:   
\[
   \frac{d\tilde{\P}_{\nu}} {d\P^{0}_{\nu }}  =
   \frac{1}{M_{\infty}},\quad \text{and}\quad
   \frac{d\tilde{\N}_{x}} {d\N^{0}_{x }}  =
   \expp{\int_0^{+\infty
    } ds\; X_s (\varphi)} . 
\]
We also have:
\begin{equation}
   \label{eq:q=NXpsi}
q(x)=\N^0_x\left[\expp{\int_0^{+\infty
    } ds\; X_s (\varphi)} -1 \right].
\end{equation}
\end{itemize}
\end{prop}

The two first points are a particular case of Theorem IV.1.6 p.252 in
\cite{PE99} on interactive drift. For the sake of completeness, we
give  a proof based on the mild form of the Log Laplace equation
\reff{mildequation} introduced in Section \ref{section0}. Notice that:
\begin{equation}
   \label{eq:psi-b0}
 \psi^{0}(\lambda)=\tilde \psi (x,\lambda+\qq(x))-\tilde \psi (x,\qq(x)).
\end{equation}
Thus,  Proposition  \ref{girsanov}  appears as  a non-homogeneous generalization  of  Corollary 4.4  in  \cite{AD09}.   We  first give  an
elementary Lemma.

\begin{lem}
   \label{lem:psi>0}
Assume  $(H2)$ and $(H3)$ hold.
The function $\varphi$ defined by  \reff{eq:deftildephi} is non-negative.
\end{lem}
\begin{proof}
The following computation:
\begin{align*}
\varphi(x)
=\tilde \psi (x,\qq(x))- \tilde{\mL}(\qq)(x) 
&= \qq(x)^2 + \tilde{\beta} \qq(x) - \tilde{\mL}(\qq)(x) \\
&= \bigg(\frac{\beta_0-\tilde{\beta}(x)}{2} \bigg)^2 + \tilde{\beta}(x) \
\frac{\beta_0-\tilde{\beta}(x)}{2} - \tilde{\mL}(\qq)(x) \\ 
&= \frac{\beta_{0}^2-\tilde{\beta}^2 (x) + 2 \tilde{\mL}(\tilde{\beta})(x)}{4} \\
\end{align*}
and the definition \reff{eq:defbeta0} of $\beta_{0}$  ensure that the function 
$\varphi$ is non-negative.  
\end{proof}

\begin{proof}[Proof of Proposition \ref{girsanov}]
First observe that $M$ is $\F$-adapted. 
As the function $\qq$ also is non-negative, we deduce from Lemma
\ref{lem:psi>0}  that the process $M$ is bounded by $\expp{X_{0}(\qq)}$.

Let $f \in b\Ec^{+}$. On the one hand, we have:
\[
\tilde{\E}_{x}[M_t  \expp{-X_t(f)}]
= \tilde{\E}_{x}[ \expp{ \qq(x)
  -X_t(\qq+f)-\int_{0}^{t} ds \ X_s(\varphi)} ] 
=  \expp{\qq(x)-r_t(x)},
\]
where, according to Theorem \ref{theo1a}, $r_t(x)$ is bounded on $[0,T]  \times E$ for all $T>0$ and satisfies:
\begin{multline}
\label{eqn10}
r_t(x) + \tilde{\rE}_{x} \bigg[ \int_{0}^{t} ds \
\tilde \psi (Y_{t-s},r_{s}(Y_{t-s})) \bigg] \\
= \tilde{\rE}_{x} \bigg[\int_{0}^{t}  ds \ \big( \tilde \psi (Y_{t-s},\qq(Y_{t-s})) - \tilde{\mL}(\qq) (Y_{t-s}) \big) + (\qq+f) (Y_t) \bigg].
\end{multline}
On the other hand, we have
\[ 
\E^{0}_{x}[ \expp{-X_t(f)}]=  \expp{-w_t(x)},
\] 
where $w_t(x)$ is bounded on $[0,T]  \times E$ for all $T>0$ and satisfies:
\[
w_t(x) + \tilde{\rE}_{x} \bigg[ \int_{0}^{t} ds \
\psi^{0}(Y_{t-s},w_{s}(Y_{t-s})) \bigg] =\tilde{\rE}_{x}[f(Y_t)].
\]
Using \reff{eq:psi-b0}, 
rewrite the previous equation under the form:
\begin{equation}
\label{girsanovmartingale4}
w_t(x) + \tilde{\rE}_{x}\bigg[ \int_{0}^{t} ds \
\tilde \psi (Y_{t-s},(w_{s} +\qq)(Y_{t-s})) \bigg] 
=  \tilde{\rE}_{x} \bigg[ \int_{0}^{t}  ds \ \tilde \psi (Y_{t-s},\qq(Y_{t-s}))  + f(Y_t) \bigg].
\end{equation}
We now make use of the Dynkin's formula with $(H3)$:
\begin{equation}
\label{girsanovmartingale5}
\qq(x) =- \tilde{\rE}_{x}\bigg[\int_{0}^{t} \tilde{\mL}(\qq)(Y_{s})\bigg] + \tilde{\rE}_{x}[\qq(Y_t)],
\end{equation}
and sum the equations \reff{girsanovmartingale4} and
\reff{girsanovmartingale5} term by term to get: 
\begin{multline}
\label{eqn10a}
(w_t+\qq)(x) + \tilde{\rE}_{x}\bigg[ \int_{0}^{t} ds \
\tilde \psi (Y_{t-s},(w_{s}+\qq)(Y_{t-s}))\bigg]  \\ 
=  \tilde{\rE}_{x} \bigg[\int_{0}^{t}  ds \
\big(\tilde \psi (Y_{t-s},\qq(Y_{t-s})) - \tilde{\mL}(\qq)
(Y_{t-s})\big) + (\qq+f) (Y_t) \bigg]. 
\end{multline} 

The functions $r_t(x)$ and $w_t(x)+\qq(x)$ are bounded on $[0,T] \times
E$ for all $T>0$ and satisfy the same equation, see equations
\reff{eqn10} and \reff{eqn10a}. By uniqueness, see Theorem \ref{theo1},
we finally get that $w_t+q=r_t$. This gives:
\begin{equation}
\label{girsanovmartingale3}
\tilde{\E}_{x}[M_t \expp{-X_t(f)}] = \E^{0}_{x}[ \expp{-X_t(f)}].
\end{equation}
The Poissonian decomposition of the superprocesses, see Theorem
\ref{canonical}, and the exponential formula enable us to extend this
relation to arbitrary initial measures $\nu$: 
\begin{equation}
\label{girsanovmartingale23}
\tilde{\E}_{\nu}[M_t  \expp{-X_t(f)}] = \E^{0}_{\nu}[ \expp{-X_t(f)}].
\end{equation}
This equality with $f=1$ and the  Markov property of $X$  proves  the
first part of item (i). 

Now, a direct induction based on the Markov property yields that, for
all positive integer $n$, and  $f_1, \ldots, f_n \in b\Ec^{+}$, $0 \leq
s_1 \leq \ldots \leq s_n \leq t$: 
\begin{equation}
\tilde{\E}_{\nu}[M_t \expp{- \sum_{1 \leq i \leq n} X_{s_i}(f_i)}] =
\E^{0}_{\nu}[ \expp{-\sum_{1 \leq i \leq n} X_{s_i}(f_i)}]. 
\end{equation}
And we conclude with an application of the monotone class theorem that,
for all non-negative $\F_t$-measurable random variable $Z$:
\begin{equation*}
\tilde{\E}_{\nu}[M_t Z] = \E^{0}_{\nu}[Z].
\end{equation*}
The martingale $M$ is bounded and thus converges a.s.  to a limit
$M_\infty $. 
We deduce that for all non-negative $\F_t$-measurable random variable $Z$:
\begin{equation}
\label{eq:EMi=E0}
\tilde{\E}_{\nu}[M_{\infty}Z] = \E^{0}_{\nu}[Z].
\end{equation}
This also holds for any non-negative $\F_\infty $-measurable random
variable $Z$. 
This gives the second item (ii). 

On  $\{\Hm<+\infty \}$,  then clearly  $M_t$ converges  to $\expp{X_0(q)
  -\int_0^{+\infty   }  ds\;   X_s   (\varphi)}   $.  Notice   that
$\P^0_\nu(\Hm=+\infty  )=0$.  We deduce  from  \reff{eq:EMi=E0} with
$Z=\ind_{\{\Hm=+\infty \}}$ 
that $\tilde \P_\nu$-a.s.  on
$\{\Hm=+\infty \}$, $M_\infty =0$. This gives the last part of item (i).

Now, we prove the third item (iii). Notice that \reff{eq:EMi=E0} implies
that $\P^0_\nu$-a.s. $M_\infty >0$. Thanks to $(H1)$, we also have that
$\tilde \P_{\nu}$-a.s.  $M_\infty >0$. Let $Z$ be a non-negative
$\F_\infty $-measurable random variable. Applying \reff{eq:EMi=E0} with  $Z$
replaced by $\ind_{\{M_\infty >0\}}Z/M_\infty $, we get:
\[
\tilde{\E}_{\nu}[Z]= \tilde{\E}_{\nu}\bigg[M_{\infty} \ind_{\{M_\infty
  >0\}} \frac{Z}{M_{\infty}}\bigg] =
\E^{0}_{\nu}\bigg[\frac{Z}{M_{\infty}}\ind_{\{M_\infty >0\}} \bigg] =
\E^{0}_{\nu}\bigg[\frac{Z}{M_{\infty}} \bigg].
\]
This gives the first part of item (iii). 

Notice  that for all  positive integer  $n$, and  $f_1, \ldots,  f_n \in
b\Ec^{+}$, $0 \leq s_1 \leq \ldots \leq s_n $, we have
\begin{align*}
   \tilde \N_x \left[1- \expp{\sum_{1 \leq i \leq n}
       X_{s_i}(f_i)}\right]
&=-\log \left(\tilde \rE_x \left[\expp{\sum_{1 \leq i \leq n}
       X_{s_i}(f_i)}\right]\right)\\
&=-\log \left( \rE^0_x \left[\expp{\sum_{1 \leq i \leq n}
       X_{s_i}(f_i)+ \int_0^{+\infty
    } ds\; X_s (\varphi) }\right]\right)+ q(x)\\
&=  \N^0_x \left[1- \expp{\sum_{1 \leq i \leq n}
       X_{s_i}(f_i)+ \int_0^{+\infty
    } ds\; X_s (\varphi)}\right] +q(x).
\end{align*}
Taking $f_i=0$ for all $i$ gives \reff{eq:q=NXpsi}. 
This implies:
\[
\tilde \N_x \left[1- \expp{\sum_{1 \leq i \leq n}
       X_{s_i}(f_i)}\right]
=  \N^0_x \left[\expp{\int_0^{+\infty
    } ds\; X_s (\varphi)} \left(1- \expp{\sum_{1 \leq i \leq n}
       X_{s_i}(f_i)}\right)\right].
\]
The monotone class theorem gives then the last part of item (iii). 
\end{proof}

\subsection{Genealogy for superprocesses}
\label{subsectionTransport}
We now recall the  genealogy of $X$ under $\P^{0}$ given by the Brownian
snake from \cite{DLG02}. We assume $(H2)$ and $(H3)$ hold. 

Let $\W$ denote the  set of all c\`adl\`ag killed paths in  $E$.  An element $w
\in  \W$ is  a c\`adl\`ag  path:  $w: [0,\eta(w))  \rightarrow E$,  with
$\eta(w)$ the lifetime  of the path $w$. By  convention the trivial path
$\{x\}$, with $x\in E$, is a  killed path with lifetime 0 and it belongs
to $\W$.  The space $\W$ is Polish for the distance:
\[ 
d(w,w')= \delta(w(0),w(0)') + |\eta(w)-\eta(w')| + \int_{0}^{\eta(w)
  \wedge \eta(w')} ds \ d_s(w_{[0,s]},w'_{[0,s]}),
\]
where $d_s$ refers  to the Skorokhod metric on  the space $D([0, s],E)$,
and  $w_{I}$  is   the  restriction  of   $w$  on   the  interval
$I$.  Denote  $\W_{x}$  the  set  of stopped  paths  $w$  such  that
$w(0)=x$.   We work on  the canonical  space of  continuous applications
from $[0,\infty)$ to  $\W$, denoted by $\bar{\Omega} \mathrel{\mathop:}=
\C(\R^{+}, \W)$, endowed  with the Borel sigma field  $\bar{\G}$ for the
distance   $d$,   and   the   canonical  right   continuous   filtration
$\bar{\G}_t=\sigma \{W_s, s \leq t \}$, where $(W_s,s\in \R^+)$ is the
canonical coordinate process. Notice $ \bar \G=\bar \G_\infty $ by
construction.  We set $H_s = \eta(W_s)$ the lifetime of
$W_s$.

\begin{defi}[Proposition 4.1.1 and Theorem 4.1.2 of \cite{DLG02}]
\label{def:W}
  Fix $W_{0}  \in \W_{x}$. There exists a  unique $\W_{x}$-valued Markov
  process   $W=(\bar{\Omega},\bar{\G},\bar{\G}_{t},W_t,\PP^{0}_{W_{0}})$,
  called the Brownian snake, starting  at $W_{0}$ and satisfying the two
  properties:
\begin{itemize}
\item[(i)]  The lifetime  process $H=(H_s  , s\geq  0)$ is  a reflecting
  Brownian  motion with  non-positive drift  $- \beta_0$,  starting from
  $H_{0}= \eta(W_{0})$.
 \item[(ii)] Conditionally given the lifetime process $H$, the
   process $(W_s, s\geq 0) $ is distributed as an inhomogeneous Markov process, with transition kernel specified by the two following prescriptions, for $0 \leq s \leq s'$:
\begin{itemize}
 \item $W_{s'}(t)=W_s(t)$ for all $t < H_{[s,s']}$, with  $H_{[s,s']}
   = \inf_{s \leq r \leq s'} H_r$. 
 \item Conditionally on $W_{s}(H_{[s,s']}-)$, the path
   $\big(W_{s'}(H_{[s,s']}+t),  0 \leq t <  H_{s'}- H_{[s,s']} \big)$ is
   independent of $W_s$ and is distributed as $Y_{[0,
     H_{s'}-H_{[s,s']})}$ under $\tilde{\rP}_{W_{s}(H_{[s,s']}-)}$. 
\end{itemize}
\end{itemize} 
This process will be called the $\beta_{0}$-snake started at $W_{0}$,
and its law  denoted by $\PP^{0}_{W_{0}}$.
\end{defi}
We will just  write $\PP^0_{x}$ for the law of the  snake started at the
trivial path $\{x\}$.  The  corresponding excursion measure $\NN^0_x$ of
$W$  is  given as  follows:  the  lifetime  process $H$  is  distributed
according to the  It\^o measure of the positive  excursion of a reflecting
Brownian motion  with non-positive drift $-  \beta_0$, and conditionally
given  the  lifetime process  $H$,  the process  $(W_s,  s\geq  0) $  is
distributed according to (ii) of Definition \ref{def:W}. Let 
\[
\sigma=\inf\{s>0; H_s=0\}
\]
denote the length of the excursion under $\NN^0_x$. 

Let  $(l^{r}_s, r  \geq 0, s \geq 0)$ be the bicontinuous version of the
local time process of $H$; where $l^{r}_s$ refers to the local time at
level $r$ at time $s$. 
We also set
$\hat{w}=w(\eta(w)-)$ for  the left end position of the path $w$.
We consider the measure valued process $X(W)=(X_t(W), t\geq 0)$ defined
under $\NN^0_x$ by:
\begin{equation}
   \label{eq:defX(W)}
X_t(W)(dx)=\int_0^\sigma d_sl^t_s\; \delta_{\hat W_s} (dx).
\end{equation}

The $\beta_{0}$-snake gives the genealogy of the $(\tilde{\mL},\beta_0,1)$
superprocess in the following sense. 

\begin{prop}[\cite{DLG02}, Theorem 4.2.1]
We have:
\begin{itemize}
\item The process $X(W)$ is under $\NN^{0}_{x}$  distributed as $X$ under $\N^{0}_{x}$.
\item Let $\sum_{j \in \J} \delta_{(x^j,W^j)}$ be a Poisson point measure
on $E \times \bar{\Omega}$ with intensity $\nu(dx) \; 
\NN^{0}_{x}[dW]$. Then $\sum_{j \in \J} X(W^{j})$ is an 
$(\tilde \mL,\beta_{0},1)$-superprocess started at $\nu$.
\end{itemize}
\end{prop}
Notice that, under $\NN^0_x$, the extinction time of $X(W)$ is defined
by 
\[
\inf \{t; X_t(W)=0\}=\sup_{s\in [0, \sigma]} H_s,
\]
and we shall write this quantity  $\Hm$ or $\Hm(W)$ if we need to stress
the dependence in $W$. This notation is coherent with \reff{eq:Hmax}.

We now transport the genealogy of $X$ under $\N^0$ to a genealogy of $X$
under $\tilde \N$. In order to
simplify notations, we shall write $X$ for $X(W)$ when there is no
confusion. 

\begin{defi}
\label{snakemeasure}
Under (H1)-(H3), we define a  measure $\tilde{\NN}_{x}$ on $(\bar{\Omega},\bar{\G})$ by:
\[ 
\forall W \in \bar{\Omega}, \quad \frac{d \tilde{\NN}_{x }} {d \NN^{0}_{x
   } }(W) 
= \frac{d
  \tilde{\N}_{x }} {d \N^{0}_{x }} (X(W)) =
\frac{1}{M_{\infty}}=\expp{\int_0^{+\infty
    } ds\; X_s (\varphi)} . 
 \]
\end{defi}
Notice the second equality in the previous definition is the third item
of Proposition \ref{girsanov}.

At this point, the genealogy defined for $X$ under $\tilde{\N}_{x}$ will
give the genealogy of $X$ under $ \N$  up to a weight. We set 
\begin{equation}
   \label{eq:defNN}
\NN_{x }=\inv{\alpha(x)}\tilde{\NN}_{x }.
\end{equation}

\begin{prop}
\label{proplinkgenealogy}
We have:
\begin{itemize}
 \item[(i)]  $X$ under $\tilde{\NN}_{x}$ is distributed as $X$ under $\tilde{\N}_{x}$.
 \item[(ii)] The weighted process $X^\text{weight}=(X^\text{weight}_t,
   t\geq 0)$ with
\begin{equation}
   \label{eq:geneaX}
X^\text{weight}_t(dx)=\int_{0}^{\sigma} d_s l_s^{t} \
   \alpha(\hat{W}_s) \delta_{\hat{W}_s} (dx), \quad t \geq 0 ,
\end{equation} is
   under $\NN_{x}$  distributed as $X$ under $\N_{x}$.
\end{itemize}
\end{prop}
We may write $X^\text{weight}(W)$ for $X^\text{weight}$ to emphasize
the dependence in the snake $W$. 
\begin{proof}
This is a direct consequence of Definition \ref{snakemeasure} and
\reff{eqnormalisation2}. 
\end{proof}

We shall say that $W$ under ${\NN}_{x}$ provides through
\reff{eq:geneaX} a genealogy for $X$ under $\N_x$. 

% We shall consider $\tau$ the right continuous
% inverse of the local time at level $0$: 
%  \[
% \tau_{s}=\inf{\{r \geq 0 ; l^{0}_r > s \}}.
% \]

\section{A Williams' decomposition}
\label{sectionwilliams}
In Section \ref{sec:bismut}, we give a decomposition of the genealogy of
the  superprocesses $(\mL,  \beta,\alpha)$  and $(\tilde{\mL},\tilde{\beta},
1)$  with   respect  to  a  randomly  chosen   individual.   In  Section
\ref{sec:williams}, we  give a Williams' decomposition  of the genealogy
of the superprocesses $(\mL, \beta,\alpha)$ and $(\tilde{\mL},\tilde{\beta},
1)$ with respect to the last individual alive.
\subsection{Bismut's decomposition}
\label{sec:bismut}
A decomposition  of the genealogy  of the homogeneous  superprocess with
respect to a randomly chosen individual is well known in the homogeneous
case, even for a general branching mechanism (see lemmas 4.2.5 and 4.6.1
in \cite{DLG02}).

We now  explain how to decompose  the snake process  under the excursion
measure ($\tilde  \NN_x$ or  $\NN^0_x$) with respect  to its value  at a
given time.  Recall $\sigma= \inf{\{  s >0, H_s=0 \}}$ denote the length
of the  excursion.  Fix a real  number $t \in  [0,\sigma]$.  We consider
the  process  $H^{(g)}$ (on  the  left of  $t$)  defined  on $[0,t]$  by
$H^{(g)}_s  =  H_{t-s}-H_{t}$ for  all  $s  \in  [0,t]$.  The  excursion
intervals above 0 of the  process $(H^{(g)}_{s}- \inf_{0 \leq s' \leq s}
H^{(g)}_{s'},  0 \leq  s \leq  t)$ are  denoted  $\bigcup_{j\in J^{(g)}}
(c_j,d_j)$.  We  also consider  the process $H^{(d)}$  (on the  right of
$t$)  defined on $[0,  \sigma-t]$ by  $H^{(d)}_s =  H_{t+s}-H_{t}$.  The
excursion intervals  above 0 of the process $(H^{(d)}_{s}-  \inf_{0 \leq s'
  \leq  s}   H^{(d)}_{s'},  0  \leq   s  \leq  \sigma-t)$   are  denoted
$\bigcup_{j\in  J^{(d)}}  (c_j,d_j)$.   We   define  the  level  of  the
excursion $j$ as $s_j= H_{t-c_j}$ if $j\in J^{(g)}$ and $s_j= H_{t+c_j}$
if  $j\in  J^{(d)}$.    We  also  define  for  the   excursion  $j$  the
corresponding excursion of the snake: $W^j=(W^{j}_s, s\geq 0)$ as
\[
W^j_s(.)=
W_{t-(c_j+s)\wedge d_j }(.+s_j)\quad\text{if $j\in J^{(g)}$, and}  \quad
W^j_s(.)=
W_{t+(c_j+s)\wedge d_j }(.+s_j)\quad\text{if $j\in J^{(d)}$}. 
\]
We consider the following  two point
measures on $\R^+\times \bar \Omega$: for $ \varepsilon \in \{g,d \}$, 
\begin{equation}
\label{defM}
 R_t^{\varepsilon} = \sum_{j\in
  J^{(\varepsilon)}}  \delta_{(s_j,W^{j})}. 
\end{equation}
Notice that  under $\NN^0_x$ (and  under $\tilde \NN_x$ if  $(H1)$ holds),
the   process    $W$   can    be   reconstructed   from    the   triplet
$(W_t,R^g_{t},R^d_t)$.  We are interested in the probabilistic structure
of this  triplet, when  $t$ is chosen  according to
the Lebesgue measure on the excursion time interval of the snake.  Under
$ \NN^0$, this result is as a consequence of Lemmas 4.2.4 and 4.2.5 from
\cite{DLG02}. We recall this result in the next Proposition.

For a point measure $R=\sum_{j \in J} \delta_{(s_j,x_j)}$ on a space $\R
\times \cx$ and $A\subset \R$,  we shall consider the restriction of $R$
to $A\times \cx$ given by  $R_{A}= \sum_{j \in J}\ind_A(s_j)
\delta_{(s_j,x_j)}$.  

\begin{prop}[\cite{DLG02}, Lemmas 4.2.4 and 4.2.5]
   \label{prop:propbismutN0}
For every measurable non-negative function
$F$, the following formulas hold:
\begin{align}
\label{eq:tNW0MM}
{\NN}_x^0\bigg[\int_{0}^{\sigma} ds \ F(W_s, R^g_s, R^d_s) \bigg]
&= \int_{0}^{\infty} \expp{-\beta_0 r} dr \ \tilde{\rE}_{x} \bigg[
F(Y_{[0,r)}, \hat R_{[0,r)}^{B,g}, \hat R_{[0,r)}^{B,d}) \bigg], \\
\label{eq:tNW0MMloc}
{\NN}_x^0\bigg[\int_{0}^{\sigma} d_sl^t_s \ F(W_s, R_s^g, R_s^d) \bigg]
&=  \expp{-\beta_0 t}  \ \tilde{\rE}_{x} \bigg[
F(Y_{[0,t)}, \hat R_{[0,t)}^{B,g}, \hat R_{[0,t)}^{B,d}) \bigg], \quad t> 0, 
\end{align}
where  under  $\tilde \rE_x$  and  conditionally  on $Y$, $ \hat
R^{B,g}$ and $\hat R^{B,d} $ are two independent  Poisson point
measures   with   intensity   $   \hat\nu^{B}(ds,dW)=  
ds   \; 
\NN_{Y_s}^0[dW]$.
\end{prop}

The next Proposition gives a similar result in the non-homogeneous case. 

\begin{prop}
\label{propbismut2}
Under $(H1)$-$(H3)$, for every measurable non-negative function
$F$, the two formulas hold:
\begin{equation}
\label{eq:tNWMM}
\tilde{\NN}_x\bigg[\int_{0}^{\sigma} ds \ F(W_s, R^g_s, R^d_s) \bigg]
= \int_{0}^{\infty} dr \ \tilde{\rE}_{x} \bigg[  \expp{- \int_{0}^{r} ds
  \ \tilde{\beta}(Y_s) } F(Y_{[0,r)},  R_{[0,r)}^{B,g}, R_{[0,r)}^{B,d}) \bigg],
\end{equation}
where under $\tilde \rE_x$ and conditionally on $Y$, $R_{[0,r)}^{B,g}$ and
$R_{[0,r)}^{B,d}$ are two independent Poisson point measures with 
intensity 
\begin{equation}
   \label{eq:def-nuB}
\nu^{B}(ds,dW)= ds \; \tilde \NN_{Y_s}[dW]
= ds \; \alpha(Y_s)  \NN_{Y_s}[dW];  
\end{equation} 
and
\begin{equation}
\label{eq:tNWMMweight}
{\NN}_x\bigg[\int_{0}^{\sigma} ds \ \alpha(\hat{W}_s) F(W_s, R_s^g,
R_s^d) \bigg] \\ 
= \int_{0}^{\infty} dr \ {\rE}_{x} \bigg[  \expp{- \int_{0}^{r} ds
  \ {\beta}(Y_s) } F(Y_{[0,r)},  R_{[0,r)}^{B,g}, R_{[0,r)}^{B,d}) \bigg],
\end{equation}
where under $ \rE_x$ and conditionally on $Y$, $R_{[0,r)}^{B,g}$ and
$R_{[0,r)}^{B,d}$ are two independent Poisson point measures with 
intensity  $ \nu^{B}$. 
\end{prop}

Observe there  is a  weight $\alpha(\hat W_s)$  in \reff{eq:tNWMMweight}
(see also \reff{eq:geneaX} where this weight appears) which modifies the
law of the individual picked  at random, changing the modified diffusion
$\tilde{\rP}_{x}$ in \reff{eq:tNWMM} into the original one $\rP_{x}$.

We shall use the following elementary Lemma on Poisson point measure.
\begin{lem} 
   \label{lem:GirsanovPoisson}
Let  $R$ be  a Poisson point measure on a Polish space with intensity
$\nu$. Let $f$ be a non-negative measurable function $f$  such that 
$\nu(\expp{f}-1)<+\infty $. Then for any non-negative measurable function
$F$, we have:
\begin{equation}
   \label{eq:GirsanovPoisson}
\E\left[F(R)\expp{R(f)}\right]=\E\left[F(\tilde R)\right]\expp{
  \nu(\expp{f}-1)},
\end{equation}
where $\tilde R$ is a Poisson point measure with intensity $ \tilde
\nu(dx)=\expp{f(x)} \nu(dx)$. 
\end{lem}

\begin{proof}[Proof of Proposition \ref{propbismut2}]
We keep notations introduced in Propositions \ref{prop:propbismutN0} and
\ref{propbismut2}. We have: 
\begin{multline*}
\tilde{\NN}_x\bigg[\int_{0}^{\sigma} ds \ F(W_s, R_s^g, R_s^d) \bigg]\\
\begin{aligned}
&= \NN^{0}_x\bigg[  \expp{\int_0^{+\infty
    } ds\; X_s (\varphi)}   \int_{0}^{\sigma} ds \ F(W_s, R_s^g,
  R_s^d) \bigg] \\ 
&= \NN^{0}_x\bigg[  \int_{0}^{\sigma} ds \ F(W_s, R_s^g,
R_s^d)\expp{(R_s^g+R_s^d)(f)} \bigg] \\
&= \int_{0}^{\infty} \!\! \expp{- \beta_0 r} dr \ \tilde{\rE}_{x} \bigg[
F(Y_{[0,r)}, \hat R_{[0,r)}^{B,g}, \hat R_{[0,r)}^{B,d})\expp{(\hat
  R_{[0,r)}^{B,g}+\hat R_{[0,r)}^{B,d})(f)} \bigg]\\ 
&= \int_{0}^{\infty} \!\! \expp{- \beta_0 r}dr \ \tilde{\rE}_{x} \bigg[
F(Y_{[0,r)},  
R_{[0,r)}^{B,g}, R_{[0,r)}^{B,d})\expp{2 \int_0^r ds\;
  \NN^0_{Y_s}[\expp{\int_0^{+\infty 
    } X_r(W)(\varphi)} -1 ]}\bigg]\\
&= \int_{0}^{\infty} \!\! \expp{- \beta_0 r} dr \ \tilde{\rE}_{x} \bigg[
F(Y_{[0,r)},  
R_{[0,r)}^{B,g},  R_{[0,r)}^{B,d})\expp{2 \int_0^r ds\; q(Y_s) }\bigg]\\
&= \int_{0}^{\infty} dr \ \tilde{\rE}_{x} \bigg[
F(Y_{[0,r)}, 
R_{[0,r)}^{B,g}, R_{[0,r)}^{B,d})\expp{- \int_0^r ds\; \tilde \beta (Y_s) }\bigg],
\end{aligned}
\end{multline*}
where the first equality comes from  $(H1)$ and item (iii) of Proposition
\ref{girsanov}, we set $f(s,W)=\int_0^{+\infty } X_r(W)(\varphi) $
for the second equality, we use Proposition \ref{prop:propbismutN0} for
the third equality, we use Lemma \ref{lem:GirsanovPoisson} for the
fourth, we use \reff{eq:q=NXpsi}
for the fifth, and the definition \reff{eq:defbeta0} of $q$ in the last. 
This proves \reff{eq:tNWMM}. 

Then replace  $F(W_s, R_s^g, R_s^d)$ by $\alpha(\hat  W_s) F(W_s, R_s^g,
R_s^d)$  in \reff{eq:tNWMM}  and  use \reff{ACequationspec}  as well  as
\reff{eq:defNN} to get \reff{eq:tNWMMweight}.
\end{proof}

The proof of the following Proposition is similar to the proof of
Proposition \ref{propbismut2} and is not reproduced here. 
\begin{prop}
\label{propbismut3}
Under $(H1)$-$(H3)$,  for every measurable non-negative function
$F$, the two formulas hold: for fixed $t> 0$,
\begin{equation}
\label{eq:tNWMMt}
\tilde{\NN}_x\bigg[\int_{0}^{\sigma} d_sl^t_s \ F(W_s, R_s^g, R_s^d) \bigg]
= \tilde{\rE}_{x} \bigg[  \expp{- \int_{0}^{t} ds
  \ \tilde{\beta}(Y_s) } F(Y_{[0,t)},  R_{[0,t)}^{B,g}, R_{[0,t)}^{B,d}) \bigg],
\end{equation}
where under $\tilde \rE_x$ and conditionally on $Y$, $R^{B,g}$ and
$R^{B,d}$ are two independent Poisson point measures with 
intensity  $ \nu^{B}$ defined in \reff{eq:def-nuB}, and
\begin{equation}
\label{eq:tNWMMweightt}
{\NN}_x\bigg[\int_{0}^{\sigma} d_sl^t_s \ \alpha(\hat{W}_s) F(W_s,
R_s^g, R_s^d) \bigg] \\ 
= {\rE}_{x} \bigg[  \expp{- \int_{0}^{t} ds
  \ {\beta}(Y_s) } F(Y_{[0,t)}, R_{[0,t)}^{B,g},  R_{[0,t)}^{B,d}) \bigg],
\end{equation}
where under $ \rE_x$ and conditionally on $Y$, $ R^{B,g}$ and
$ R^{B,d}$ are two independent Poisson point measures with 
intensity  $\nu^{B}$. 
\end{prop}
As an example of application of this Proposition, we can recover easily the  following well known result.  
\begin{cor}
\label{propbismut4}
Under $(H1)$-$(H3)$, for every measurable non-negative functions
$f$ and $g$ on $E$, we have:
\[
\N_x \bigg[ X_{t}(f) \expp{-X_{t}(g)} \bigg] = \ \rE_{x} \bigg[ \expp{-
  \int_{0}^{t} ds \ \partial_{\lambda} \psi \big(Y_s, \ \N_{Y_s}\big[1-
  \expp{X_{t-s}(g)}\big]\big) } f(Y_t) \bigg]. 
\]
\end{cor}
In particular, we recover the so-called ``many-to-one'' formula (with
$g=0$ in Corollary \ref{propbismut4}):
\begin{equation}
\label{manytoone}
\N_x [ X_{t}(f)] = \ \rE_{x} \bigg[ \expp{-
  \int_{0}^{t} ds \ \beta(Y_s) } f(Y_t) \bigg]. 
\end{equation}

\begin{rem}
\label{rem:Bismut}
Equation \reff{manytoone} justifies the introduction of the following family of probability measures indexed by $t \geq 0$:
\begin{align}
\label{eq:bismutspine}
\frac{d \rP_{x \ | \D_{t}}^{(B,t)}}{d \rP_{x \   | \D_{t}}} = \frac{\expp{-\int_{0}^{t} ds \; \beta(Y_s)}}{\rE_{x}\left[\expp{-\int_{0}^{t} ds \; \beta(Y_s)}\right]},
\end{align}
which can be understood as the law of the ancestral lineage of an individual sampled at random at height $t$ under the excursion measure $\N_{x}$, and also correspond to Feynman Kac penalization of the original spatial motion $\rP_x$ (see \cite{RY09}). Notice that this law does not depend on the parameter $\alpha$.  
These probability measures are not compatible as $t$ varies but will be shown in Lemma \ref{lem:convbismutspine} to converge as $t \to \infty$ in restriction to $\D_{s}$, $s$ fixed, $s \leq t$, under some ergodic assumption (see $(H9)$ in Section \ref{sec:H}).
\end{rem} 

\begin{proof}
We set for  $w\in \W$ with $\eta(w)=t$ and $r_1,r_2$ two point measures
on $\R^+\times \bar \Omega$
\[
F(w,r_1,r_2) =f(\hat w)\expp{h(r_1)+ h(r_2)},
\]
where $h(\sum_{i\in I} \delta_{(s_i, W^i)})=\sum_{s_i<t}
X^\text{weight}(W^i)_{t-s_i}(g)$.  
We have:
\begin{align*}
   \N_x \bigg[ X_{t}(f) \expp{-X_{t}(g)} \bigg]
&={\NN}_x\bigg[\int_{0}^{\sigma} d_sl^t_s \ \alpha(\hat{W}_s) F(W_s,
R_s^g, R_s^d) \bigg] \\ 
&= {\rE}_{x} \bigg[  \expp{- \int_{0}^{t} ds
  \ {\beta}(Y_s) } f(Y_t) \expp{ h(R_{[0,r)}^{B,g})+h(R_{[0,r)}^{B,d})}
  \bigg]\\
&= {\rE}_{x} \bigg[  \expp{- \int_{0}^{t} ds
  \ {\beta}(Y_s) } f(Y_t) \expp{ -\int_0^t 2 \alpha(Y_s) \NN_{Y_s}[1-\expp{
    X^\text{weight}_{t-s}(g)}]}\bigg]\\ 
&= \rE_{x} \bigg[ \expp{-
  \int_{0}^{t} ds \ \partial_{\lambda} \psi \big(Y_s, \ \N_{Y_s}\big[1-
  \expp{X_{t-s}(g)}\big]\big) } f(Y_t) \bigg],
\end{align*}
where we used item (ii) of Proposition \ref{proplinkgenealogy} for the
first and last equality, 
\reff{eq:tNWMMweightt} with  $F$ previously defined for the second, 
formula for exponentials of Poisson point measure and \reff{eq:defNN}
 for the third. 
\end{proof}

\subsection{Williams' decomposition}
\label{sec:williams}
We first recall the Williams'  decomposition for the Brownian snake (see
\cite{WI74} for  Brownian excursions, \cite{SE96} for  Brownian snake or
\cite{AD08}  for  general  homogeneous  branching
mechanism without spatial motion).

Under  the excursion  measures  $\NN^0_x$, $\tilde  \NN_x$ and  $\NN_x$,
recall that  $\Hm=\sup_{[0,\sigma]} H_s$.  Because of the  continuity of
$H$, we can define  $\Tm=\inf\{s>0, H_s=\Hm\}$. Notice the properties of
the Brownian excursions implies that a.e. $H_s=\Hm$ only if $s=\Tm$. 
We set $v^0_t(x)=\N_{x}^0[\Hm >  t]$ and recall this function does not
depend on $x$. Thus, we shall write $v^0_t$ for $v^0_t(x)$. Standard
computations give: 
\[
v^0_t=\frac{\beta_0}{\expp{\beta_0 t} -1}\cdot
\]
The next
result is a straightforward adaptation from Theorem  3.3 of \cite{AD08} and gives the distribution of $(\Hm, W_\Tm, R_\Tm^g,
R_\Tm^d)$ under ${\NN}^0_x$. 

\begin{prop}[Williams' decomposition under ${\NN}^0_x$]
   \label{prop:WilliamN0}
We have:
\begin{itemize}
 \item[(i)] The distribution of $\Hm$ under ${\NN}^0_x$ is characterized
   by: ${\NN}^0_{x}[\Hm > h]= v_h^0$. 
 \item[(ii)] Conditionally on $\{\Hm=h_0\}$, the law of $W_\Tm$ under
   ${\NN}_x^0$ is distributed as $Y_{[0,h_0)}$ under $\tilde \rP_{x}$.
 \item[(iii)] Conditionally on $\{\Hm=h_0\}$ and $W_\Tm$, $R_\Tm^g$ and
   $R_\Tm^d $ are under ${\NN}^0_x$ independent Poisson point measures on
   $\R^+ \times \bar{\Omega}$ with intensity:
\[ 
\ind_{[0,h_0)}(s) ds \; \ind_{\{\Hm(W) < h_0 -s\}} \;{\NN}^0_{W_\Tm(s)}[ 
dW].
\]
\end{itemize}
In other words, for any non-negative measurable function $F$, we have
\[
\NN^{0}_x\bigg[F(\Hm, W_\Tm, R_\Tm^g, R_\Tm^d)\bigg] 
= - \int_{0}^{\infty} \partial_h v^0_h \; dh\;   \tilde{\rE}_{x} \bigg[
F(h,Y_{[0,h)}, \hat R^{W,(h),g}, \hat R^{W,(h),d}) \bigg],
\]
where under $\tilde \rE_x$ and conditionally on $Y_{[0,h)}$, $\hat R^{W,(h),g}$ and
$\hat R^{W,(h),d}$ are two independent Poisson point measures with 
intensity  $ \hat \nu^{W,(h)}(ds,dW)= \ind_{[0,h)}(s) ds \;
\ind_{\{\Hm(W) < h -s\}} \;\NN_{Y_s}^0[ dW]$. 
\end{prop}
Notice that items (ii) and (iii) in the previous Proposition implies the
existence of a measurable family $(\NN^{0,(h)}_x, h>0)$ of
probabilities on $(\bar \Omega,  \bar \G)$ such that $\NN^{0,(h)}_x$ is
the distribution of $W$ (more precisely of $( W_\Tm, R_\Tm^g, R_\Tm^d)$)
under ${\NN}^0_x$ conditionally on  $\{\Hm=h\}$. 

% We take the opportunity to note the following time reversal property of the height process of the homogeneous snake, and its consequence on the following remark.
% \begin{cor}
% \label{Esty}
% We have:
% 
% \end{cor}

% \begin{proof}
% Note the process $\left( H_s, 0 \leq s \leq \Tm \right)$ is a Bessel process up to its first hitting time of $t$: the reversibility property (i) is well known (see Exercice 0.2.5 of \cite{PIT}) and can be established by approximating this process by a suitably conditioned simple random walk. The second point readily follows.
% \end{proof}

\begin{rem}
\label{Esty}
In Klebaner \& al \cite{KRS07}, the Esty time reversal ``is obtained by conditioning a [discrete time] Galton Watson process in negative time upon entering state $0$ (extinction) at time $0$ when starting at state $1$ at time $-n$ and letting $n$ tend to infinity''. The authors then observe that in the linear fractional case (modified geometric offspring distribution) the Esty time reversal has the law of the same Galton Watson process conditioned on non extinction. 
Notice that in our continuous setting, the process $\left( H_s, 0 \leq s \leq \Tm \right)$ is under $\NN_{x}^{0,(h)}$ a Bessel process up to its first hitting time of $h$, and thus is reversible: $\left( H_s, 0 \leq s \leq \Tm \right)$ under $\NN_x^{0,(h)}$ is distributed as $\left( h-H_{\Tm-s}, 0 \leq s \leq \Tm \right)$ under $\NN_x^{0,(h)}$.
It is also well known (see Corollary 3.1.6 of \cite{DLG02}) that $\left( H_{\sigma-s}, 0 \leq s \leq \sigma-\Tm \right)$ under $\NN_x^{0,(h)}$ is distributed as $\left( H_{s}, 0 \leq s \leq \Tm \right)$ under $\NN_x^{0,(h)}$. We deduce from these two points that $\left( X_s(1), 0 \leq s \leq h \right)$ under $\NN_x^{0,(h)}$ is distributed as $\left( X_{h-s}(1), 0 \leq s \leq h \right)$ under $\NN_x^{0,(h)}$.
This result, which holds at fixed $h$, gives a pre-limiting version of the Esty time reversal in continuous time. Passing to the limit as $h \to \infty$, see Section \ref{backward}, we  get the equivalent of the Esty time reversal in a continuous setting.
\end{rem}
% \begin{proof}
% Decomposing  $\left( H_s, 0 \leq s \leq \Tm \right)$ into the excursions above its future infimum yields the William's decomposition stated above. It is also possible to decompose $\left( H_s, 0 \leq s \leq \Tm \right)$ into the excursions $H_i$ under its current supremum; let us denote $(H_i, i\in \I)$ these excursions and $(s_i,i \in \I)$ the associated values of the supremum. The measure $\sum_{i \in \I} \delta_{(s_i, H_i)}$ is still a Poisson point measure on $\R^{+} \times \Omega$, but with intensity $\ind_{(0,h]}(s) ds \; \ind_{\{ \Hm \leq s\}} \NN^{0}(-H \in .)$. Reversing the time, we find out that the excursions above its future infimum of the process $\left( h-H_{\Tm-s}, 0 \leq s \leq \Tm \right)$ have the same law than that of the process $\left( H_s, 0 \leq s \leq \Tm \right)$. This yields point (i). The same argument applies to the right hand side of the excursion, $\left( H_{\sigma-s}, 0 \leq s \leq \Tm-\sigma \right)$, and point (ii) readily follows.
% \end{proof}

Before stating the  Williams' decomposition, Theorem \ref{theowilliams},
let us prove  some properties for the functions $  v_t(x) = \N_{x}[\Hm >
t]=\N_x[X_t\neq 0]$  and $\tilde v_t(x)=  \tilde \N_{x}[\Hm >  t]$ which
will play a significant r\^ole in the next Section. Recall
\reff{eqnormalisation} states that 
\[
\alpha v_t=\tilde v_t. 
\]

Notice also that \reff{eq:defbeta0} implies that $q$ is bounded from
above by $(\beta_0+\norm{\tilde
  \beta}_{\infty })/2$.

\begin{lem}
\label{EDPV}
Assume $(H1)$-$(H3)$. We have:
\begin{equation}
   \label{eq:majovtilde}
q(x)+ v^0_t\geq  \tilde v_t(x) \geq  v^0_t. 
\end{equation}
Furthermore for fixed $x \in E$, $\tilde v_t(x)$ is of class $\C^1$ in
$t$ and we have:
\begin{equation}
   \label{eq:derivvtilde}
\partial_t \tilde v_t(x)=\tilde \rE_x
\left[\expp{\int_0^t \Sigma_r(Y_{t-r}) \; dr}\right]\; \partial_t v^0_t,
\end{equation}
where the function $\Sigma$ defined by:
\begin{equation}
   \label{eq:defSigma}
\Sigma_t(x)= 2(v^0_t + q(x) - \tilde v_t(x)) \\
=\partial_{\lambda} \psi^{0}(v^{0}_t) - \partial_{\lambda} \tilde \psi
(x,\tilde{v}_{t}(x))  
\end{equation}
satisfies:
\begin{equation}
   \label{eq:ineqSigma}
0\leq \Sigma_t(x) \leq 2 q(x)\leq  \beta_0+\norm{\tilde
  \beta}_{\infty }. 
\end{equation}
\end{lem}

\begin{proof}
 We deduce from item (iii) of Proposition \ref{girsanov} that, as
 $\varphi\geq 0$ (see Lemma \ref{lem:psi>0}), 
\[
\tilde v_t(x)=\tilde \N_x[X_t\neq 0]=\N^0_x \left[\ind_{\{X_t\neq 0\}}
\expp{\int_0^ {+\infty } ds\; X_s(\varphi) }\right]
\geq \N_x^0[X_t\neq 0]=v^0_t.
\]
We also have
\begin{align*}
\tilde v_t(x)
&=\N^0_x \left[\ind_{\{X_t\neq 0\}}
\expp{\int_0^ {+\infty } ds\; X_s(\varphi) }\right]\\
&=\N^0_x \left[\expp{\int_0^ {+\infty } ds\; X_s(\varphi) } -1\right] 
+ \N^0_x \left[1- 
\ind_{\{X_t= 0\}}
\expp{\int_0^ {+\infty } ds\; X_s(\varphi) }\right]\\
&= q(x) + \N^0_x \left[1- 
\ind_{\{X_t= 0\}}
\expp{\int_0^ {+\infty } ds\; X_s(\varphi) }\right]\\
&\leq q(x) + \N^0_x \left[1- 
\ind_{\{X_t= 0\}}\right]\\
&=q(x) + v_t^0,
\end{align*}
where we used \reff{eq:q=NXpsi} for the third equality. This proves
\reff{eq:majovtilde}. 

Using the Williams' decomposition under $\NN^0_x$, we get:
\[
\tilde v_t(x)=-\int_t^{+\infty }\partial_r v_r^0 \; dr\; \NN^{0,(r)}_x
\left[ \expp{\int_0^{+\infty } ds\; X_s(\varphi)}\right]. 
\]
Using again the Williams' decomposition under $\NN^0_x$, we have
\begin{align}
\nonumber
\NN^{0,(r)}_x
\left[ \expp{\int_0^{+\infty } ds\; X_s(\varphi)}\right]
&=\tilde \rE_x    \left[\expp{2 \int _0^ r ds \; \NN^0_{Y_{r-s}} \left[
      (\expp{\int_0^{+\infty } dt\; X_t(\varphi)}
      -1)\ind_{\{X_s=0\}} \right]}\right]\\
\label{eq:NN0r}
&=\tilde \rE_x    \left[\expp{2 \int _0^ r ds \; \NN^0_{Y_{s}} \left[
      (\expp{\int_0^{+\infty } dt\; X_t(\varphi)}
      -1)\ind_{\{X_{r-s}=0\}} \right]}\right].
\end{align}
We deduce that, for fixed $x$,  $r\mapsto \NN^{0,(r)}_x
\left[ \expp{\int_0^{+\infty } ds\; X_s(\varphi)}\right]$ is
non-decreasing and continuous as $\NN^0_y[\Hm=t]=0$ for $t>0$. 
Therefore, we deduce that for fixed $x$, $\tilde v_t(x)$ is of class
$\C^1$ in $t$:
\[
\partial_t \tilde v_t(x)=\NN^{0,(t)}_x
\left[ \expp{\int_0^{+\infty } ds\; X_s(\varphi)}\right]
\; \partial_t v^0_t.
\]
We have thanks to item (iii) from Proposition \ref{girsanov}: 
\begin{multline}
\label{eq:NN0-Sigma}
\NN^0_{y} \left[
      (\expp{\int_0^{+\infty } dt\; X_t(\varphi)}
      -1)\ind_{\{X_{s}=0\}} \right]
\\
\begin{aligned}
 &  = \NN^0_{y} \left[X_s \neq 0 \right] + 
\NN^0_{y} \left[\expp{\int_0^{+\infty } dt\; X_t(\varphi)}
      -1\right] - \NN^0_{y} \left[\expp{\int_0^{+\infty } dt\;
        X_t(\varphi)} \ind_{\{X_s\neq 0\}}  \right]\\
&= v^0_s + q(y) - \tilde v_s(y)\\
&= \frac{1}{2} \bigg[ \partial_{\lambda} \psi^{0}(v_s^{0})- \partial_{\lambda} \tilde \psi (y,\tilde v_{s}(y)) \bigg],
\end{aligned}   
\end{multline}
where the last equality follows from \reff{tildepsi}, \reff{eq:defbeta0} and \reff{psi0}. Thus, with $\Sigma_{s}(y)=\partial_{\lambda} \psi^{0}(v_s^{0})- \partial_{\lambda} \tilde \psi (y,\tilde v_{s}(y))$,  we deduce that: 
\[
\NN^{0,(t)}_x
\left[ \expp{\int_0^{+\infty } ds\; X_s(\varphi)}\right]
= \tilde \rE_x    \left[\expp{\int _0^ t ds \;
    \Sigma_{s}(Y_{t-s})}\right].
\]
This implies \reff{eq:derivvtilde}. Notice that, thanks to
\reff{eq:majovtilde},  $\Sigma$ is non-negative
 and bounded from above by $2q$. 
\end{proof}

Fix $h>0$. We define the probability measures $\rP^{(h)}$ absolutely continuous
with respect to $\rP$ and $\tilde \rP$  on $\D_h$ with Radon-Nikodym
derivative: 
\begin{equation}
\label{defPh}
\frac{d\rP^{(h)}_{x \ |\D_h}} {d\tilde \rP_{x  \ |\D_h}} 
=\frac{\expp{\int_0^h \Sigma_{h-r}(Y_{r}) \; dr}}
{\tilde \rE_x
\left[\expp{\int_0^h \Sigma_{h-s}(Y_{s}) \; dr}\right]}\cdot
\end{equation}

Notice  this Radon-Nikodym derivative  is 1  if the  branching mechanism
$\psi$ is homogeneous.  We deduce from \reff{eq:derivvtilde} and
\reff{eq:defSigma} that: 
\[
\frac{d\rP^{(h)}_{x \ |\D_h}} {d\tilde \rP_{x  \ |\D_h}} 
= \frac{\partial_h v^0_h}{\partial_h \tilde v_h(x)}
\expp{-\int_0^h dr\; \left( \partial_\lambda \tilde \psi (Y_r,\tilde
  v_{h-r}(Y_r)) - \partial_\lambda \psi^{0}(v^{0}_{h-r}) \right)}
\]
and, using \reff{ACequationspec}:
\begin{equation}
   \label{eq:RNderivPhP}
\frac{d\rP^{(h)}_{x \ |\D_h}} {d \rP_{x  \ |\D_h}} 
= \inv{\alpha(Y_h)} \frac{\partial_h v^0_h}{\partial_h
  v_h(x)} \expp{-\int_0^h dr\; \left( \partial_\lambda \psi(Y_r,
  v_{h-r}(Y_r))- \partial_\lambda \psi^{0}(v^{0}_{h-r}) \right)}.
\end{equation}

In the next Lemma, we give an intrinsic representation of the
Radon-Nikodym derivatives \reff{defPh} and \reff{eq:RNderivPhP}, which does not
involve $\beta_0$ or $v^0$. 
\begin{lem}
\label{AC2}
Assume $(H1)$-$(H3)$. Fix $h>0$. The processes $M^{(h)}=(M^{(h)}_t, t\in
[0,h))$ and  $\tilde M^{(h)}=(\tilde M^{(h)}_t, t\in
[0,h))$, with:
\[
 M^{(h)}_t= \frac{\partial_{h} v_{h-t}(Y_t)}{\partial_{h} v_{h}(x)}
\expp{- \int_{0}^{t}  ds \ \partial_{\lambda} \psi(Y_s,v_{h-s}(Y_s)) }
\quad\text{and}\quad
\tilde  M^{(h)}_t= \frac{\partial_{h} \tilde v_{h-t}(Y_t)}{\partial_{h}
  \tilde v_{h}(x)}
\expp{- \int_{0}^{t}  ds \ \partial_{\lambda} \tilde \psi 
  (Y_s,\tilde v_{h-s}(Y_s)) },
\]
are  non-negative bounded $\D_t$-martingales respectively under
$\rP_{x}$ and $\tilde \rP_x$. Furthermore, we have for $0 \leq t < h$:  
\begin{equation}
\label{defPh2}
\frac{d\rP^{(h)}_{x \ |\D_t}} {d\rP_{x  \ |\D_t}} =M^{(h)}_t
\quad\text{and}\quad
\frac{d\rP^{(h)}_{x \ |\D_t}} {d\tilde \rP_{x  \ |\D_t}} =\tilde M^{(h)}_t
. 
\end{equation}
\end{lem}
Notice the limit $M^{(h)}_h$ of $M^{(h)}$ and 
the limit $\tilde M^{(h)}_h$ of $\tilde M^{(h)}$ are respectively given by the right-handside
of \reff{eq:RNderivPhP} and \reff{defPh}. 

\begin{rem}
Comparing \reff{ACequation2} and \reff{defPh2}, we have that $\rP_{x}^{(h)}=\rP_{x}^{g}$ with $g(t,x)=\partial_{h} v_{h-t}(x)$, if $g$ satisfies the assumptions of Remark \ref{remPgt}.
\end{rem}

\begin{proof}
First of all, the process $\tilde M^{(h)}$ is clearly
$\D_t$-adapted. Using \reff{eq:derivvtilde}, we get:
\[
\tilde \rE_y \left[ \expp{\int_0^{h-t} \Sigma_{h-t-r}(Y_r) \; dr} \right]
=\frac{\partial_h \tilde v_{h-t} (y)}{\partial_h  v^0_{h-t} }\cdot
\]
We set: 
\[
\tilde M^{(h)}_h=\frac{\expp{\int_0^h \Sigma_{h-r}(Y_{r}) \; dr}}
{\tilde \rE_x
\left[\expp{\int_0^h \Sigma_{h-s}(Y_{s}) \; dr}\right]}\cdot
\]
We have:
\begin{align*}
   \tilde \rE_x [\tilde M^{(h)}_h | \D_t]
&=\frac{\expp{\int_0^t \Sigma_{h-r}(Y_{r}) \; dr}}
{\tilde \rE_x
\left[\expp{\int_0^h \Sigma_{h-s}(Y_{s}) \; dr}\right]}\tilde \rE_{Y_t} 
\left[ \expp{\int_0^{h-t} \Sigma_{h-t-r}(Y_r) \; dr} \right]\\
&=\frac{\partial_h \tilde v_{h-t} (Y_t)}{\partial_h \tilde v_{h} (x)}
\frac{ \partial_h  v^0_{h} }{\partial_h  v^0_{h-t} } \expp{\int_0^t
  \Sigma_{h-r}(Y_{r}) \; dr}\\ 
&=\frac{\partial_h \tilde v_{h-t} (Y_t)}{\partial_h \tilde v_{h} (x)}
\expp{- \int_{0}^{t} \partial_{\lambda} \tilde \psi \big(Y_s,\tilde v_{h-s}(Y_s) \big)  \; ds}
\frac{ \partial_h  v^0_{h} }{\partial_h  v^0_{h-t} } 
\expp{ \int_{0}^{t} \partial_{\lambda} \psi^{0} \big(v^{0}_{h-s}\big)  \; ds} 
\end{align*}
In the homogeneous setting, $v^{0}$ simply solves the ordinary
differential equation: 
\begin{align*}
\partial_h v^{0}_h = - \psi^{0}(v^{0}_h). 
\end{align*}
This implies that
\[
\partial_h \log(\partial_h v^0_h)=\frac{\partial^2_h v_h^0}{\partial_h
  v_h^0} =-\partial_{\lambda} \psi^0(v_h^0)
\]
and thus 
\begin{equation}
\label{eq:v0}
\frac{ \partial_h  v^0_{h} }{\partial_h  v^0_{h-t} } 
\expp{ \int_{0}^{t} \partial_{\lambda} \psi^{0} \big(v^{0}_{h-s}\big)
  \; ds}=1 . 
\end{equation}
We deduce that
\[
  \tilde \rE_x [\tilde M^{(h)}_h | \D_t]= 
\frac{\partial_h \tilde v_{h-t} (Y_t)}{\partial_h \tilde v_{h} (x)}
 \expp{-\int_0^t dr\; \partial_\lambda \tilde \psi  (Y_r, \tilde
v_{h-r}(Y_r)) }=\tilde M^{(h)}_t. 
\]
Therefore, $\tilde M^{(h)}$ is  a $\D_t$-martingale under $\tilde \rP_x$
and   the   second  part   of   \reff{defPh2}   is   a  consequence   of
\reff{defPh}. Then, use \reff{ACequationspec} to get that $M^{(h)}$ is a
$\D_t$-martingale under $\rP_x$ and the first part of \reff{defPh2}.
\end{proof}

We now give the Williams'  decomposition:  the distribution of
$(\Hm, W_\Tm, R_\Tm^g, R_\Tm^d)$  under $ {\NN}_x$ or equivalently under
$\tilde  \NN_x/\alpha(x)$.  Recall  the  distribution  $  \rP_{x}^{(h)}$
defined in \reff{defPh} or \reff{eq:RNderivPhP}.

\begin{theo}[Williams' decomposition under $\NN_x$]
\label{theowilliams}
Assume $(H1)$-$(H3)$. We have:
\begin{itemize}
 \item[(i)] The distribution of $\Hm$ under ${\NN}_x$ is characterized
   by: ${\NN}_{x}[\Hm > h]= v_h(x)$. 
 \item[(ii)] Conditionally on $\{\Hm=h_0\}$, the law of $W_\Tm$ under
   ${\NN}_x$ is distributed as $Y_{[0,h_0)}$ under $ \rP_{x}^{(h_0)}$.
 \item[(iii)] Conditionally on $ \{\Hm=h_0\}$ and $W_\Tm$, $R_\Tm^g$ and
   $R_\Tm^d $ are under ${\NN}_x$ independent Poisson point measures on
   $\R^+ \times \bar{\Omega}$ with intensity:
\[ 
\ind_{[0,h_0)}(s) ds \;\ind_{\{\Hm(W')
  < h_0 -s\}} \alpha(W_\Tm(s))\;  {\NN}_{W_\Tm(s)}[ 
dW'].
\]
\end{itemize}
In other words, for any non-negative measurable function $F$, we have
\[
\NN_x\bigg[F(\Hm, W_\Tm, R_\Tm^g, R_\Tm^d)\bigg] 
= - \int_{0}^{\infty} \partial_h v_h(x) \; dh\;   {\rE}_{x}^{(h)} \bigg[
F(h, Y_{[0,h)},   R^{W,(h),g},   R^{W,(h),d}) \bigg],
\]
where under $ \rE_x^{(h)}$ and conditionally on $Y_{[0,h)}$, $
R^{W,(h),g}$ and 
$R^{W,(h),d}$ are two independent Poisson point measures with 
intensity:
\begin{equation}
   \label{eq:def-nuW}
  \nu^{W,(h)}(ds,dW)= \ind_{[0,h)}(s) ds \; \ind_{\{\Hm(W) < h -s\}} 
\alpha(Y_s) \; \NN_{Y_s}[ dW].
\end{equation}
\end{theo}
Notice that items (ii) and (iii) in the previous Proposition imply the
existence of a measurable family $(\NN^{(h)}_x, h>0)$ of
probabilities on $(\bar \Omega,  \bar \G)$ such that $\NN^{(h)}_x$ is
the distribution of $W$  (more precisely of $( W_\Tm, R_\Tm^g,
R_\Tm^d)$) under ${\NN}_x$ conditionally on  $\{\Hm=h\}$. 

\begin{proof}
 We keep notations introduced in Proposition \ref{prop:WilliamN0} and
 Theorem \ref{theowilliams}.  We have:  
\begin{multline*}
\tilde{\NN}_x\bigg[ F(\Hm, W_\Tm, R_\Tm^g, R_\Tm^d) \bigg]\\
\begin{aligned}
&= \NN^{0}_x\bigg[  \expp{\int_0^{+\infty
    } ds\; X_s (\varphi)}  F(\Hm, W_\Tm, R_\Tm^g, R_\Tm^d)\bigg] \\ 
&= \NN^{0}_x\bigg[  F(\Hm, W_\Tm, R_\Tm^g, R_\Tm^d)
\expp{(R_\Tm^g+R_\Tm^d)(f)} \bigg] \\ 
&= - \int_{0}^{\infty} \partial_h v^0_h \; dh\;   \tilde{\rE}_{x} \bigg[
F(h,Y_{[0,h)}, \hat R^{W,(h),g}, \hat R^{W,(h),d})\expp{(\hat 
  R^{W,(h),g}+\hat R^{W,(h),d})(f)} \bigg]\\ 
&= - \int_{0}^{\infty} \partial_h v^0_h \; dh\;   \tilde{\rE}_{x} \bigg[
F(h,Y_{[0,h)},  R^{W,(h),g},  R^{W,(h),d})\expp{2\int_0^h
   ds\;  \NN_{Y_s}^0\left[ (\expp{\int_0^{+\infty } dt \;
       X_t(\varphi)} -1) \ind_{\{X_{r-s}=0\}} \right]}  \bigg]\\ 
&= - \int_{0}^{\infty} \partial_h v^0_h \; dh\;   \tilde{\rE}_{x} \bigg[
F(h,Y_{[0,h)}, R^{W,(h),g},   R^{W,(h),d})\expp{\int_0^h
  \Sigma_{h-s} (Y_s)  \; ds}  \bigg]\\ 
&= - \int_{0}^{\infty} \partial_h v^0_h \; dh\;  
 \tilde{\rE}_{x} \bigg[
\expp{\int_0^h
  \Sigma_{h-s} (Y_s)  \; ds}  \bigg] {\rE}_{x}^{(h)} \bigg[
F(h,Y_{[0,h)},  R^{W,(h),g},   R^{W,(h),d}) \bigg]\\ 
&= - \int_{0}^{\infty} \partial_h \tilde v_h(x) \; dh\;  
 {\rE}_{x}^{(h)} \bigg[
F(h,Y_{[0,h)},  R^{W,(h),g},  R^{W,(h),d}) \bigg],
\end{aligned}   
\end{multline*}
where the first equality comes from  $(H1)$ and item (iii) of Proposition
\ref{girsanov}; we set $f(s,W)=\int_0^{+\infty } X_r(W)(\varphi) $
for the second equality; we use Proposition \ref{prop:WilliamN0} for
the third equality; we use Lemma \ref{lem:GirsanovPoisson} for the
fourth with  $
R^{W,(h),g}$ and 
$R^{W,(h),d}$ which under $\tilde  \rE_x^{(h)}$ and conditionally on
$Y_{[0,h)}$ are two independent Poisson point measures with 
intensity  $  \nu^{W,(h)}$; we use  \reff{eq:NN0-Sigma}
for the fifth, definition \reff{defPh} of $\rE^{(h)}_x$
for the sixth, and \reff{eq:derivvtilde} for the seventh. 
Then use \reff{eq:defNN} and \reff{eqnormalisation} to conclude. 
\end{proof}

The definition of $\NN^{(h)}_x$ gives in turn  sense to the conditional law
$\N^{(h)}_x=\N_{x}(.|\Hm=h)$ of the 
$(\mL,\beta,\alpha)$ superprocess conditioned to die at time $h$, for all $h>0$. 
The next Corollary is then a straightforward consequence of Theorem
\ref{theowilliams}. 

\begin{cor}
   \label{cor:Nh}
Assume $(H1)$-$(H3)$. Let $h>0$. Let  $x\in  E$  and
$Y_{[0,h)}$  be distributed  according  to 
  $\rP^{(h)}_{x}$.  Consider  the Poisson point  measure $\cn=\sum_{j\in
    J}  \delta_{(s_j,X^{j})}$ on  $[0,h) \times  \Omega$  with intensity:  
\[
2\ind_{[0,  h)}(s)   ds  \;  \ind_{   \{  \Hm(X)  <   h-s  \}}
  \alpha(Y_s)\; \N_{Y_s}[ dX].
\]
 The process $X^{(h)}=(X^{(h)}_t, t\geq
  0)$, which is  defined  for all $t\geq 0$ by:
\[
X^{(h)}_t= \sum_{j\in J, \, s_j<t} X_{t-s_j }^{j},
\]
is distributed according to $\N^{(h)}_x$. 
\end{cor}

We now give the superprocess counterpart of Theorem
\ref{theowilliams}. 

\begin{cor}[Williams' decomposition under $\P_\nu$]
\label{corwilliams}
Assume $(H1)$-$(H3)$. We have the following result. 
\begin{itemize}
 \item[(i)] Sample a positive number $h_{0}$ according to the law of
   $\Hm$ under $\P_\nu$: $\P_{\nu}(\Hm \leq h)=  \expp{-\nu(v_{h})}$.
 \item[(ii)] Conditionally on $h_0$, sample $x_0 \in E$ according to
   the probability measure  
\[
\frac{\partial_{h} v_{h_{0}} (x)}{\nu(\partial_{h} v_{h_{0}})} \nu(dx).
\]
\item[(iii)]  Conditionally  on  $h_0$  and  $x_0$,  sample  $X^{(h_0)}$
  according to the probability measure $\N_{x_0}^{(h_0)} $.
\item[(iv)] Conditionally on $h_0$, sample $X'$, independent of
  $x_0$ and $X^{(h_0)}$, 
 according to the probability measure $\P_{\nu}(. |\Hm < h_{0})$.
\end{itemize}
Then the measure valued process $X'+X^{(h_0)}$  has distribution
$\P_{\nu}$. 
\end{cor}
In particular the distribution  of $X'+X^{(h_0)}$ conditionally on $h_0$
(which  is given  by (ii)-(iv)  from Corollary  \ref{corwilliams})  is a
regular   version  of   the  distribution   of   the  $(\mL,\beta,\alpha)$
superprocess conditioned  to die at a  fixed time $h_0$,  which we shall
write $\P_\nu^{(h_0)}$.

\begin{proof}
   Let $\mu$ be a finite  measure on $\R^+$ and $f$ a non-negative
   measurable function defined on $\R^+\times E$. For a measure-valued
   process $Z=(Z_t, t\geq 0)$ on $E$, we set $Z(f\mu)=\int f(t,x)\;
   Z_t(dx) \mu(dt)$. We also write $f_s(t,x)=f(s+t,x)$. 

   Let  $X'$ and  $X^{(h_0)}  $ be  defined as  in Corollary
   \ref{corwilliams}.  In order to characterized the distribution of 
the process $X'+X^{(h_0)}$, we  shall  compute 
\[
A=\E[\expp{ - X'(f\mu) - X^{(h_0) }(f\mu)}].
\]
We shall use notations from Corollary \ref{cor:Nh}.  We have:
\begin{align*}
   A
&= - \int_0^{+\infty }  \nu(\partial_h v_h) \expp{-\nu(h)}dh\;
\int_E \frac{\partial_{h} v_{h} (x)}{\nu(\partial_{h} v_{h})}\;  \nu(dx) \\
& \hspace{2cm}\rE_x^{(h)}\left[\E[\expp{ - \sum_{j\in J} X^j
    (f_{s_j}\mu)}|Y_{[0,h)}]\right]
 \E_\nu\left[\expp{ - X(f\mu)}| \Hm< h\right]\\
&= - \int_E\nu(dx) \; 
 \int_0^{+\infty }  \partial_{h} v_{h} (x)\; dh\\
& \hspace{2cm}
\rE_x^{(h)}\left[\E[\expp{ - \sum_{j\in J} X^j
    (f_{s_j}\mu)}|Y_{[0,h)}]\right]
 \E_\nu\left[\expp{ - X(f\mu)}\ind_{\{\Hm< h\}}\right],
\end{align*}
where we used  the definition of $X'$ and $\cn$  for the first equality,
and the equality $\P_{\nu}(\Hm < h)=\P_{\nu}(\Hm \leq h)=\expp{-\nu(h)}$
for the  second. Recall  notations from Theorem  \ref{theowilliams}.  We
set:
\[
G\left(\sum_{i\in I} \delta_{(s_i, W^i)},\sum_{i'\in I'} \delta_{(s_{i'},
    W^{i'})} \right)=
\expp{- \sum_{j \in I\cup I'} X^j(W^j)(f_{s_j}\mu)}
\]
and $g(h)=\E_\nu\left[\expp{ - X(f\mu)}\ind_{\{\Hm< h\}}\right]$. We
have:
\begin{align*}
   A
&= - \int_E\nu(dx) \; 
 \int_0^{+\infty }  \partial_{h} v_{h} (x)\; dh\; 
\rE_x^{(h)}\left[G(R^{W,(h),g}, R^{W,(h),d} )g(h) \right]\\
&=\int_E\nu(dx) \;  \NN_x\left[G(R_g^{\Tm}, R_d^{\Tm}) g(\Hm)\right]\\
&=\int_E\nu(dx) \;  \N_x\left[\expp{-X(f\mu)}
  \E_\nu\left[\expp{-X(f\mu)}\ind_{\{\Hm<h\}}  \right]_{|h=\Hm}\right]\\
&=\E\left[\sum_{i\in I} \expp{ - X^i(f\mu)} \prod_{j\in I; \;
      j\neq i}\expp{- X ^j(f\mu)} \ind_{ \{\Hm^j<\Hm^i\}}\right]\\
&=\E\left[\expp{ - \sum_{i\in I} X^i(f\mu)}\right]\\
&=\E_\nu \left[\expp{ -  X(f\mu)}\right],
\end{align*}
where we used  the definition of $G$ and $g$ for the first and third
equalities, Theorem
\ref{theowilliams} for the second equality, the master formula for
Poisson point measure $\sum_{i\in I} \delta_{X^i}$ with intensity 
$\nu(dx)\; \N_x[dX]$ for the fourth  equality (and the
obvious notation $\Hm^i=\inf\{t\geq 0; X^i_t=0\}$) and Theorem
\ref{canonical} for 
the last equality. Thus we get:
\[
\E[\expp{ - X'(f\mu) -  X^{(h_0) }(f\mu)}]=\E_\nu
\left[\expp{ -  X(f\mu)}\right].
\]
This readily implies that the  process $X'+ X^{(h_0)}$ is distributed as
$X$ under $\P_\nu$.
\end{proof}

\section{Some applications}

\subsection{The law of the Q-process}

Recall $\P^{(h)}_\nu$  defined after Corollary  \ref{corwilliams} is the
distribution of the $(\mL,\beta,\alpha)$-superprocess started at $\nu\in
\cm_f(E)$  conditionally on  $\{\Hm=h\}$.  We  consider  also $\P^{(\geq
  h)}_\nu=\P_\nu(\;\cdot\;      |\Hm\geq     h)$      and     $\N^{(\geq
  h)}_x=\N_x(\;\cdot\;   |\Hm\geq   h)$   the   distributions   of   the
$(\mL,\beta,\alpha)$-superprocess conditionally on $\{\Hm\geq h\}$.

The  distribution  of the  Q-process,  when  it  exists, is defined  as the  weak  limit of  $\P^{(\geq  h)}_\nu$ when  $h$ goes  to
infinity. The next Lemma insures that if $\P^{(h)}_\nu$ weakly converges
to  a  limit   $\P_\nu^{(\infty  )}$,  then  this  limit   is  also  the
distribution of the Q-process.

\begin{lem}
\label{lemmaconv}
Fix $t>0$.  If $\P_{\nu}^{(h)}$ (resp.  $\N_{x}^{(h)}$) converges weakly
to  $\P_{\nu}^{(\infty)}$  (resp.  $\N_{x}^{(\infty  )}$)  on  $(\Omega,
\F_t)$,  then   $\P_{\nu}^{(  \geq  h)}$   (resp.  $\N_{x}^{(\geq  h)}$)
converges weakly  to $\P_{\nu}^{(\infty)}$ (resp.  $\N_{x}^{(\infty )}$)
on $(\Omega, \F_t)$.
\end{lem}

\begin{proof}
Let $Z=\ind_{A}$  with $A\in\F_t$ such that $\P_{\nu}^{(\infty)}(\partial A)=0$. Using the Williams' decomposition under $\P_{\nu}$ given by Corollary
\ref{corwilliams}, we have for $h>t$:
\[
\E_{\nu}^{(\geq h)} [Z]
= \expp{\nu(v_h)}\int_{h}^{\infty} \E_{\nu}^{(h')} [Z] \; f(h') dh', 
\]
where $f(h)=- \nu(\partial_h v_h)\exp (-\nu(v_h))$. 
We write down the difference:
\[
\E_{\nu}^{(\geq h)}[Z] - \E_{\nu}^{(\infty)}[Z]
=  \expp{\nu(v_h)}\int_{h}^{\infty} \big(\E_{\nu}^{(h')}
[Z]-\E_{\nu}^{(\infty)}[Z]\big) \; f(h') dh'.
\]
Since $\P_{\nu}^{(h')}$    weakly     converges    to
$\P_{\nu}^{(\infty)}$ on  $(\Omega, \F_t)$ and since
$\P_{\nu}^{(\infty)}(\partial A)=0$, we deduce  that  
$\lim_{h'\rightarrow+\infty  } \E_{\nu}^{(h')}[Z]  -
\E_{\nu}^{(\infty)}[Z]=0$. We
conclude  that  $\lim_{h\rightarrow+\infty  } \E_{\nu}^{(\geq  h)}[Z]  -
\E_{\nu}^{(\infty)}[Z]=0$, which  gives the result. The proof  is
similar for the conditioned excursion measures.
\end{proof}

We now address the question  of convergence of the family of probability
measures $(\P_{x}^{(h)}, h \geq 0)$.  

Recall from \reff{defPh2} that for
all $0 \leq t < h$:
\[
\frac{d \rP^{(h)}_{x \ |\D_t}} {d \rP_{x  \ |\D_t}}=M_{t}^{(h)} .
\]
We shall consider the following assumption on the convergence in law of
the spine.

\medskip
 \noindent
$(H4)$ \textbf{For all $t\geq 0$, $\rP_{x}$-a.s. $(M^{(h)}_t, h>t)$
  converges to a limit say $M^{(\infty )}_t$, and
  $\rE_{x}[M_{t}^{(\infty)}]=1$.} 
\medskip

% We shall consider the following assumption on the convergence in law of
% the spine:
% 
% \medskip
%  \noindent
% $(H4)$ \textbf{For all $x \in E$, there exists a probability measure $\rP^{(\infty)}_{x}$ such that 
% for all $h >0$, $\rP_{x}^{(h)}$ is absolutely continuous with respect to $\rP_{x}^{(\infty)}$ on $\D_t$ for each $0 \leq t < h$ with Radon-Nikodym derivative $M_{t}^{(h),\infty}(x)  = \frac{d \rP^{(h)}_{x \ |\D_t}} {d \rP^{(\infty)}_{x  \ |\D_t}}$ satisfying:}
% \begin{equation}
% \label{defPinfty}
% \forall t \geq 0, M_{t}^{(h),\infty} \to 1  \mbox{ in } \rL^1(\rP^{(\infty)}_{x}) \mbox{ as } h \to \infty. 
% \end{equation}
% \medskip

% $t\geq 0$, $\rP_{x}$-a.s. $(M^{(h)}_t, h>t)$
%   converges to a limit say $M^{(\infty )}_t$, and
%   $\rE_{x}[M_{t}^{(\infty)}]=1$.} 

Note that Scheff\'e's  lemma implies that the convergence  also holds in
$\rL^{1}(\rP_{x})$.  
Furthermore, since  $(M^{(h)}_t,  t\in [0,  h))$ is  a
non-negative martingale,  there exists a  version of $(M^{(\infty
  )}_t, t\geq 0)$ which is a non-negative martingale.

% \begin{rem}
% The   assumptions   $(H1)$-$(H4)$   will   be   checked   in   Section
%  \ref{sec:twoexamples} in two special cases:  the multitype Feller diffusion and the superdiffusion, and under some additional assumptions, involving the
%   generalized eigenvalue $\lambda_{0}$ defined by:
% 
% and the associated eigenfunction. We also note that $(H4)$ obviously holds in the homogeneous setting for any spatial motion.
% \end{rem}

\begin{rem}
We provide in Section \ref{sec:twoexamples} sufficient conditions for $(H1)$-$(H4)$ to hold in the case of the multitype
Feller diffusion and the superdiffusion. These conditions are stated in term of the generalized eigenvalue
$\lambda_{0}$ defined by
\begin{equation}
 \label{eq:eigenvalue}
\lambda_0 = \sup{\{ \ell \in \R, \exists u\in \D(\mL), u > 0 \mbox{ such that } (\beta - \mL) u = \ell\  u\}},
\end{equation}
and its associated eigenfunction.
\end{rem}

% Nevertheless, these are just sufficient conditions, and $(H4)$ obviously holds in the homogeneous setting for any spatial motion. 
% Denoting by $\phi_{0}$ the eigenfunction associated to $\lambda_0$, we will prove that, if $\lambda_{0} \geq 0$ and the spatial motion $\rP^{\phi_{0}}$ (recall the definition \reff{ACequation}) satisfies an ergodic assumption, see $(H9)$ in Section \ref{sec:H}, assumptions  $(H1)$-$(H4)$  hold, see Lemmas \ref{lem1234} and \ref{lem123}.
% 
% However, this ergodic assumption is not necessary to get $(H4)$ as $(H4)$ holds in the homogeneous case for any spatial motion for instance.

% M_{t}^{(h),\infty}= \frac{\partial_{t} v_{h-t}(Y_{t})\expp{\lambda t} }{\partial_{t} v_{h}(Y_{0})}  \frac{\phi(Y_0)}{\phi(Y_t)} \expp{- \int_{0}^{t} ds \ 2 \alpha(Y_s) v_{h-s}(Y_s)},
% \] 
% with $Y$ the canonical process under $\rP_{x}$. In those two cases $Y$
% is ergodic, 

\begin{rem}
The family $(\rP_{x}^{(h)}, h \geq 0)$ and the family $(\rP_{x}^{(B,h)}, h \geq 0)$ defined in Remark \ref{rem:Bismut} will be shown in Lemma \ref{lem:convbismutspine} to converge to the same limiting probability measure.
\end{rem}
%under some additional ergodic assumption $(H9)$ defined in Section \ref{sec:H}}.

Under $(H4)$, we define the probability measure $\rP^{(\infty)}_{x}$ on
$(D,\D)$ by its Radon Nikodym derivative, for all $t\geq 0$:
\begin{align}
\label{defPinfty}
 \frac{d \rP_{x \ | \D_t}^{(\infty)}} {d \rP_{x \ | \D_t}} =
 M_{t}^{(\infty)}. 
\end{align}
By construction, the probability measure $\rP^{(h)}_{x}$ converges weakly
to $\rP^{(\infty)}_{x}$ on $\D_t$, for all $t \geq 0$. 

Let $\nu \in \M_{f}(E)$. We shall consider the following assumption:

\medskip
 \noindent
$(H5)_\nu$ \textbf{There exists a measurable function
$\rho$ such that the following convergence holds in  $\rL^1(\nu)$:}
\[
\frac{\partial_{h} v_{h}}{
  \nu(\partial_{h} v_{h})}
\; \xrightarrow[h \to +\infty]{}  \; \rho.
\] 
\medskip

In particular, we have $\nu(\rho)=1$. 
Let $\nu \in \M_{f}(E)$. Under $(H4)$ and $(H5)_\nu$, we set:
\[
 \rP^{(\infty )}_\nu(dY)=\int_E  \nu(dx) \rho(x)  \:  \rP^{(\infty)}_{x}(dY).
\]
Notice then that $ \int_E
\nu(dx)\; \frac{\partial_{h} v_{h}(x) }{
  \nu(\partial_{h} v_{h})} \rP^{(h)}_x(dY)$ converges weakly to
$\rP^{(\infty )}_\nu(dY)$ on $\D_t$,
for all $t\geq 0$.

\begin{rem}
   \label{rem:nu=Dirac}
If $\nu$ a constant times the Dirac mass  $\delta_x$, for some $x\in E$,
then $(H5)_\nu$ holds if $(H4)$ holds and  in this case we have
 $\rP^{(\infty )}_{\nu}=\rP^{(\infty)}_{x}$. 
\end{rem}

We can now state the result on the convergence of $\NN^{(h)}_x$. 

\begin{theo}
\label{theo-convN}
Assume $(H1)$-$(H4)$. Let  $t\geq 0$. The   triplet  $((W_{\Tm})_{[0,t]}, (R_\Tm^g)_{[0,t]},  (R_\Tm^d  )_{
  [0,t]})$ under  $\NN_x^{(h)}$ converges weakly to  the distribution of
the triplet  $(Y_{[0,t]}, R^{B,g}_{[0,t]}, R^{B,d}_{[0,t]}  )$ where $Y$
has   distribution  $\rP_{x}^{(\infty)}$   and  conditionally   on  $Y$,
$R^{B,g}$ and $R^{B,d}$ are  two independent Poisson point measures with
intensity $\nu^B$ given by  \reff{eq:def-nuB}. We even have the slightly
stronger result. For any bounded measurable function $F$, we have:
\begin{align}
\label{convNh}
\NN_x^{(h)}\bigg[F\big((W_{\Tm})_{[0,t]}, (R_\Tm^g)_{[0,t]},
  (R_\Tm^d )_{ [0,t]}\big)\bigg] 
\; \xrightarrow[h \to +\infty]{}  \; 
 \rE_x^{(\infty)} \bigg[
F\big(Y_{[0,t]}, R^{B,g}_{[0,t]},  R^{B,d}_{[0,t]} \big) \bigg].
\end{align}
\end{theo}

\begin{proof}
Let $h>t$. We use notations from Theorems \ref{theo-convN} and
\ref{theowilliams}. 
Let $F$ be a bounded measurable function on $\W\times (\R^{+} \times
\bar{\Omega})^2$. 
From the Williams' decomposition, Theorem \ref{theowilliams}, we have:
\begin{align*}
\NN_x^{(h)}\bigg[F((W_{\Tm})_{[0,t]}, (R_\Tm^g)_{[0,t]},
  (R_\Tm^d )_{ [0,t]})\bigg] 
& = \rE_x^{(h)} \bigg[
F\big(Y_{[0,t]}, R^{W,g,(h)}_{[0,t]},  R^{W,d,(h)}_{[0,t]} \big) \bigg]\\
& = \rE_{x}^{(h)} \big[  \varphi^{h}(Y_{[0,t]})\big],
\end{align*}
where $\varphi^{h}$ is defined by:
\[
\varphi^{h}(y_{[0,t]})
=\rE_x^{(h)} \bigg[
F\big(y_{[0,t]}, R^{W,g,(h)}_{[0,t]},  R^{W,d,(h)}_{[0,t]} \big)
\bigg|Y=y \bigg] .
\]
We also set:
\[
\varphi^{\infty }(y_{[0,t]})
= \rE_x^{(\infty)} \bigg[
F\big(y_{[0,t]}, R^{B,g}_{[0,t]},  R^{B,d}_{[0,t]} \big)
\bigg|Y=y \bigg] .
\]
We want to control:
\[
\Delta_{h}= \NN_x^{(h)}\bigg[F((W_{\Tm})_{[0,t]}, (R_\Tm^g)_{[0,t]},
  (R_\Tm^d )_{ [0,t]})\bigg]  - 
 \rE_x^{(\infty)} \bigg[
F(Y_{[0,t]}, R^{B,g}_{[0,t]},  R^{B,d}_{[0,t]} ) \bigg].
\]
Notice that:
\begin{align}
\nonumber
\Delta_{h}
&=  \rE_{x}^{(h)} \big[  \varphi^{h}(Y_{[0,t]})\big]  - \rE_{x}^{(\infty)}
\big[  \varphi^{\infty}(Y_{[0,t]})\big] \\ 
\label{eq:Delta}
&     =    \big(    \rE_{x}^{(h)}     \big[    \varphi^{h}(Y_{[0,t]})\big]-
\rE_{x}^{(\infty)} \big[ \varphi^{h}(Y_{[0,t]})\big]\big)+ \rE_{x}^{(\infty
  )} \big[(\varphi^{h}- \varphi^{\infty})(Y_{[0,t]})\big] .
\end{align} 
We prove the first  term of the right hand-side
of \reff{eq:Delta} converges to $0$. We have:
\[
 \rE_{x}^{(h)} \big[ \varphi^{h}(Y_{[0,t]})\big]- \rE_{x}^{(\infty)} \big[
 \varphi^{h}(Y_{[0,t]})\big]  = \rE_{x}^{(\infty)} \big[ (M_t^{(h)} -
 M_{t}^{(\infty)}) \: \varphi^{h}(Y_{[0,t]})\big].  
\]
Then use  that $\varphi^{h}$ is  bounded by $\norm{F}_\infty  $ and
the  convergence of $(M_{t}^{(h)},  h>t)$ towards  $M_{t}^{(\infty)}$ in
$\rL^1(\rP_x^{(\infty)})$ to get:
\begin{equation}
   \label{eq:lim2=0}
 \lim_{h \to \infty}\rE_{x}^{(h)} \big[  \varphi^{h}(Y_{[0,t]})\big]-
 \rE_{x}^{(\infty)} \big[  \varphi^{h}(Y_{[0,t]})\big]= 0. 
\end{equation}

\medskip
We then  prove the second term of the right hand-side
of \reff{eq:Delta} converges to $0$. 
Conditionally on $Y$, $R^{W,g,(h)}_{[0,t]}$ and
$R^{W,d,(h)}_{[0,t]} $
(resp. $R^{B,g}_{[0,t]}$ and $R^{B,d}_{[0,t]}$) are independent
Poisson point measures  with intensity  $\ind_{[0,t]}(s)\;
\nu^{W,(h)}(ds,dW)$ 
where $\nu^{W,(h)}$ is given by
\reff{eq:def-nuW} 
(resp. $\ind_{[0,t]}(s)\; \nu^B(ds,dW)$ where $\nu^B$ is given by
\reff{eq:def-nuB}). And we have:
\[
 \ind_{[0,t]}(s)\;
\nu^{W,(h)}(ds,dW) = \ind_{\{\Hm(W) < h -s\}} \ind_{[0,t]}(s)\;
\nu^B(ds,dW).
\]
Thanks to 
\reff{eqnormalisation} and \reff{eq:majovtilde}, we get that:
\[
\int \ind_{\{\Hm(W) \geq  h -s\}} \ind_{[0,t]}(s)\;
\nu^B(ds,dW)
= \int_0^t ds\; \alpha(y_s)\N_{y_s}[\Hm \geq  h -s]
= \int_0^t ds\; v_{h-s}(y_s)<+\infty .
\]

The proof of the next Lemma is postponed to the end of this Section. 
\begin{lem}
   \label{lem:cvPoisson}
Let  $R$ and $\tilde R$ be two  Poisson point measures on a Polish space
with respective intensity $\nu$ and $\tilde \nu$. 
Assume that $ \tilde 
\nu(dx)=\ind_A(x) \nu(dx)$, where $A$ is measurable and 
$\nu(A^c)<+\infty$. 
Then for any  bounded measurable function
$F$, we have:
\[
\val{\E[F(R)]-\E[F(\tilde R)]} \leq  2\norm{F}_\infty \nu(A^c).
\]
\end{lem}

Using this Lemma with $\nu$ given by $\ind_{[0,t]}(s)\;
\nu^B(ds,dW)$ and $A$ given by $\{\Hm(W) < h -s\}$, we deduce that:
\[
\val{(\varphi^{h}-\varphi^{\infty })(y_{[0,t]}) }
\leq  4 \norm{F}_\infty \int_0^t ds\; v_{h-s}(y_s).
\]
We deduce that:
\[
\val{\rE_{x}^{(\infty )} \big[(\varphi^{h}- \varphi^{\infty})(Y_{[0,t]})\big]}
\leq  4 \norm{F}_\infty \rE_{x}^{(\infty )} \left[\int_0^t ds\; 
  v_{h-s}(Y_s)\right] .
\]
Recall that $(H1)$ implies that $v_{h-s}(x) $ converges to 0 as $h$ goes
to infinity.  Since $v$ is bounded (use 
\reff{eqnormalisation} and \reff{eq:majovtilde}), by dominated
convergence, we get:
\begin{equation}
   \label{eq:lim1=0}
 \lim_{h \to \infty} \rE_{x}^{(\infty )} \big[(\varphi^{h}-
 \varphi^{\infty})(Y_{[0,t]})\big] =   0. 
\end{equation}
Therefore, we deduce from \reff{eq:Delta} that
$\lim_{h\rightarrow+\infty } \Delta_{h}=0$, which gives \reff{convNh}.
\end{proof}

We now define a superprocess with spine
distribution $ \rP^{(\infty )}_\nu$. 

\begin{defi}
   \label{defi:Q} 
Let $\nu\in \cm_f(E)$. Assume    $\rP^{(\infty)}_{\nu}$ 
 is well defined.  Let $Y$ be distributed according to $\rP^{(\infty)}_{\nu}$, and, conditionally on
   $Y$, let  $\cn=\sum_{j\in J} \delta_{(s_j,X^{j})}$ be  a Poisson point
   measure   with  intensity:
\[
2  \ind_{\R^+}(s)   ds   \;  \alpha(Y_s)
   \N_{Y_s}[dX]. 
\]
Consider the process $X^{(\infty )}=(X^{(\infty )}_t, t\geq
  0)$, which is  defined  for all $t\geq 0$ by:
\[
X^{(\infty )}_t= \sum_{j\in J, \, s_j<t} X_{t-s_j }^{j}. 
\]
\begin{itemize}
\item[(i)]  Let $X'$ independent of $X^{(\infty )}$ and distributed according
  to $\P_{\nu}$. Then, we write $ \P_\nu^{(\infty )}$ for the
  distribution of $X'+X^{(\infty )}$. 
   \item[(ii)] If $\nu$ is the Dirac mass at $x$, we write $\N^{(\infty )}_x$
     for the distribution of $X^{(\infty )}$.
\end{itemize}
\end{defi}

As a  consequence of Theorem \ref{theo-convN}, we get  the convergence of  $\P_\nu^{(h)}$. We shall 
write $\P_x^{(h)}$ when $\nu$ is the Dirac mass at $x$.

\begin{cor}
   \label{cor:cvNh}
Under $(H1)$-$(H4)$, we have that, for all $t\geq 0$: 
\begin{itemize}
\item[(i)]  The distribution $\N_{x}^{(h)}$ converges weakly to
  $\N_{x}^{(\infty)}$ on $(\Omega, \F_t)$.
\item[(ii)] The distribution $\P_{x}^{(h)}$ converges weakly to
  $\P_{x}^{(\infty)}$ on $(\Omega, \F_t)$.

\item[(iii)] Let $\nu\in \cm_f(E)$.  If furthermore $(H5)_\nu$ holds, then
  the     distribution     $\P_{\nu}^{(h)}$     converges    weakly     to
  $\P_{\nu}^{(\infty)}$ on $(\Omega, \F_t)$.
\end{itemize}
\end{cor}

\begin{proof}
Point (i) is a direct consequence of    Theorem \ref{theo-convN},
Definition \ref{defi:Q} and Proposition \ref{proplinkgenealogy}. 

Point (ii) is a direct consequence of point (i), Corollary
\ref{corwilliams} and the weak convergence of  $\P_{x}^{(\leq h)}$ to
$\P_x$ as $h$ goes to infinity. 

According to  Corollary \ref{corwilliams}, under  $\P_\nu^{(h)}$, $X$ is
distributed  according  to $X'+X^{(h)}$  where  $X'$  and $X^{(h)}$  are
independent, $X'$ is distributed  according to $\P_{\nu}^{(\leq h)}$ and
$X^{(h)}$ is distributed according to
\[
\int_E \nu(dx)\; \frac{\partial_{h} v_{h} (x)}{\nu(\partial_{h} v_{h})}
\; \N^{(h)}_x[dX]. 
\]
Assumption $(H5)_\nu$ implies this distribution converges weakly to:
\[
\int_E \nu(dx)\; \rho(x) 
\; \N^{(\infty )}_x[dX]
\]
(because of the convergence of the densities in $\rL^1(\nu)$) on $(\Omega,
\F_t)$ as $h$ goes to infinity.  This and the weak convergence of
$\P_{\nu}^{(\leq h)}$ to 
$\P_\nu$ as $h$ goes to infinity gives point (iii). 
\end{proof}

\begin{proof}[Proof of Lemma \ref{lem:cvPoisson}]
Similarly to Lemma \ref{lem:GirsanovPoisson} (formally take $f=-\infty
\ind_{A^c}$), we have:
\[
\E\left[F(R)\ind_{\{R(A^c)=0\}} \right]=\E\left[F(\tilde R)\right]\expp{-
  \nu(A^c)}.
\]
We deduce that:
\begin{align*}
\val{\E[F(R)]-\E[F(\tilde R)]}
&=\val{\E[F(R)]-\E[F(R)\ind_{\{R(A^c)=0\}}]\expp{  \nu(A^c)}}\\
&\leq \val{\E[F(R)]-\E[F(R)\ind_{\{R(A^c)=0\}}]}
+ \val{\E[F(R)\ind_{\{R(A^c)=0\}}] (1-\expp{  \nu(A^c)})}\\
&\leq  \norm{F}_\infty (1- \P(R(A^c)=0))
+ \norm{F}_\infty \P(R(A^c)=0) (\expp{  \nu(A^c)}-1)\\
&  = 2\norm{F}_\infty  (1-\expp{-\nu(A^c)})\\
&\leq 2\norm{F}_\infty \nu(A^c).
\end{align*}   
This gives the result.
\end{proof}

\subsection{Backward from the extinction time}
\label{backward}

We  shall work  in this  section with  the space  $D^{-}  = D(\R^{-},E)$
equipped   with   the  Skorokhod   topology.   We   also  consider   the
$\sigma$-fields $\D_{I}  = \sigma(Y_r, r\in I)$  for $I$ an interval on
$(-\infty , 0]$.

Let us denote by $\theta$ the translation operator, which maps any process $R$ to the shifted process $\theta_{h}(R)$ defined by:
\[
\theta_{h}(R)_{\cdot}= R_{\cdot+h} \ .
\] 
The process  $R$ may be  a path,  a killed path  or a point  measure, in
which  case   we  set,   for  $R=\sum_{j  \in   J}  \delta_{(s_j,x_j)}$,
$\theta_{h}(R)  = \sum_{j  \in J}  \delta_{(h+s_j,x_j)}$.   We also
denote $\rP^{(-h)}$ the push  forward probability measure of $\rP^{(h)}$
by $\theta_{h}$, defined on $\D_{[-h,0]}$ by:
\begin{align}
\label{defP-h}
\rP^{(-h)}(Y  \in  \bullet)  =  \rP^{(h)} (\theta_{h}(Y)  \in  \bullet)=
\rP^{(h)} ((Y_{h+s}, s \in [-h,0]) \in \bullet).
\end{align}

We introduce the following assumptions.

\medskip
\noindent 
$(H6)$ \textbf{There exists a probability measure on
  $(D^{-},\D_{(-\infty ,0]})$ denoted 
  $\rP^{(-\infty)}$ such that for all $x \in E$, $t\geq 0$, and $f$
  bounded and $\D_{[-t,0]}$ measurable:}  
\begin{align*}
  \rE_{x}^{(-h)}    \big[f(Y_{[-t,0]})\big]   \;    \xrightarrow[h   \to
  +\infty]{} \; \rE^{(-\infty)} \big[f(Y_{[-t,0]})\big].
\end{align*}
\medskip

\noindent
\noindent
$(H7)$ \textbf{For all $t>0$, there exists a non negative function $g$ such that  for all $x \in E$, for all $ h>0$:} 
\begin{align*}
v_{h}(x)-v_{h+t}(x) \leq g(h) \quad\text{and}\quad  \int_{1}^{\infty} dr \; g(r) < \infty.
\end{align*}
\medskip

Note that  the probability measure $\rP^{(-\infty)}$ in  $(H6)$ does not
depend on the starting point $x$.

\medskip

We can now state the result on the convergence of the superprocess
backward from the extinction time.

\begin{theo}
\label{convtop}
Under $(H1)$-$(H4)$ and $(H6)$. 
\begin{itemize}
\item[(i)]     The     distribution     of    the     triplet     $\big(
  \theta_{h}(W_{\Tm})_{[-t,0]},         \theta_{h}(R_{g}^{\Tm})_{[-t,0]},
  \theta_{h}(R_{d}^{\Tm})_{[-t,0]}\big)$  under  $\NN_x^{(h)}$ converges
  weakly   to   the    distribution   of   the   triplet   $(Y_{[-t,0]},
  R_{[-t,0]}^{W,g},  R_{[-t,0]}^{W,d})$  where   $Y$  has  distribution $\rP^{(-
    \infty)}$ and conditionally on  $Y$, $R^{W,g}$ and $R^{W,d}$ are two
  independent Poisson point measures with intensity:
\[
\ind_{\{s<0\}}  
\alpha(Y_s) \: ds \; \ind_{\{\Hm(W) < -s\}} \;  \NN_{Y_s}[dW].
\] 
We even have the slightly
stronger result. For any bounded measurable function $F$, we have:
\begin{multline}
\label{convNhtop}
\NN_x^{(h)}\bigg[F\big(
  \theta_{h}(W_{\Tm})_{[-t,0]},         \theta_{h}(R_{g}^{\Tm})_{[-t,0]},
  \theta_{h}(R_{d}^{\Tm})_{[-t,0]}\big)\bigg] \\
\; \xrightarrow[h \to +\infty]{}  \; 
 \rE^{(-\infty)} \bigg[
F\big(Y_{[-t,0]}, R^{W,g}_{[-t,0]},  R^{W,d}_{[-t,0]} \big) \bigg].
\end{multline}
\item[(ii)]   If    furthermore   $(H7)$   holds,    then   the   process
  $\theta_{h}(X)_{[-t,0]}=(X_{h+s},  s\in [-t,0])$  under $\N_{x}^{(h)}$ weakly  converges towards  $X^{(-\infty  )}_{[-t,0]}$, where  for
  $s\leq 0$:
\[
X^{(-\infty )}_s=\sum_{j\in J,  \; s_j<s} X^j_{s-s_j},
\]
and conditionally on $Y$ with distribution $\rP^{(-\infty)}$, $\sum_{j\in J}
\delta_{(s_j,X^{j})}$ is a Poisson point measure with
intensity:
 \[ 
2 \;  \ind_{\{s<0\}} \alpha(Y_s) \: ds \; \ind_{\{\Hm(X) < -s\}} \;
\N_{Y_s}[ dX].
\] 
\end{itemize}
\end{theo}

\begin{rem}
We provide in Lemmas \ref{H6H7multi} and \ref{lem123} sufficient conditions for $(H6)$ and $(H7)$ to hold in the case of the multitype
Feller diffusion and the superdiffusion. These conditions are stated in term of the generalized eigenvalue
$\lambda_{0}$ defined in \reff{eq:eigenvalue} and its associated eigenfunction.
\end{rem}
% The assumptions $(H6)$ and $(H7)$ will be checked in Section \ref{sec:twoexamples}  for the multitype
% Feller diffusion and the superdiffusion, see Lemmas  in which
%  is positive, and the spatial motion $\rP^{\phi_{0}}$ (recall the definition \reff{ACequation}) satisfies an ergodic assumption. 
% be  the probability measure on $(D, \D)$
% defined by \reff{ACequation} with $g$ replaced by $\phi_{0}$:

\begin{proof}
  Let  $0<t<h$.  We  use  notations  from  Theorems  \ref{theowilliams},
  \ref{theo-convN} and  \ref{convtop}.  Let $F$ be  a bounded measurable
  function  on  $\W^{-}  \times  (\R^{-}  \times  \bar{\Omega})^2$  with
  $\W^{-}$ the set of killed paths indexed by negative times. We want to
  control  $\delta_{h}$ defined by:
\begin{multline*}
\delta_{h}=
\NN_{x}^{(h)}\bigg[F\big(\theta_{h}(W_{\Tm})_{[-t,0]},
\theta_{h}(R_\Tm^g)_{[-t,0]},   
 \theta_{h}(R_\Tm^d )_{[-t,0]}\big)\bigg] \\
-  \rE^{(-\infty)} \bigg[
F\big(Y_{[-t,0]}, R^{W,g}_{[-t,0]},  R^{W,d}_{[-t,0]} \big) \bigg].
\end{multline*}
We set:
\[
\Upsilon(y_{[-t,0]})
=\rE^{(-\infty )} \bigg[
F\big(y_{[-t,0]}, R^{W,g}_{[-t,0]},  R^{W,d}_{[-t,0]} \big)
\bigg|Y=y \bigg].
\]
We deduce from  Williams' decomposition, Theorem \ref{theowilliams},
and the definition of $R^{W,g}$ and $R^{W,d}$,  that:
\[
\NN_{x}^{(h)}\bigg[F\big(\theta_{h}(W_{\Tm})_{[-t,0]},
\theta_{h}(R_\Tm^g)_{[-t,0]},  
 \theta_{h}(R_\Tm^d )_{[-t,0]}\big)\bigg]
= \rE_{x}^{(-h)} \big[  \Upsilon(Y_{[-t,0]})\big].
\]
We thus can rewrite $\delta_{h}$ as:
\begin{align*}
\delta_{h}
&=  \rE_{x}^{(-h)} \big[  \Upsilon(Y_{[-t,0]})\big]- \rE^{(-\infty)} \big[  \Upsilon(Y_{[-t,0]})\big].  
\end{align*}
The   function  $\Upsilon$   being  bounded   by   $\|F\|_{\infty}$  and
measurable,   we    may   conclude   under    assumption   $(H6)$   that
$\lim_{h\rightarrow+\infty } \delta_{h}=0$. This proves point (i).

\medskip 
We now  prove point  (ii). Let $t>0$ and $\varepsilon>0$ be fixed. 
Let  $F$ be a bounded measurable function
on  the  space  of  continuous measure-valued  applications  indexed  by
negative   times.   For a
point measure on $\R^{-}\times \bar{\Omega}$,  $M=\sum_{i \in
  \I} \delta_{(s_i,W_i)}$,  we set:
\begin{align*}
\tilde{F}(M)=F\left(\big(\sum_{i \in \I} \theta_{s_i}
  (X(W_i)) \big)_{[-t,0]} \right).
\end{align*}
For $h>t$, we want a control of $\bar \delta_h$ defined by:
\[
\bar{\delta}_{h} =
\N_{x}^{(h)}\bigg[F\big(\theta_{h}(X)_{[-t,0]}\big)\bigg]  -
\rE^{(-\infty)} \bigg[ 
\tilde{F} \big( R^{W,g}+  R^{W,d} \big) \bigg] . 
\]
By Corollary \ref{cor:Nh}, we have:
\[
\N_{x}^{(h)}\bigg[F\big(\theta_{h}(X)_{[-t,0]}\big)\bigg] 
 = \NN_{x}^{(h)} \bigg[\tilde{F} \big(\theta_{h}(R_\Tm^g+R_\Tm^d ) \big)\bigg].
\]
Thus, we get:
\begin{equation}
\label{eq:bardelta}
\bar{\delta}_{h}  = \NN_{x}^{(h)} \bigg[\tilde{F}
\big(\theta_{h}(R_\Tm^g+ R_\Tm^d ) \big)\bigg]  -  \rE^{(-\infty)} \bigg[
\tilde{F} \big(R^{W,g}+  R^{W,d} \big) \bigg].
\end{equation}
For $a>s$ fixed, we introduce $\bar{\delta}^{a}_{h}$, for $h>a$, defined
by:
\begin{equation}
\label{eq:bardeltaa}
\bar{\delta}^{a}_{h}        =        \NN_{x}^{(h)}       \bigg[\tilde{F}
\big(\theta_{h}(R_\Tm^g+ R_\Tm^d)_{[-a,0]}
\big)\bigg]  -  \rE^{(-\infty)}  \bigg[  \tilde{F}
\big((R^{W,g}+R^{W,d}) _{[-a,0]} \big) \bigg].
\end{equation}
Notice the restriction of the point measures to $[-a,0]$.
Point (i) directly yields that $\lim_{h\rightarrow+\infty }
\bar{\delta}^{a}_{h}=0 $. Thus, there exists $h_a>0$ such that for all
$h \geq h_a$, 
\begin{align*}
\bar{\delta}^{a}_{h} \leq \varepsilon/2.
\end{align*}
We now consider the difference $\bar{\delta}_{h}-
\bar{\delta}^{a}_{h}$. We associate to  the point measures $M$
introduced above the most recent common ancestor of the 
population alive at time $-t$:
\[
 A(M) =\sup\{ s>0 ; \sum_{i\in \I} \ind_{\{s_i<-s\}}
   \ind_{\{\Hm(W_i)>-t-s_i\}} 
 \neq 0 \}. 
\]
Let us observe that:
\begin{equation}
\label{eq:A1}
\NN_{x}^{(h)} \mbox{ a.s., }  \  \tilde{F} \big(\theta_{h}(R_\Tm^g + R_\Tm^d )\big) \ind_{\{A \leq a \}}
= \tilde{F} \big(\theta_{h}(R_\Tm^g+ R_\Tm^d )_{[-a,0]}\big) \ind_{\{A
  \leq a \}} , 
\end{equation}
with $A=A(\theta_{h}(R_\Tm^g+R_\Tm^d )_{[-h,0]})$ in the left and in the
right hand side. Similarly, we have:
\begin{equation}
\label{eq:A2}
\rP^{(-\infty)} \mbox{ a.s., } \ \tilde{F} \big(R^{W,g}+  R^{W, d}
\big)  \ind_{\{A \leq a \}}  =  
\tilde{F} \big( (R^{W,g}+R^{W,d})_{[-a,0]} \big)  \ind_{\{A
  \leq a \}} , 
\end{equation}
with $A=A\big(R^{W,g}+  R^{W,d} \big)$ in the left and in the right hand side.
We thus deduce the following bound on $\bar{\delta}_{h}- \bar{\delta}^{a}_{h}$:
\begin{align*}
|\bar{\delta}_{h}- \bar{\delta}^{a}_{h}| & \leq  2 \|F\|_{\infty} \bigg[\NN_{x}^{(h)}\big[A > a\big] +
  \rP^{(-\infty)}\big[A > a\big] \bigg]\\ 
& =  2 \|F\|_{\infty} \bigg[ \rE^{(-h)}_{x}\big[1-\expp{-\int_{a}^{h} \!\!dr
  \ 2\alpha(v_{r-t}-v_{r})(Y_{-r})}] +
\rE^{(-\infty)}\big[1-\expp{-\int_{a}^{\infty} \!\!dr \
  2\alpha(v_{r-t}-v_{r})(Y_{-r})}] \bigg]  \\ 
& \leq 8 \|F\|_{\infty} \|\alpha\|_\infty \int_{a-t}^{\infty} dr \ g(r),
\end{align*}
where we used  \reff{eq:bardelta}, \reff{eq:bardeltaa}, \reff{eq:A1} and
\reff{eq:A2} for  the first  inequality, the definition  of $A$  for the
first equality, as well as $(H7)$ and the fact that $1-\expp{-x} \leq x$
if  $x \geq 0$ for the last inequality.  From $(H7)$, we can choose $a$
large enough  such that: $|\bar{\delta}_{h}-  \bar{\delta}^{a}_{h}| \leq
\varepsilon/2$.    We  deduce  that   for  all   $h\geq  \max(a,h_{a})$:
$|\bar{\delta}_{h}|  \leq   |\bar{\delta}_{h}-  \bar{\delta}^{a}_{h}|  +
|\bar{\delta}^{a}_{h}| \leq \varepsilon$.  This proves point (ii).
\end{proof}

\section{The assumptions $(H4)$, $(H5)_\nu$ and  $(H6)$} 
\label{sec:H}

\textbf{We  assume in  all  this section  that  $\rP$ is  the  distribution of  a
  diffusion in $\R^K$ for $K$ integer or the  law  of a  finite  state space  Markov Chain,}  see
Section \ref{sec:twoexamples} and the references therein. In particular,
the   generalized  eigenvalue   $\lambda_{0}$   of  $(\beta - \mL)$ (see
\reff{defeigenvalue1}  or \reff{defeigenvalue2}) is  known to  exist. We
will denote by $\phi_{0}$ the associated right eigenvector. We shall consider
the assumption:

\medskip
 \noindent
$(H8)$ \textbf{There exist two positive constants $C_1$ and $C_2$ such
  that $\forall x \in E, \ C_1 \leq  \phi_{0}(x) \leq C_2$; and $\phi_{0}\in \D(\mL)$.}
\medskip

Under $(H8)$, let $\rP^{\phi_{0}}_x$ be  the probability measure on $(D, \D)$
defined by \reff{ACequation} with $g$ replaced by $\phi_{0}$:
\begin{equation}
\label{defPphi}
\forall t \geq 0, \quad \frac{d \rP_{x \ | \D_t}^{\phi_{0}}} {d \rP_{x \ |
    \D_t}} = \frac{\phi_{0}(Y_t)}{\phi_{0}(Y_0)} \expp{- \int_{0}^{t} ds \
  (\beta(Y_s)- \lambda_{0})}. 
\end{equation}

We shall also consider
the assumption:

\medskip
 \noindent
$(H9)$ \textbf{The probability measure $\rP^{\phi_{0}}$ admits a
  stationary measure $\pi$, and we have:}
\begin{equation}
\label{H6''}
\sup_{f \in b\Ec, \|f \|_\infty  \leq 1} |\rE_{x}^{\phi_{0}}[f(Y_t)]-\pi(f)| \;
\xrightarrow[t \to +\infty]{}  \; 0.
\end{equation}

Notice the  two hypotheses $(H8)$ and $(H9)$ hold for the examples of
Section \ref{sec:twoexamples}, see Lemmas \ref{rem:H9-Feller} and
\ref{rem:H9-diffu}. 

Let us mention at this point that we will check that
$\rP^{\phi_{0}}_{x}=\rP^{(\infty)}_{x}$ with $\rP_{x}^{(\infty)}$ defined by \reff{defPinfty}, see Proposition  \ref{prop:H4+Pf=Pinf}.

\subsection{Proof of $(H4)$-$(H6)$}
\label{sec:H4+H6}

Notice $(H9)$ implies that the probability measure $\rP^{\phi_{0}}_{\pi}$ admits a
stationary version on $D(\R,E)$, which we still denote by
$\rP_{\pi}^{\phi_{0}}$. 

We  introduce  a  specific   $h$-transform  of  the  superprocess.  From
Proposition \ref{prophtransform}  and the definition  of the generalized
eigenvalue \reff{defeigenvalue1} and \reff{defeigenvalue2}, we have that
the $h$-transform given by  Definition \ref{htransform} with $g=\phi_{0}$ of
the        $(\mL,\beta,\alpha)$        superprocess        is        the
$(\mL^{\phi_{0}},\lambda_{0},\alpha \phi_{0})$  superprocess.  We define
$v^{\phi_{0}}$ for all $t>0$ and $x \in E$ by:
\begin{equation}
 \label{eq:vphi}
v^{\phi_{0}}_t(x)=\N_{x}^{(\mL^{\phi_{0}},\lambda_{0},\alpha \phi_{0})}[\Hm>t].
\end{equation}
Observe that, as in \reff{eqnormalisation}, the following normalization
holds between $v^{\phi_{0}}$ and $v$:
\begin{equation}
   \label{eq:normalisation-vf}
v^{\phi_{0}}_t(x)=\frac{v_t(x)}{\phi_{0}(x)}\cdot
 \end{equation} 
Our first task is to give precise bounds on the decay of $v_{t}^{\phi_{0}}$
as $t$ goes to $\infty$. 

We first offer bounds
for the case $\lambda_{0}=0$ in Lemma \ref{boundcritical}, relying on a
coupling argument. This in turn gives sufficient condition under which $(H1)$ holds in Lemma \ref{nonnegativelambdaH1}  
We then give  Feynman-Kac representation
formulae, Lemma \ref{FK}, which yield  exponential bounds in the case
$\lambda_{0}>0$, see Lemma \ref{bound}. 
We finally strengthen in Lemma \ref{FK+} the bound of Lemma \ref{bound} by proving the exponential behavior of $v_{t}^{\phi_{0}}$ in the case $\lambda_{0}>0$.
The proofs of Lemmas \ref{boundcritical}, \ref{nonnegativelambdaH1}, \ref{FK}, \ref{bound} and
\ref{FK+} are given in Section \ref{sec:preuve-annexe}.

We first give  a bound in the case $\lambda_{0}=0$.
% We set $\underline{\alpha \phi_{0}} = \inf_{x \in E} \{ \alpha(x)
% \phi_{0}(x) \}$, and $\overline{\alpha \phi_{0}} = \sup_{x \in E} \{ \alpha(x)
% \phi_{0}(x) \}=\norm{\alpha\phi_{0}}_\infty $.  
\begin{lem}
\label{boundcritical}
Assume $\lambda_{0}=0$, $(H2)$ and $(H8)$. Then for all $t >0$:
\begin{align*}
\alpha \phi_{0}(x) \frac{1}{\left\| \alpha \phi_{0} \right\|^2 _{\infty}} \leq  t\; v^{\phi_{0}}_{t}(x) \leq  \alpha \phi_0 (x) \left\| \frac{1}{ \alpha \phi_{0} } \right\|^2 _{\infty}.
\end{align*} 
\end{lem}
A coupling argument then implies that $(H1)$ holds:
\begin{lem}
\label{nonnegativelambdaH1}
Assume $\lambda_{0} \geq 0$, $(H2)$ and $(H8)$. Then $(H1)$ holds.
\end{lem}
We give a Feynman-Kac's formula for $v^\phi_{0}$ and $\partial
v^\phi_{0}$. 
\begin{lem} 
\label{FK}
Assume $\lambda_{0} \geq 0$, $(H2)$-$(H3)$ and $(H8)$. Let $\varepsilon >0$. We have:
\begin{align}
\label{eq:FKv}
v_{h+\varepsilon}^{\phi_{0}}(x) 
&=  \expp{-\lambda_{0} h}  \rE^{\phi_{0}}_{x}
\bigg[ \expp{- \int_{0}^{h} ds \ \alpha(Y_s) \: \phi_{0}(Y_s) \:
  v^{\phi_{0}}_{h+\varepsilon-s}(Y_s)} v_{\varepsilon}^{\phi_{0}}(Y_h)\bigg], \\ 
\label{eq:FKdv}
\partial_{h} v_{h+\varepsilon}^{\phi_{0}}(x) 
&= \expp{-\lambda_{0} h} \rE^{\phi_{0}}_{x} \bigg[ \expp{- 2\int_{0}^{h} ds \
   \alpha(Y_s) \: \phi_{0}(Y_s) \: v^{\phi_{0}}_{h+\varepsilon-s}(Y_s)} \partial_{h}
v_{\varepsilon}^{\phi_{0}}(Y_h)\bigg]. 
\end{align}
\end{lem}

We give exponential bounds for $v^\phi_{0}$ and $\partial_t
v^\phi_{0}$ in the subcritical case. 
\begin{lem}
\label{bound}
Assume $\lambda_{0}>0$, $(H2)$-$(H3)$ and $(H8)$. Fix $t_{0}>0$. There exists $C_3$ and $C_4$ two positive constants such that, for all $x \in E$, $t>t_{0}$:
\begin{align}
\label{eq:encadr-v}
C_3    & \leq v^{\phi_{0}}_{t}(x) \expp{\lambda_{0} t} \leq C_4. 
\end{align}
There exists $C_5$ and $C_6$ two positive constants such that, for all $x \in E$, $t>t_{0}$:
\begin{align}
\label{eq:encadr-dv}
C_5   &   \leq |\partial_{t} v_{t}^{\phi_{0}}(x)| \expp{\lambda_{0} t} \leq C_6. 
\end{align}
\end{lem}

As a direct consequence of \reff{eq:encadr-v}, we get the following Lemma. 
\begin{lem}
\label{lem:H7}
Assume $\lambda_{0}>0$, $(H2)$-$(H3)$ and  $(H8)$. Then $(H7)$ holds.
\end{lem}
In what follows,  the notation $o_h(1)$ refers to  any function $F_h$
such  that $\lim_{h\rightarrow+\infty }  \norm{F_h}_\infty =0$.   We now
improve on Lemma \ref{bound}, by using the ergodic formula \reff{H6''}. 
\begin{lem}
\label{FK+}
Assume $\lambda_0\geq 0$, $(H2)$-$(H3)$ and $(H8)$-$(H9)$ hold.  Then for all $\varepsilon >0$,  we have:
\begin{align}
\label{eq:dev-dv1-enh}
 \partial_{t} v^{\phi_{0}}_{h+\varepsilon}(x) \; \expp{\lambda_{0} h}
= \rE_{\pi}^{\phi_{0}} \big[\expp{-2\int_{0}^{h} ds \  \alpha\: \phi_{0} \:
  v^{\phi_{0}}_{s+\varepsilon}(Y_{-s})} \partial_{t}
v^{\phi_{0}}_{\varepsilon}(Y_0) \big] (1+o_h(1)). 
\end{align}
In addition, for $\lambda_{0} >0$, we have that: 
\[
\rE_{\pi}^{\phi_{0}} \big[\expp{-2\int_{0}^{\infty} ds \  \alpha\: \phi_{0} \:
  v^{\phi_{0}}_{s+\varepsilon}(Y_{-s})} \partial_{t}
v^{\phi_{0}}_{\varepsilon}(Y_0) \big] 
\]
is finite (notice the integration is up to $+\infty $) and:
\begin{align}
\label{eq:dev-dv2-enh}
\partial_{t} v^{\phi_{0}}_{h+\varepsilon}(x) \; \expp{\lambda_{0} h} 
= \rE_{\pi}^{\phi_{0}} \big[\expp{-2 \int_{0}^{\infty} ds \  \alpha\: \phi_{0} \:
  v^{\phi_{0}}_{s+\varepsilon}(Y_{-s})} \partial_{t}
v^{\phi_{0}}_{\varepsilon}(Y_0) \big] + o_h(1). 
\end{align}
\end{lem}

\medskip

Our next goal is to prove $(H4)$ from $(H8)$-$(H9)$, see Proposition \ref{prop:H4+Pf=Pinf}. 

Fix $x \in E$. We observe from \reff{defPh2} and \reff{defPphi} that
$\rP^{(h)}_{x}$ is absolutely continuous with respect to
$\rP^{\phi_{0}}_{x}$ on $\D_{[0,t]}$ for $0 \leq t <h$. We define
$M_{t}^{(h),\phi_{0}}$  the corresponding Radon-Nikodym derivative:
\[
M_{t}^{(h),\phi_{0}} = \frac{d \rP_{x \ | \D_{[0,t]}}^{(h)}} {d
  \rP_{x \ | \D_{[0,t]}}^{\phi_{0}}} \cdot
\]
Using  \reff{defPh2}, \reff{defPphi} and 
the normalization $ v(x) =  v^{\phi_{0}}(x)\;
\phi_{0}(x)$, we get:  
\begin{align}
\label{defPphih+} 
M_{t}^{(h),\phi_{0}} 
&= \frac{\partial_{t} v_{h-t}(Y_{t})\expp{-\lambda_{0} t} }{\partial_{t}
  v_{h}(Y_{0})}  \frac{\phi_{0}(Y_0)}{\phi_{0}(Y_t)} \expp{- 2\int_{0}^{t} ds \
  \alpha(Y_s) \; v_{h-s}(Y_s)} \\
&=\frac{\partial_{t} v^{\phi_{0}}_{h-t}(Y_{t})\expp{-\lambda_{0} t}
}{\partial_{t} v^{\phi_{0}}_{h}(Y_{0})}  \expp{- 2\int_{0}^{t} ds \ 
  \alpha(Y_s) \; \phi_{0}(Y_s) \; v^{\phi_{0}}_{h-s}(Y_s)}. \nonumber
\end{align}

We have the following result on the convergence of $M_{t}^{(h),\phi_{0}} $. 
\begin{lem}
\label{lem:H4}
Assume $(H2)$-$(H3)$ and $(H8)$-$(H9)$.
For $\lambda_{0} \geq 0$, we have:
\begin{align*}
M_{t}^{(h),\phi_{0}}  
\; \xrightarrow[h \to +\infty]{}  \; 
1 \quad \text{$\rP_{x}^{\phi_{0}}$-a.s. and in $\rL^{1}(\rP_{x}^{\phi_{0}})$},
\end{align*}
and for $\lambda_{0}>0$, we have:
\begin{align*}
M_{h/2}^{(h),\phi_{0}}  
\; \xrightarrow[h \to +\infty]{}  \; 
1 \quad \text{$\rP_{x}^{\phi_{0}}$-a.s. and in $\rL^{1}(\rP_{x}^{\phi_{0}})$}.
\end{align*}
\end{lem}

\begin{proof}
We compute:
\begin{align*}
M_{t}^{(h),\phi_{0}} 
&= \frac{\partial_{t} v^{\phi_{0}}_{h-t}(Y_{t})\expp{\lambda_{0} (h-t)}
}{\partial_{t} v_{h}^{\phi_{0}}(Y_{0}) \expp{\lambda_{0} h}}  \expp{-
2  \int_{0}^{t} ds \  \alpha(Y_s) \phi_{0}(Y_s)  v^{\phi_{0}}_{h-s}(Y_s)}\\  
&=\frac{\rE^{\phi_{0}}_{\pi} \big[   \expp{- 2\int_{-(h-t-\varepsilon)}^{0} ds \ 
    \alpha(Y_s) \: \phi_{0}(Y_s) \: v^{\phi_{0}}_{\varepsilon-s}( Y_s)} \partial_{h}
  v_{\varepsilon}^{\phi_{0}}({Y}_{0}) \big]\big(1+o_h(1)
  \big)}{\rE^{\phi_{0}}_{\pi} \big[   \expp{- 2\int_{-h-\varepsilon}^{0} ds \ 
    \alpha(Y_s) \: \phi_{0}(Y_s) \: v^{\phi_{0}}_{\varepsilon-s}( Y_s)} \partial_{h}
  v_{\varepsilon}^{\phi_{0}}({Y}_{0}) \big]\big(1+o_h(1)\big)}
\big(1+o_h(1)\big)\\ 
&=\frac{\rE^{\phi_{0}}_{\pi} \big[   \expp{- 2\int_{-h-\varepsilon}^{0} ds \  \alpha(Y_s)
    \: \phi_{0}(Y_s) \: v^{\phi_{0}}_{\varepsilon-s}( Y_s)} \partial_{h}
  v_{\varepsilon}^{\phi_{0}}({Y}_{0}) \big]}
{\rE^{\phi_{0}}_{\pi} \big[   \expp{- 2\int_{-h-\varepsilon}^{0} ds \ 
    \alpha(Y_s) \: \phi_{0}(Y_s) \: v^{\phi_{0}}_{\varepsilon-s}( Y_s)} \partial_{h}
  v_{\varepsilon}^{\phi_{0}}({Y}_{0}) \big]}\big(1+o_h(1)\big)\\
&=1+o_h(1),
\end{align*}
where we  used \reff{defPphih+} for the first  equality, \reff{eq:dev-dv1-enh}
twice  and  the  boundedness of  $\alpha$  and  $\phi_{0}$  as well  as  the
convergence of $v_{h}$ to $0$  for the second, and Lemma \ref{bound} (if
$\lambda_0>0$) or  Lemma \ref{boundcritical} (if  $\lambda_0=0$) for the
fourth.  Since $o_h(1)$ is bounded  and converges uniformly to $0$, we get
that   the   convergence  of   $M_{t}^{(h),\phi_{0}}$   towards  $1$   holds
$\rP^{\phi_{0}}_{x}$-a.s. and in $L^1(\rP^{\phi_{0}}_{x})$.

Similar   arguments   relying   on  \reff{eq:dev-dv2-enh}   instead   of
\reff{eq:dev-dv1-enh}  imply  that  $M_{h/2}^{(h),\phi_{0}}  =1+o_h(1)$  for
$\lambda_0>0$.  Since $o_h(1)$ is  bounded and converges uniformly to $0$,
we get  that the convergence  of $M_{h/2}^{(h),\phi_{0}}$ towards  $1$ holds
$\rP^{\phi_{0}}_{x}$-a.s. and in $L^1(\rP^{\phi_{0}}_{x})$.
\end{proof}

The previous Lemma enables us to conclude about $(H4)$. 

\begin{prop}
\label{prop:H4+Pf=Pinf}
Assume  $\lambda_{0}  \geq 0$, $(H2)$-$(H3)$  and $(H8)$-$(H9)$.  Then
$(H4)$ holds with $\rP_{x}^{(\infty)} = \rP_{x}^{\phi_{0}}$.
\end{prop}

\begin{proof}
Notice that:
\[
M_{t}^{(h)} = \frac{d \rP_{x \ | \D_{[0,t]}}^{(h)}} {d
  \rP_{x \ | \D_{[0,t]}}}
= M_{t}^{(h),\phi_{0}} \frac{d \rP_{x \ | \D_{[0,t]}}^{\phi_{0}}} {d
  \rP_{x \ | \D_{[0,t]}}} \cdot
\]
The convergence  $\lim_{h\rightarrow+\infty }M_{t}^{(h),\phi_{0}}  =1$
$\rP_{x}^{\phi_{0}}$-a.s. and in $\rL^{1}(\rP_{x}^{\phi_{0}})$ readily implies 
$(H4)$. Then, use \reff{defPinfty} to get  $\rP^{(\infty)} = \rP^{\phi_{0}}$.
\end{proof}

Notice that $(H5)_\nu$ is a direct consequence of Lemma \ref{FK+}.

\begin{cor}
   \label{cor:H5}
Assume  $\lambda_{0}  \geq 0$, $(H2)$-$(H3)$  and $(H8)$-$(H9)$.  Then
$(H5)_\nu$ holds  with
$\rho=\phi_{0}/\nu(\phi_{0})$. 
\end{cor}
\begin{proof}
We deduce from \reff{eq:normalisation-vf} and \reff{eq:dev-dv1-enh} that:
\[
\partial_{t} v_{h}(x) = f(h) \phi_{0}(x) \; (1+o_h(1)) \expp{-\lambda_{0} h}  ,
\]
for some positive function $f$ of $h$. Then we get:
\[
\frac{\partial_{h} v_{h}(x)}{
  \nu(\partial_{h} v_{h})}
= \frac{ \phi_{0}(x)} {\nu(\phi_{0})} (1+o_h(1)).
\] 
This gives $(H5)_\nu$, as $o_h(1)$ is bounded,  with
$\rho=\phi_{0}/\nu(\phi_{0})$. 
\end{proof}

\medskip

Our next goal is to prove $(H6)$ from $(H8)$-$(H9)$, see Proposition \ref{lem:H6}. 

Observe from \reff{eq:RNderivPhP}, \reff{defP-h} and \reff{defPphi} that
$\rP_{\pi}^{(-h)}$ is absolutely continuous with respect to
$\rP_{\pi}^{\phi_{0}}$ on $\D_{[-h,0]}$. We define
$L^{(-h)}$ the corresponding Radon-Nikodym derivative:
\begin{equation}
   \label{defPphih0}
L^{(-h)} 
= \frac{d \rP_{\pi \ | \D_{[-h,0]}}^{(-h)}} {d \rP_{\pi \ |
    \D_{[-h,0]}}^{\phi_{0}}} 
=\inv{\alpha(Y_0) \phi_{0}(Y_0)} \frac{\partial_h v^0_h\expp{\beta_0 h}}{\partial_{h}
  v^{\phi_{0}}_h(Y_{-h})\expp{\lambda_0 h} } \expp{-2\int_{-h}^0 
  (\alpha(Y_s) v_{-s}(Y_s) - v^0_{-s})\; ds}.
\end{equation}
The next Lemma insures the convergence of $L^{(-h)}$ to
a limit, say $L^{(-\infty)}$.

\begin{lem}
\label{lem:Lh-Linfty}
Assume $\lambda_{0}>0$, $(H2)$-$(H3)$ and $(H8)$-$(H9)$. We have:
\begin{align*}
L^{(-h)}
\; \xrightarrow[h \to +\infty]{}  \; 
L^{(-\infty )} 
 \quad \text{$\rP_{\pi}^{\phi_{0}}$-a.s. and in $\rL^{1}(\rP_{\pi}^{\phi_{0}})$}.
\end{align*}
\end{lem}

\begin{proof}
  Notice that $\lim_{h\rightarrow+\infty } \partial_h v^0_h\expp{\beta_0
    h}=  - \beta_0^2$.   We    also   deduce   from   \reff{eq:defSigma},
  \reff{eq:ineqSigma}   and  \reff{eq:encadr-v}   that   $\int_{-h}^0  
  (\alpha(Y_s) v_{-s}(Y_s)  - v^0_{-s})\; ds$ increases, as  $h$ goes to
  infinity;   to  $\int_{-\infty  }^0     (\alpha(Y_s)   v_{-s}(Y_s)  -
  v^0_{-s})\; ds$ which is finite.  For fixed $t>0$, we also deduce from
  \reff{eq:dev-dv2-enh} (with $h$ replaced by $h-t$ and $\varepsilon$ by
  $t$) that $\rP_{\pi}^{\phi_{0}}$ a.s.:
\[
\lim_{h\rightarrow+\infty } 
\partial_{t} v^{\phi_{0}}_{h}(Y_{-h}) \expp{\lambda_{0} h}
= \expp{\lambda_0 t} \rE_{\pi}^{\phi_{0}}\big[ \expp{- 2\int_{-\infty}^{-t}
  ds \ 
  \alpha(Y_s) \ \phi_{0}(Y_s) \, v^{\phi_{0}}_{-s}(Y_s)} \partial_{t}
v^{\phi_{0}}_{t}(Y_{-t})  \big].
\]
We deduce  from \reff{defPphih0} the  $\rP_{\pi}^{\phi_{0}}$ a.s. convergence
of    $(L^{(-h)},   h    >0)$   to    $L^{(-\infty)}$.     Notice   from
\reff{eq:encadr-dv} that,  for  fixed  $t$,  the
sequence $(L^{(-h)},  h >t)$ is bounded. Hence  the previous convergence
holds also in $\rL^{1}(\rP^{\phi_{0}}_{\pi})$.
\end{proof}

As $\rE^{\phi_{0}}_{\pi} \big[ L^{(-h)} \big] = 1$, we deduce that 
 $\rE^{\phi_{0}}_{\pi} \big[ L^{(-\infty)}\big] = 1$.  
We define the probability measure $\rP^{(-\infty),\phi_{0}}_{\pi}$ on
$(D^-,\D_{(-\infty ,0]})$ by its Radon Nikodym derivative: 
\begin{align}
\label{defPphiinfty}
\frac{d \rP_{\pi \ | \D_{(-\infty ,0]}}^{(-\infty), \phi_{0}}}{d \rP_{\pi \
    | \D_{(-\infty ,0]}}^{\phi_{0}}} = 
L^{(-\infty)}.
\end{align}

\begin{rem}
   \label{rem:Linfini}
Assume $\lambda_0 >0$, $(H2)$-$(H3)$ and $(H8)$-$(H9)$.
Define for $h>t>0$:
\begin{align*}
L^{(-h)}_ {-t} \!\! &= \rE^{\phi_{0}}_\pi [L^{(-h)}|\D_{(-\infty ,-t]}]
= \frac{d \rP_{\pi \ | \D_{[-h,-t]}}^{(-h)}} {d \rP_{\pi \ |
    \D_{[-h,-t]}}^{\phi_{0}}} \\
L^{(-\infty )}_ {-t}\!\! &=\rE^{\phi_{0}}_\pi [L^{(-\infty )}|\D_{(-\infty ,-t]}]
= \frac{d \rP_{\pi \ | \D_{(-\infty ,-t]}}^{(-\infty ),\phi_{0}}} {d \rP_{\pi \ |
    \D_{(-\infty ,-t]}}^{\phi_{0}}}\cdot
\end{align*}
Using \reff{eq:v0} and Lemma \ref{FK}, we get:
\begin{align*}
L^{(-h)}_{-t} 
&= \frac{\partial_{t} v_{t}(Y_{-t})}{\partial_{t} v_{h}(Y_{-h})}
\frac{\phi_{0}(Y_{-h})}{\phi_{0}(Y_{-t})} \expp{-\lambda_{0} (h-t)}\expp{-
  2\int_{-h}^{-t} ds \  \alpha(Y_s) \, v_{-s}(Y_s)} \\ 
&=\frac{ \expp{- 2\int_{-h}^{-t} ds \  \alpha(Y_s) \ \phi_{0}(Y_s) \,
    v^{\phi_{0}}_{-s}(Y_s)} \partial_{t} v^{\phi_{0}}_{t}(Y_{-t}) }
{\rE_{Y_{-h}}^{\phi_{0}}\big[ \expp{- 2\int_{-h}^{-t} ds \  \alpha(Y_s) \
    \phi_{0}(Y_s) \, v^{\phi_{0}}_{-s}(Y_s)} \partial_{t} v^{\phi_{0}}_{t}(Y_{-t})
  \big]}\cdot 
\end{align*}
Using Lemma \ref{FK+} and convergence of $(L^{(-h)}_{-t} , h>t)$ to $L^{(-\infty)}_{-t}$, which
is a consequence of Lemma \ref{lem:Lh-Linfty}, we also get that for $t> 0$:
\begin{align*}
L^{(-\infty)}_{-t} =   \frac{ \expp{- 2\int_{-\infty}^{-t} ds \ 
    \alpha(Y_s) \ \phi_{0}(Y_s) \, v^{\phi_{0}}_{-s}(Y_s)} \partial_{t}
  v^{\phi_{0}}_{t}(Y_{-t}) } {\rE_{\pi}^{\phi_{0}}\big[ \expp{-
    2\int_{-\infty}^{-t} ds \ \alpha(Y_s) \ \phi_{0}(Y_s) \,
    v^{\phi_{0}}_{-s}(Y_s)} \partial_{t} v^{\phi_{0}}_{t}(Y_{-t})  \big]}\cdot
\end{align*}
Those  formulas are  more  self-contained than  \reff{defPphih0} and  the
definition of  $L^{(-\infty )}$  as a  limit, but they  only  hold  for
$t>0$.
\end{rem}

The following Proposition gives that $(H6)$ holds. 
\begin{prop}
\label{lem:H6}
Assume $\lambda_{0}>0$, $(H2)$-$(H3)$ and  $(H8)$-$(H9)$. Then   $(H6)$
holds with  $\rP^{(-\infty)} = \rP^{(-\infty), \phi_{0}}$. 
\end{prop}
\begin{proof}
Let $0 <  t$ and $F$ be a bounded and $\D_{[-t,0]}$ measurable function. For $h$ large enough, we have:
\begin{align*}
\rE_{x}^{(-h)} \big[  F(Y_{[-t,0]})\big] 
&= \rE_{x}^{(h)} \big[  \rE_{Y_{h/2}}^{(h/2)} \big[
F(\theta_{h/2}({Y})_{[-t,-s]}) \big] \big] \\ 
&= \rE^{\phi_{0}}_{x} \big[ M_{h/2}^{(h), \phi_{0} } \ \rE_{Y_{h/2}}^{(h/2)}
\big[ F(\theta_{h/2}({Y})_{[-t,0]}) \big]  \big] \\ 
&= \rE^{\phi_{0}}_{x} \big[ \rE_{Y_{h/2}}^{(h/2)} \big[
F(\theta_{h/2}({Y})_{[-t,0]}) \big]  \big]+o_h(1) \\ 
&=\rE_{\pi}^{(h/2)} \big[ F(\theta_{h/2}({Y})_{[-t,0]}) \big] +
o_h(1) \\ 
&=\rE_{\pi}^{(-h/2)} \big[ F({Y}_{[-t,0]}) \big] + o_h(1),
\end{align*}
where we used the definition of $\rP^{(-h)}$ and the Markov property for
the  first equality, Lemma \ref{lem:H4} together with $F$ bounded by
$\|F\|_{\infty}$ for the third, and assumption $(H9)$ for  the fourth. 
We continue  the computations as follows:
\begin{align*}
\rE_{x}^{(-h)} \big[  F(Y_{[-t,0]})\big] 
&=\rE_{\pi}^{\phi_{0}} \big[ L^{(-h/2)}   F({Y}_{[-t,0]}) \big] +
o_h(1) \\ 
&=\rE_{\pi}^{\phi_{0}} \big[ L^{(-\infty)}
F({Y}_{[-t,0]}) \big] + o_h(1) \\ 
&=\rE_{\pi}^{(-\infty),\phi_{0}} \big[ F({Y}_{[-t,0]}) \big] + o_h(1) ,
\end{align*}
where we used Lemma \ref{lem:Lh-Linfty} for the second equality.
This gives $(H6)$ with  $\rP^{(-\infty)} = \rP^{(-\infty), \phi_{0}}$. 
\end{proof}
 
\subsection{Proof of Lemmas \ref{boundcritical}, \ref{nonnegativelambdaH1}, \ref{FK},
\ref{bound} and \ref{FK+}}
\label{sec:preuve-annexe}

\begin{proof}[Proof of Lemma \ref{boundcritical}]

From $(H2)$ and $(H8)$, there exist $m, M \in \R$ such that
$$\forall x \in E, 0 < m \leq \alpha \phi_{0}(x) \leq M <\infty .$$
Let $W$ be a $(\frac{M}{\alpha \phi_{0}} \mL,0,M)$ Brownian snake 
and define the time change $\Phi$ for every $w \in \W$ by
$\Phi_{t}(w)=\int_{0}^{t} ds \ \frac{M}{\alpha \phi_{0}}(w(s)).$
As $\partial_{t} \Phi_{t}(w) \geq 1$, we have that $t \to \Phi_{t}(w)$ is strictly increasing. Let $t \to \Phi^{(-1)}_t(w)$ denote its inverse. Then, using Proposition 12 of \cite{DS00}, first step of the proof, we have that the time changed snake $W \circ \Phi^{-1}$, with value 
$$(W \circ \Phi^{-1})_s=(W_s(\Phi^{-1}_t(W_s)), t \in [0,\Phi^{-1}(W_s,H_s)])$$ 
at time $s$, is a $(\mL,0,\alpha \phi_{0})$ Brownian snake. 
Noting the obvious bound on the time change $\Phi_{t}^{-1}(w) \leq t$, we have, according to Theorem 14 of \cite{DS00}:
\[ \PP^{\big( \frac{M}{\alpha \phi_{0}} \ \mL^{\phi_{0}},0,M\big)}_{\frac{\alpha \phi_0(x)}{M}\delta_{x}} \left( \Hm \leq t \right) \geq \PP^{(\mL^{\phi_{0}},0,\alpha \phi_{0})}_{\delta_{x}}\left( \Hm \leq t\right)\]
which implies:
\[ \frac{\alpha \phi_0 (x)}{M} \  \NN^{\big( \frac{M}{\alpha \phi_{0}} \mL^{\phi_{0}},0,M\big)}_{x} \left( \Hm > t \right) \leq \NN^{(\mL^{\phi_{0}},0,\alpha \phi_{0})}_{x}\left( \Hm > t\right)\]
from the exponential formula for Poisson point measures.
% % \[\N^{(\frac{(M)}{\alpha \phi_{0}} \mL,0,(M))}_{x} (H_{max} > t) \leq \N^{(\mL,0,\alpha \phi_{0})}_{x}(H_{max} > t)\]
% from the exponential formula \reff{eq:masterformula}.
Now, the left hand side of this inequality can be computed explicitly:
\[
\NN^{\big(\frac{M}{\alpha \phi_{0}} \mL^{\phi_{0}},0,M\big)}_{x} \left(\Hm >t \right)=\N^{\big(\frac{M}{\alpha \phi_{0}} \mL^{\phi_{0}},0,M\big)}_{x} \left( \Hm >t\right)=\frac{1}{Mt}
\]
and the right hand side of this inequality is $v^{\phi_{0}}_{t}(x)$ from \reff{eq:vphi}.
We thus have proved that:
\[  \frac{\alpha \phi_0 (x)}{M^2 t} \leq v^{\phi_{0}}_{t}(x), \]
and this yields the first part of the inequality of Lemma \ref{boundcritical}. The second part is obtained in the same way using the coupling with the $\big(\frac{m}{\alpha \phi_{0}} \mL^{\phi_{0}},0,m\big)$ Brownian snake. 
\end{proof}
\begin{proof}[Proof of lemma \ref{nonnegativelambdaH1}]
Assumption $(H2)$ and $(H8)$ allow us to apply Lemma \ref{boundcritical} for the case $\lambda_{0}=0$, which yields that $v^{\phi_0}_{\infty}=0$, and then $v_{\infty}=0$ thanks to \reff{eq:normalisation-vf}. This in turn implies that $(H1)$ holds in the case $\lambda_{0}=0$ according to Lemma \ref{lemmaextinction}.
For $\lambda_{0}>0$, we may use item 5 of Proposition 13 of \cite{DS00} (which itself relies on a Girsanov theorem) with $\P^{(\mL,0,\alpha \phi_{0})}$ in the r\^ole of $\P^{c}$ and $\P^{(\mL^{\phi{0}},\lambda_{0},\alpha \phi_{0})}$ in the r\^ole of $\P^{b,c}$ to conclude that the extinction property $(H1)$ holds under $\P^{(\mL^{\phi{0}},\lambda_{0},\alpha \phi_{0})}$. 
\end{proof}

\begin{proof}[Proof of Lemma \ref{FK}]
Let $ \varepsilon>0$.
The function $v^{\phi_{0}}$ is known to solve the following mild form of the Laplace equation, see equation \reff{mildequation}:
\begin{equation*}
\quad v_{t+s}^{\phi_{0}}(x) + \rE^{\phi_{0}}_{x}\bigg[\int_{0}^{t} dr \     
\big(\lambda_{0} \: v^{\phi_{0}}_{t+s-r}(Y_r) + \alpha(Y_r) \phi_{0}(Y_r)
(v^{\phi_{0}}_{t+s-r}(Y_r))^2 \big)
\bigg]=\rE_{x}^{\phi_{0}}\big[v^{\phi_{0}}_s(Y_{t})\big].  
\end{equation*}
By differentiating  with respect  to $s$ and  taking $t=t-s$,  we deduce
from   dominated  convergence   and  the   bounds  \reff{eq:majovtilde},
\reff{eq:derivvtilde}  and \reff{eq:ineqSigma} on  $v^{\phi_{0}}=v/\phi_{0}$ and
its  time derivative (valid under the assumptions $(H1)$-$(H3)$) the following  mild  form on  the time  derivative
$\partial_{t} v^{\phi_{0}}$:
\begin{equation*}
\quad \partial_{t} v_t^{\phi_{0}}(x) + \rE^{\phi_{0}}_{x}\bigg[\int_{0}^{t-s} dr \ 
\big(\lambda_{0} \:  + 2 \alpha(Y_r) \phi_{0}(Y_r) v^{\phi_{0}}_{t-r}(Y_r) \big) \partial_{t} v^{\phi_{0}}_{t-r}(Y_r)
\bigg]=\rE^{\phi_{0}}_{x}\big[\partial_{t} v^{\phi_{0}}_s(Y_{t-s})\big].
\end{equation*}
From the Markov property, for fixed $t>0$, the two following processes:
\begin{equation*}
\bigg(v^{\phi_{0}}_{t-s}(Y_{s})- \int_{0}^{s} dr \     
\big(\lambda_{0}  + \alpha(Y_r) \phi_{0}(Y_r) v^{\phi_{0}}_{t-r}(Y_r)\big) \: v^{\phi_{0}}_{t-r}(Y_r) , 0 \leq s < t \bigg)
\end{equation*}
and
\begin{equation*}
\bigg(\partial_{t} v^{\phi_{0}}_{t-s}(Y_{s})- \int_{0}^{s} dr \ 
\big(\lambda_{0} \:  + 2 \alpha(Y_r) \phi_{0}(Y_r) v^{\phi_{0}}_{t-r}(Y_r)
\big) \partial_{t} v^{\phi_{0}}_{t-r}(Y_r) , 0 \leq s < t \bigg)
\end{equation*}
are  $\D_{s}$-martingale under $\rP_{\pi}^{\phi_{0}}$.
A Feynman-Kac manipulation, as done in the proof of Lemma \ref{AClemma},
enables us to conclude that for fixed $t>0$:
\begin{equation*}
\bigg( v^{\phi_{0}}_{t-s}(Y_{s}) \expp{- \int_{0}^{s} dr \ \big( 
\lambda_{0}  + \alpha(Y_r) \phi_{0}(Y_r) v^{\phi_{0}}_{t-r}(Y_r) \big)}, 0 \leq
s <t \bigg) 
\end{equation*}
and
\begin{equation*}
\bigg( \partial_{t} v^{\phi_{0}}_{t-s}(Y_{s}) \expp{-\int_{0}^{s} dr \
  \big(\lambda_{0} \:  + 2 \alpha(Y_r) \phi_{0}(Y_r) v^{\phi_{0}}_{t-r}(Y_r)
  \big)}, 0 \leq s < t \bigg) 
\end{equation*}
are $\D_{s}$-martingale under $\rP_{\pi}^{\phi_{0}}$. Taking expectations
at time $s=0$ and $s=h$ with $t=h+\varepsilon$, we get the representations
formulae stated in the Lemma: 
\begin{align*}
v_{h+\varepsilon}^{\phi_{0}}(x) 
&=  \expp{-\lambda_{0} h}  \rE^{\phi_{0}}_{x} \bigg[ \expp{- \int_{0}^{h} ds
  \ \alpha(Y_s) \: \phi_{0}(Y_s) \: v^{\phi_{0}}_{h+\varepsilon-s}(Y_s)}
v_{\varepsilon}^{\phi_{0}}(Y_h)\bigg],\\ 
\partial_{h} v_{h+\varepsilon}^{\phi_{0}}(x) 
&= \expp{-\lambda_{0} h} \rE^{\phi_{0}}_{x} \bigg[ \expp{- 2\int_{0}^{h} ds \
  \alpha(Y_s) \: \phi_{0}(Y_s) \: v^{\phi_{0}}_{h+\varepsilon-s}(Y_s)} \partial_{h}
v_{\varepsilon}^{\phi_{0}}(Y_h)\bigg]. 
\end{align*}

\end{proof}
\begin{proof}[Proof of Lemma \ref{bound}]
  Since
  $v^{\phi_{0}}_{\varepsilon}=v_{\varepsilon}/\phi_{0}=\tilde{v}_{\varepsilon}/(\alpha
  \phi_{0})$, we  can conclude from \reff{eq:majovtilde},  $(H2)$ and $(H8)$
  that $v^{\phi_{0}}_{\varepsilon}$ is bounded  from above and from below by
  positive     constants.     Similarly,     we     also    get     from
  \reff{eq:derivvtilde},   \reff{eq:defSigma}  and  \reff{eq:ineqSigma}
  that  $|\partial_{h} \tilde{v}_{\varepsilon}|$  is bounded  from above
  and from below by two positive constants. Thus, we have  the existence
  of four positive constants,  $D_1$,  $D_2$, $D_3$ and 
  $D_4$,  such that, for all $x
  \in E$:
\begin{align}
\label{eq:initialbound1}
D_1  & \leq v_{\varepsilon}^{\phi_{0}}(x) \leq D_2,\\
\label{eq:initialbound2}
D_3  & \leq |\partial_{t} v_{\varepsilon}^{\phi_{0}}(x)| \leq D_4.
\end{align}
From equations \reff{eq:FKv}, \reff{eq:initialbound1} and the positivity
of $v^{\phi_{0}}$, we deduce that: 
\begin{align}
\label{eq:boundv}
v_{h+\varepsilon}^{\phi_{0}}(x) \leq  D_2 \expp{-\lambda_{0} h}. 
\end{align}
Putting back \reff{eq:boundv} into \reff{eq:FKv}, we have the converse
inequality $D_5 \expp{-\lambda_{0} h}  \leq v_{h+\varepsilon}^{\phi_{0}}(x)
$
with $D_5 = D_1 \exp{\{- D_2 \norm{\alpha}_\infty \norm{\phi_{0}}_\infty /
  \lambda_{0}\}}>0$. This gives \reff{eq:encadr-v}. 

Similar  arguments   using  \reff{eq:FKdv}  and  \reff{eq:initialbound2}
instead    of   \reff{eq:FKv}    and    \reff{eq:initialbound1},   gives
\reff{eq:encadr-dv}.
\end{proof}

\begin{proof}[Proof of Lemma \ref{FK+}]
Using  the Feynman-Kac representation of $\partial_{h}
v_{h+\varepsilon}^{\phi_{0}}$ from  \reff{eq:FKv} as well as the Markov property, we have:
\begin{align*}
\partial_{h} v_{h+\varepsilon}^{\phi_{0}}(x) \expp{\lambda_{0} h} 
&=
\rE^{\phi_{0}}_{x} \bigg[ \expp{- 2\int_{0}^{h} ds \  \alpha(Y_s) \: \phi_{0}(Y_s) \:
  v^{\phi_{0}}_{h+\varepsilon-s}(Y_s)} \partial_{h}
v_{\varepsilon}^{\phi_{0}}(Y_h)\bigg] \\ 
&= \rE^{\phi_{0}}_{x} \bigg[ \expp{- 2 \int_{0}^{\sqrt{h}} \!\!ds \:
  \alpha \: 
  \phi_{0} \: v^{\phi_{0}}_{h+\varepsilon-s}(Y_s)} \rE^{\phi_{0}}_{Y_{\sqrt{h}}} \big[
\expp{- 2 \int_{0}^{h-\sqrt{h}} \!\! ds \: \alpha \: \phi_{0} \:
  v^{\phi_{0}}_{h-\sqrt{h}+\varepsilon-s}( Y_s)} \partial_{h}
v_{\varepsilon}^{\phi_{0}}({Y}_{h-\sqrt{h}}) \big]\bigg].
\end{align*}
Notice that
\begin{equation}
   \label{eq:majo-int-h}
\bigg| \int_{0}^{\sqrt{h}} ds \  \alpha \:
  \phi_{0} \: v^{\phi_{0}}_{h+\varepsilon-s}(Y_s) \bigg|
\leq  \norm{\alpha}_\infty \norm{\phi_{0}}_\infty  \sqrt{h}
\norm{v^{\phi_{0}}_{h+\varepsilon- \sqrt{h}}}_\infty =o_h(1),
\end{equation}
according   to   Lemma    \ref{bound}   if   $\lambda_0>0$   and   Lemma
\ref{boundcritical} if $\lambda_0=0$.
We get:
\begin{align*}
\partial_{h} v_{h+\varepsilon}^{\phi_{0}}(x) \expp{\lambda_{0} h} 
&= \rE^{\phi_{0}}_{x} \bigg[ \rE^{\phi_{0}}_{Y_{\sqrt{h}}} \big[  \expp{-
  2\int_{0}^{h-\sqrt{h}} ds \  \alpha(Y_s) \: \phi_{0}(Y_s) \:
  v^{\phi_{0}}_{h-\sqrt{h}+\varepsilon-s}( Y_s)} \partial_{h}
v_{\varepsilon}^{\phi_{0}}({Y}_{h-\sqrt{h}}) \big]\bigg]
\big(1+o_h(1)\big) \\ 
&= \rE^{\phi_{0}}_{\pi} \bigg[   \expp{- 2\int_{0}^{h-\sqrt{h}} ds \ 
  \alpha(Y_s) \: \phi_{0}(Y_s) \: v^{\phi_{0}}_{h-\sqrt{h}+\varepsilon-s}(
  Y_s)} \partial_{h} v_{\varepsilon}^{\phi_{0}}({Y}_{h-\sqrt{h}})
\bigg] \big(1+o_h(1)\big) \\ 
&= \rE^{\phi_{0}}_{\pi} \bigg[   \expp{- 2 \int_{-(h-\sqrt{h})}^{0} ds \
  \alpha(Y_s) \: \phi_{0}(Y_s) \: v^{\phi_{0}}_{\varepsilon-s}( Y_s)} \partial_{h}
v_{\varepsilon}^{\phi_{0}}({Y}_{0}) \bigg] \big(1+o_h(1)\big)\\ 
&= \rE^{\phi_{0}}_{\pi} \bigg[   \expp{- 2\int_{-h}^{0} ds \ \alpha(Y_s) \: \phi_{0}(Y_s)
  \: v^{\phi_{0}}_{\varepsilon-s}( Y_s)} \partial_{h}
v_{\varepsilon}^{\phi_{0}}({Y}_{0}) \bigg] \big(1+o_h(1)\big),
\end{align*}
where we used \reff{eq:majo-int-h} for the first equality, $(H9)$ for
the second, stationarity of $Y$ under $\rP^{\phi_{0}}_{\pi}$ for the third
and \reff{eq:majo-int-h} again for the last. 
This gives \reff{eq:dev-dv1-enh}.

Moreover, if $\lambda_{0}>0$, we get that:
\[
\rE^{\phi_{0}}_{\pi} \bigg[   \expp{- 2 \int_{-\infty}^{0} ds \  \alpha \:
  \phi_{0} \: v^{\phi_{0}}_{\varepsilon-s}( Y_s)} \partial_{h}
v_{\varepsilon}^{\phi_{0}}({Y}_{0}) \bigg] 
\]
is finite and that:
\[
\lim_{h'\rightarrow+\infty }
\rE^{\phi_{0}}_{\pi} \bigg[   \expp{- 2 \int_{-h'}^{0} ds \  \alpha \:
  \phi_{0} \: v^{\phi_{0}}_{\varepsilon-s}( Y_s)} \partial_{h}
v_{\varepsilon}^{\phi_{0}}({Y}_{0}) \bigg] 
=\rE^{\phi_{0}}_{\pi} \bigg[   \expp{- 2 \int_{-\infty}^{0} ds \  \alpha \:
  \phi_{0} \: v^{\phi_{0}}_{\varepsilon-s}( Y_s)} \partial_{h}
v_{\varepsilon}^{\phi_{0}}({Y}_{0}) \bigg] .
\]
Therefore, we deduce \reff{eq:dev-dv2-enh} from \reff{eq:dev-dv1-enh}.
\end{proof}

\subsection{About the Bismut spine.}

Choosing uniformly an individual at random at height $t$ under $\NN_x$ and letting $t \to \infty$, we will see that the law of the ancestral lineage should converge in some sense to the law of the oldest ancestral lineage which itself converges to $\rP_{x}^{(\infty)}$ defined in \reff{defPinfty}, according to Lemma \ref{prop:H4+Pf=Pinf}.

We have defined in \reff{eq:bismutspine} the following family of probability measure indexed by $t \geq 0$:

\begin{align*}
\frac{d \rP_{x \ | \D_{t}}^{(B,t)}}{d \rP_{x \   | \D_{t}}} = \frac{\expp{-\int_{0}^{t} ds \; \beta(Y_s)}}{\rE_{x} \left[\expp{-\int_{0}^{t} ds \; \beta(Y_s)} \right]} \cdot
\end{align*}

\begin{lem}
\label{lem:convbismutspine}
Assume $(H8)$-$(H9)$. We have, for every $0 \leq t_0 \leq t$:
\[\frac{d \rP_{x \ | \D_{t_0}}^{(B,t)}}{d \rP_{x \   | \D_{t_0}}} \; \xrightarrow[t \to +\infty]{}  \; \frac{d \rP_{x \ | \D_{t_0}}^{(\infty)}}{d \rP_{x \   | \D_{t_0}}} \quad \text{$\rP_{x}$-a.s. and in $\rL^{1}(\rP_{x})$}. \]
\end{lem}

Note that there is no restriction on the sign of $\lambda_0$ for this Lemma to hold. 

\begin{rem}
This result correspond to the so called globular phase in the random polymers literature (see \cite{CK10}, Theorem 8.3).
\end{rem}

\begin{proof}
We have:
\begin{align*}
\frac{d \rP_{x \ | \D_{t_0}}^{(B,t)}}{d \rP_{x \   | \D_{t_0}}} &= \expp{-\int_{0}^{t_0} ds \; \beta(Y_s)} \frac{\rE_{Y_{t_0}}\left[\expp{-\int_{0}^{t-t_0} ds \; \beta(Y_s)}\right]}{\rE_{x}\left[\expp{-\int_{0}^{t} ds \; \beta(Y_s)}\right]} \\
&= \expp{-\int_{0}^{t_0} ds \;  (\beta(Y_s) - \lambda_0)} \frac{\rE_{Y_{t_0}}\left[\expp{-\int_{0}^{t-t_0} ds \; (\beta(Y_s)-\lambda_0)}\right]}{\rE_{x}\left[\expp{-\int_{0}^{t} ds \; (\beta(Y_s)-\lambda_0)}\right]}\\
&= \expp{-\int_{0}^{t_0} ds \;  (\beta(Y_s) - \lambda_0)} \frac{\phi_0(Y_{t_0})}{\phi_0(Y_0)} \frac{\rE_{Y_{t_0}}\left[\expp{-\int_{0}^{t-t_0} ds \; (\beta(Y_s)-\lambda_0)} \frac{\phi_0(Y_{t-t_0})}{\phi_0(Y_0)} \frac{1}{\phi_0(Y_{t-t_0})} \right]}{\rE_{x}\left[\expp{-\int_{0}^{t} ds \; (\beta(Y_s)-\lambda_0)} \frac{\phi_0(Y_{t})}{\phi_0(x)} \frac{1}{\phi_0(Y_t)}\right]} \\
&=  \frac{d \rP^{\phi_0}_{x \ | \D_{t_0}}}{d \rP_{x \   | \D_{t_0}}} \frac{\rE^{\phi_0}_{Y_{t_0}} \left[1/\phi_0(Y_{t-t_0})\right]}{\rE^{\phi_0}_{x}\left[ 1/\phi_0(Y_t) \right]} \\
& \underset{t \rightarrow \infty }{\longrightarrow} \frac{d \rP_{x \ | 
\D_{t_0}}^{\phi_0}}{d \rP_{x \   | \D_{t_0}}} \frac{\pi(\frac{1}{\phi_0})}{\pi(\frac{1}{\phi_0})} 
=  \frac{d \rP_{x \ | \D_{t_0}}^{\phi_0}}{d \rP_{x \   | \D_{t_0}}},
\end{align*}
where we use the Markov property at the first equality, we force the apparition of $\lambda_0$ at the second equality and we force the apparition of $\phi_0$ at the third equality in order to obtain the Radon Nikodym derivative of $\rP^{\phi_0}_{x}$ with respect to $\rP_{x}$: this observation gives the fourth equality. The ergodic assumption $(H9)$ ensures the $\rP_{x}$-a.s. convergence to $1$ of the fraction in the fourth equality as $t$ goes to $\infty$. Since $$\left((t,y) \to \rE^{\phi_0}_{y}\left[1/\phi_0(Y_{t-t_0}) \right]/\rE^{\phi_0}_{x}\left[1/\phi_0(Y_t)\right]\right)$$ is bounded according to $(H8)$, we conclude that the convergence also holds in $\rL^{1}(\rP_{x})$. Then use Lemma \ref{prop:H4+Pf=Pinf} to get that $\rP^{\phi_0}_{x}=\rP^{(\infty)}_{x}$.
\end{proof}

\section{Two examples}
\label{sec:twoexamples}

In this section, we specialize the results of the previous sections to the case of the multitype Feller process and of the superdiffusion.

\subsection{The multitype Feller diffusion}
\label{sec:Ex-multi}

The  multitype Feller diffusion  is the  superprocess with  finite state
space: $E=\{1,  \dots, K\}$ for $K$  integer. In this  case, the spatial
motion  is   a  pure  jump   Markov  process,  which  will   be  assumed
irreducible. Its  generator $\mL$ is  a square matrix  $(q_{ij})_{1 \leq
  i,j \leq K}$ of size $K$  with lines summing up to $0$, where $q_{ij}$
gives the transition rate from $i$  to $j$ for $i \neq j$. The functions
$\beta$    and    $\alpha$     defining    the    branching    mechanism
\reff{eq:branchingmechanism} are vectors of  size $K$: this implies that
$(H2)$  and  $(H3)$ automatically  hold.   For  more  details about  the
construction   of  finite   state  space   superprocess,  we   refer  to
\cite{DY94}, example 2,  p. 10, and to \cite{CH08}  for investigation of
the Q-process.

The generalized eigenvalue $\lambda_{0}$ is defined by:
\begin{equation}
 \label{defeigenvalue1}
\lambda_0 = \sup{\{ \ell \in \R, \exists u > 0 \mbox{ such that } (\mbox{Diag}(\beta)-\mL) u = \ell u \}},
\end{equation}
where Diag$(\beta)$ is the diagonal $K \times K$ matrix with diagonal coefficients derived from the vector $\beta$. We stress that the generalized eigenvalue is also the Perron Frobenius eigenvalue, \textit{i.e.} the eigenvalue with the maximum real part, which is real by Perron Frobenius theorem, see \cite{SE06}, Exercise 2.11. Moreover, the associated eigenspace is one-dimensional. We will denote by $\phi_{0}$ and $\tilde{\phi}_{0}$ its generating left, resp. right, eigenvectors, normalized so that $\sum_{i=1}^{K} \phi_{0}(i) \tilde{\phi}_{0}(i)=1$, and the coordinates of $\phi_{0}$ and $\tilde{\phi}_{0}$ are positive.

We first check that the two assumptions we made in Section 6 are satisfied.
\begin{lem}
\label{rem:H9-Feller}
Assumptions $(H8)$ and $(H9)$ hold with $\pi= \phi_{0} \: \tilde{\phi}_{0}$.
\end{lem}

\begin{proof}
Assumption $(H8)$ is obvious in the finite state space setting. Assumption $(H9)$ is a classical statement about irreducible finite state space Markov Chains.
\end{proof}

\begin{lem}
\label{lem1234}
Assume $\lambda_{0} \geq 0$. Then $(H1)$, $(H4)$ and $(H5)_{\nu}$ hold.
\end{lem}
\begin{proof}
Assumption $(H2)$ and $(H8)$ hold according to Lemma \ref{rem:H9-Feller}. Together with $\lambda_{0} \geq 0$, this allows us to apply Lemma \ref{nonnegativelambdaH1} to obtain $(H1)$.
Then use Proposition \ref{prop:H4+Pf=Pinf} to get $(H4)$ and Corollary \ref{cor:H5} to get $(H5)_{\nu}$.
\end{proof}

% observe that the global extinction property $(H1)$ is equivalent to the non negativity of $\lambda_{0}$.

\begin{lem}
\label{H6H7multi}
Assume $\lambda_{0} > 0$. Then $(H6)$ and $(H7)$ holds.
\end{lem}

\begin{proof}
We apply Proposition \ref{lem:H6} to prove $(H6)$ and Lemma \ref{lem:H7} to prove $(H7)$.
\end{proof}

Recall that $\rP_{x}^{(h)}$ and $\rP_{x}^{(\infty)}$ were defined in  \reff{defPh2} and \reff{defPinfty} respectively.
\begin{lem}
\label{lem:spinemulti}
We have:
\begin{itemize}
 \item[(i)] $\rP_{x}^{(h)}$ is a continuous time inhomogeneous Markov chain on $[0,h)$ issued from $x$ with transition rates from $i$ to $j$, $i \neq j$, equal to $\frac{\partial_{h} v_{h-t}(j)}{\partial_{h} v_{h-t}(i)} q_{ij}$ at time $t$, $0 \leq t< h$.
 \item[(ii)] $\rP_{x}^{(\infty)}$ is a continuous time homogeneous Markov chain on $[0,\infty)$ issued from $x$ with transition rates from $i$ to $j$, $i \neq j$, equal to $\frac{\phi_{0}(j)}{\phi_{0}(i)} q_{ij}$.
\end{itemize}
\end{lem}

\begin{proof}
The first item is a consequence of a small adaptation of Lemma \ref{generatorLh} for time dependent function.
Namely, let $g_t(x)$ be a time dependent function. Consider the law of process $(t,Y_t)$ and consider the probability measure $\rP^{g}$ defined by \reff{ACequation} with $g(t,Y_t) = g_{t}(Y_t)$. Denoting by $\mL^{g}_t$ the generator of (the inhomogeneous Markov process) $Y_t$ under $\rP^{g}$, we have that:
\begin{align}
\label{Lgt}
\forall u \in \D_{g}(\mL), \quad  \mL^{g}_t(u)= \frac{\mL(g_t u)-\mL(g_t)u}{g_t}\cdot 
\end{align}
% \label{defLh}
Recall that for all vector $u$, $ \mL(u)(i)= \sum_{j \neq i} q_{ij} \big(u(j)-u(i)\big)$. Then apply \reff{Lgt} to the time dependent function $g_t(x)=\partial_{t}v_{h-t}(x)$, and note that $\rP^{g}=\rP^{(h)}$ thanks to \reff{defPh2}.
For the second item, observe that Proposition \ref{prop:H4+Pf=Pinf} identifies $\rP_{x}^{(\infty)}$ with $\rP^{\phi_{0}}$. Use then Lemma \ref{generatorLh} to conclude.
\end{proof}

% convergence of $\P^{(\geq h)}$ to the Q-process (Lemma \ref{lemmaconv})

William's decomposition under $\N_{x}^{(h)}$ (Propositions \ref{corwilliams}) together with the convergence of this decomposition (Theorem \ref{theo-convN}) then hold under the assumption $\lambda_{0} \geq 0$.
Convergence of the distribution of the superprocess near its extinction time under $\N_{x}^{(h)}$ (Proposition \ref{convtop}) holds under the stronger assumption $\lambda_{0}>0$. We were unable to derive an easier formula for $\rP^{(-\infty)}$ in this context.

Remark that Lemma \ref{lemmaconv}, Definition \ref{defi:Q} and Corollary \ref{cor:cvNh} give a precise meaning to the ``interactive immigration'' suggested by the authors in Remark 2.8. of \cite{CH08}.

\subsection{The superdiffusion}
\label{sec:Ex-diffu}

% , see \cite{EP99} for a detailed analysis.
% The following elements of definition are taken from 
The superprocess associated to a diffusion is called superdiffusion. We first define the diffusion and the relevant quantities associated to it, and take for that the general setup from \cite{PI95}.
Here $E$ is an arbitrary domain of $\R^{K}$ for $K$ integer. Let $a_{ij}$ and ${b_i}$ be in $\C^{1,\mu}(E)$, the usual H\"older space of order $\mu \in [0,1)$, which consists of functions whose first order derivatives are locally H\"older continuous with exponent $\mu$, for each $i,j$ in $\{1, \ldots, K\}$. Moreover, assume that the functions $a_{i,j}$ are such that the matrix $(a_{ij})_{(i,j) \in \{1...K\}^2}$ is positive definite. Define now the generator $\mL$ of the diffusion to be the elliptic operator:
\[ \mL(u) = \sum_{i=1}^{K} b_i \; \partial_{x_i}u  + \frac{1}{2} \sum_{i,j=1}^{K} a_{ij} \; \partial_{x_i,x_j}u.\]
The generalized eigenvalue $\lambda_{0}$ of the operator $\beta - \mL$ is defined by:
\begin{align}
\label{defeigenvalue2}
\lambda_{0} = \sup{\{ \ell \in \R, \exists u\in \D(\mL), u > 0 \mbox{ such that } (\beta - \mL) u = \ell \  u \}} \cdot
\end{align}
Denoting $\rE$ the expectation operator associated to the process with generator $\mL$, we recall an equivalent probabilistic definition of the generalized eigenvalue $\lambda_{0}$:
\begin{align*}
 \lambda_{0}= - \sup_{A \subset \R^K} \lim_{t \to \infty}  \frac{1}{t} \log \rE_{x}\big[ \expp{-\int_{0}^{t}  ds \; \beta(Y_s)} \ind_{\{\tau_{A^c} >t\}} \big], 
\end{align*}
for any $x \in \R^K$, where $\tau_{A^c} = \inf{ \{t > 0: Y(t) \notin A\}}$ and the supremum runs over the compactly embedded subsets $A$ of $\R^K$.
We assume that the operator $(\beta-\lambda_{0})- \mL$ is critical in the sense that the space of positive harmonic functions for $(\beta-\lambda_{0}) - \mL$ is one dimensional, generated by $\phi_{0}$. In that case, the space of positive harmonic functions of the adjoint of $(\beta-\lambda_{0})- \mL$ is also one dimensional, and we denote by $\tilde{\phi_{0}}$ a generator of this space.
We further assume that the operator $(\beta-\lambda_{0})- \mL$ is \textbf{product-critical}, \textit{i.e.} $\int_{E} dx \: \phi_{0}(x) \: \tilde{\phi_{0}}(x) < \infty$, in which case we can normalize the eigenvectors in such a way that $\int_{E} dx \: \phi_{0}(x) \: \tilde{\phi_{0}}(x) =1$. This assumption (already appearing in \cite{EK10}) is a rather strong one and implies in particular that $\rP^{\phi_0}$ is the law recurrent Markov process, see Lemma \ref{rem:H9-diffu} below.

% , which will be assumed bounded from below and from above by two positive constants (this is our assumption $(H8)$)
Concerning the branching mechanism, we will assume, in addition to the conditions stated in section \ref{section0}, that $\alpha \in \C^4(E)$.

% We first check that assumption $(H9)$ holds.
\begin{lem}
\label{rem:H9-diffu}
Assume $(H8)$. Assumption $(H9)$ holds with $\pi(dx)=\phi_{0}(x) \: \tilde{\phi}_{0}(x) \:dx$.
\end{lem}

\begin{proof}
  We  assume that  $(\beta-\lambda_{0}) -  \mL$ is  a  critical operator
  which is product critical.  Note that $-\mL^{\phi_{0}}$ is the (usual)
  $h$-transform of  the operator $(\beta  -\lambda_{0})- \mL$ with  $h =
  \phi_{0}$,   where  the  $h$-transform   of  $\mL(\cdot)$   is  $\mL(h
  \cdot)/h$.     Then   Remark   5    of   \cite{EK10}    implies   that
  $-\mL^{\phi_{0}}$ is  again a critical operator which  is also product
  critical with corresponding $\phi_0$ and $\tilde \phi_0$ given by 1 and
  $\phi_{0} \; \tilde{\phi}_{0}$. Then Theorem 9.9 p.192 of \cite{PI95},
  see (9.14), states that $(H9)$ holds.
\end{proof}

% , see Theorem 3.3 p. 148 of \cite{PI95} But

Note that the non negativity of the generalized eigenvalue of the operator $(\beta - \mL)$ now characterizes in general the local extinction property (the superprocess $X$ suffers local extinction if its restrictions to compact domains of $E$ suffers global extinction); see \cite{EK04} for more details on this topic. However, under the boundedness assumption we just made on $\alpha$ and $\phi_{0}$, the extinction property $(H1)$ holds, as will be proved (among other things) in the following Lemma.

\begin{lem}
\label{lem123}
Assume $\lambda_{0} \geq 0$ and $(H8)$. Then $(H1)$-$(H4)$ and $(H5)_{\nu}$ hold. If moreover $\lambda_{0}>0$, then $(H6)$ and $(H7)$ holds.
\end{lem}
\begin{proof}
The assumption $\alpha \in \C^{4}(E)$ ensures that $(H2)$ and $(H3)$ hold. 
Then the end of the proof is similar to the end of the proof of Lemma \ref{lem1234} and the proof of Lemma \ref{H6H7multi}. 
\end{proof}

Recall that $\rP_{x}^{(h)}$ and $\rP_{x}^{(\infty)}$ were defined in  \reff{defPh2} and \reff{defPinfty}.

\begin{lem}
\label{lem:spinediff}
We have:
\begin{itemize}
 \item $\rP^{(h)}_{x}$ is an inhomogeneous diffusion on $[0,h)$ issued
   from $x$ with generator at time $t\in [0,h)$:  $\displaystyle
   \left(\mL+a \frac{\nabla \partial_{h} v_{h-t}}{\partial_{h} v_{h-t}}
   \nabla .\right)$. 
 \item $\rP^{(\infty)}_{x}$ is an homogeneous diffusion on $[0,\infty)$ issued from $x$ with generator 
$
   \left(\mL+a \frac{\nabla \phi_{0}}{\phi_{0}} \nabla . \right)$.
\end{itemize}
\end{lem}
\begin{proof}
The proof is similar to the proof of Lemma \ref{lem:spinemulti}.
\end{proof}

% Informally, a drift pointing to the high values of $\partial_{h} v_{h-t}$ (resp. $\phi_{0}$) is added to the initial process with generator $\mL$ in the construction of $\rP_{x}^{(h)}$ (resp. $\rP_{x}^{(\infty)}$). 

William's decomposition under $\N_{x}^{(h)}$ (Propositions \ref{corwilliams}) together with the convergence of this decomposition (Theorem \ref{theo-convN}) then hold under the assumption $\lambda_{0} \geq 0$ and $(H8)$.
Convergence of the distribution of the superprocess near its extinction time under $\N_{x}^{(h)}$ (Proposition \ref{convtop}) holds under the stronger assumption $\lambda_{0}>0$. 

\begin{rem}
\label{backbone-spine}
Engl\"ander and Pinsky offer in \cite{EP99} a decomposition of
supercritical non-homoge\-neous superdiffusion using immigration on the
backbone formed by the prolific individuals (as denominated further in
Bertoin, Fontbona and Martinez \cite{BFM08}). It is interesting to note
that the generator of the backbone is $\mL^{w}$ where $w$ formally
satisfies the evolution equation $\mL w = \psi(w)$, whereas the generator
of the spine of the Q process investigated in Theorem \ref{theo-convN}
is $\mL^{\phi_{0}}$ where $\phi_{0}$ formally satisfies $\mL \phi_0 =
\beta \phi_0$. 
In particular, we notice that the generator of the backbone $\mL^{w}$
depends on both $\beta$ and $\alpha$ and that the generator of our spine
$\mL^{\phi_0}$ does not depend on $\alpha$.  
\end{rem}

%\bibliographystyle{abbrv}
%\bibliography{doku}

\end{document}